\documentclass[12pt,letterpaper]{article}

\usepackage{renstyles}

\newlength{\tempdima}

\newcommand{\rowname}[1]%
{\rotatebox{90}{\makebox[\tempdima][c]{\rotatebox{-90}{#1}}}}

\newcommand{\rownamerev}[1]%
{\rotatebox{90}{\makebox[\tempdima][c]{{#1}}}}

\renewcommand{\RED}[1]{#1}
\newcommand{\edited}[1]{{\color{red} {#1}}}
\renewcommand{\edited}[1]{{#1}}

\title{A three-stage method for reconstructing multiple coefficients in coupled photoacoustic and diffuse optical imaging}

\author{
	Yinxi Pan\thanks{Department of Applied Physics and Applied Mathematics, Columbia University, New York, NY 10027; yp2586@columbia.edu}
    \and	Kui Ren\thanks{Department of Applied Physics and Applied Mathematics, Columbia University, New York, NY 10027; kr2002@columbia.edu}
    \and
        Shanyin Tong\thanks{Department of Applied Physics and Applied Mathematics, Columbia University, New York, NY 10027; st3503@columbia.edu (corresponding author)}
}

\begin{document}

\maketitle

\begin{abstract}

This paper studies inverse problems in quantitative photoacoustic tomography with additional optical current data supplemented from diffuse optical tomography. We propose a three-stage image reconstruction method for the simultaneous recovery of the absorption, diffusion, and Gr\"uneisen coefficients. We demonstrate, through numerical simulations, that: (i) when the Gr\"uneisen coefficient is known, the addition of the optical measurements allows a more accurate reconstruction of the scattering and absorption coefficients; and (ii) when the Gr\"uneisen coefficient is not known, the addition of optical current measurements allows us to reconstruct uniquely the Gr\"uneisen, the scattering and absorption coefficients. Numerical simulations based on synthetic data are presented to demonstrate the effectiveness of the proposed idea.

\end{abstract}



\section{Introduction}
\label{SEC:Intro}

Diffuse optical tomography (DOT) is a medical imaging method that utilizes near-infrared (NIR) light to probe scattering and absorption properties of optically heterogeneous mediums such as biological tissues. In DOT, we send NIR photons into the medium to be probed. The photons travel inside the medium following a given physical law of absorption and scattering. We measure the current of the photons that exit the medium. From such measured data, we are interested in reconstructing the scattering and absorption properties of the medium. We refer interested readers to ~\cite{Arridge-IP99, ArSc-IP09,Bal-IP09,Ren-CiCP10,AsBeCaCaNa-DCDS24,Zimmermann-IP23,BeCaQu-arXiv24,KaMoPuTa-BOE24,YiYaWaWaGuCaZhHe-BOE24,SiJaPr-EPJP24,RaFa-PRR23} and references therein for overviews on the physics and mathematics of this imaging modality as well as some random samples of recent developments in the field.

Besides the fact that the NIR photons used in DOT are less harmful to biological tissues than the x-rays used in traditional x-ray tomography, the main advantage of DOT over x-ray tomography is the ability to image optical contrast of biological tissues (in terms of optical properties of healthy and unhealthy tissues) that x-ray tomography can not offer. The main disadvantage of DOT is its poor imaging resolution. Mathematically speaking, high-frequency information in the scattering and absorption properties of tissues are damped heavily, making it almost impossible to recover that information from measured DOT data. 

Photoacoustic tomography (PAT) is a hybrid imaging method that attempts to overcome the low-resolution disadvantage DOT. PAT couples ultrasound imaging and diffuse optical tomography via the photoacoustic effect, enabling high-resolution imaging of optical contrasts of heterogeneous mediums. A typical way to induce the photoacoustic effect is to send a short pulse of near-infrared (NIR) light into an optically heterogeneous medium, such as a piece of biological tissue. In the light propagation process, a portion of the photons are absorbed by the medium. The absorbed energy causes the medium to heat up slightly, and then cool down after the rest of the photons exit. The heating and cooling of the medium forces the medium to expand and then contract. This expansion and contraction generates a pressure field inside the medium, which then propagates outwards in the form of ultrasound. One can then measure the ultrasound signals on the surface of the medium. \RED{Interested readers are referred to~\cite{Beard-IF11,CoLaBe-SPIE09,Kuchment-MLLE12,KuKu-HMMI10,LiWa-PMB09,PaSc-IP07,Scherzer-Book10,Wang-DM04,Wang-IEEE08} references therein for overviews of the physical principles as well as the practical applications of PAT.} 

\RED{The objective of PAT is to use the measured ultrasound data to recover information on the interior optical properties of the underlying medium. This is often done in a two-phase manner. In the first phase, one reconstructs the initial pressure field inside the medium from measured ultrasound signals. This is a well-studied process (see~\cite{Kuchment-MLLE12} and references therein) and is often what PAT reconstruction means. The initial pressure field is related to the optical properties of the medium, the light intensity, and the so-called Gr\"unesen coefficient of the medium; see~\eqref{EQ:QPAT Data} in the next section. %
In the second phase, we reconstruct the optical properties of the medium from the internal initial pressure field data we reconstructed in the first phase. This phase is often called quantitative photoacoustic (QPAT)~\cite{GaOsZh-LNM12,BaRe-IP11, Zemp-AO10,HaNeRa-IP15,SaTaCoAr-IP13,CoLaBe-SPIE09,HaPuArTa-IPI24,TaCo-JBO24,SuHoSuQi-PMB23}.}

In this paper, we study inverse problems in the coupled system of QPAT and DOT. To be more precise, we study image reconstruction problems in QPAT with additional boundary optical current measurements given from DOT. The coupling idea was also presented in \cite{NyPuTa-BOE17}.
Nevertheless, we propose a novel three-stage image reconstruction method for the combined modality and present numerical simulation results based on synthetic data to demonstrate the performance of the proposed method. We show that the addition of optical measurements allows us to uniquely reconstruct three optical properties, the Gr\"uneisen coefficient $\Gamma$, the absorption coefficient $\sigma$, and the scattering coefficient $\gamma$, simultaneously. When the Gr\"uneisen coefficient is known, we show that the addition of the optical measurements provides a more balanced reconstruction of the absorption and the scattering coefficients.

The rest of the paper is organized as follows. In Section \ref{SEC:Model}, we briefly review the mathematical models of the coupled QPAT and DOT problem, followed by the compatibility condition. We then review the theoretical results on the QPAT and DOT problem, and discuss the uniqueness result under the coupled problem in Section \ref{SEC:Theory}. In Section \ref{SEC:Alg}, we present the three-stage method for reconstructions based on the coupled system.
Numerical simulations based on synthetic data in two dimensions will be provided in Section \ref{SEC:Num} to validate the method and demonstrate the quality of the reconstructions.
Concluding remarks are then offered in Section \ref{SEC:Con}.

\section{The mathematical models}
\label{SEC:Model}

The propagation of NIR photons in optically heterogeneous mediums is often modeled by diffusion processes. Let us denote by $\Omega\subset \bbR^d\ (d=2,3)$ the medium of interest with smooth boundary $\partial\Omega$.

\paragraph{The diffusion model in DOT.} Let $u(t, \bx)$ be the density of photons at location $\bx$ and time $t$. Then $u$ solves the following diffusion equation~\cite{Arridge-IP99}
\begin{equation}\label{EQ:Diff DOT Time}
	\begin{array}{rcll}
  	\dfrac{1}{c}\dfrac{\partial u}{\partial t}-\nabla\cdot \gamma \nabla u +\sigma u &=& 0 &\mbox{in}\ \ \bbR_+\times \Omega,\\
     u+  \kappa\gamma \dfrac{\partial u}{\partial n} &=& g(t,\bx)& \mbox{on}\ \ \bbR_+\times\partial\Omega,\\
    u(0,\bx) &=& 0 & \mbox{in}\ \ \Omega,
	\end{array}	
\end{equation}
where $c$ is the speed of light in the medium, $\gamma$ is the scattering coefficient, $\sigma$ is the absorption coefficient, $\kappa$ is the rescaled (nondimensionalized) extrapolation length, and $g$ is a given source at the boundary $\partial\Omega$ probing the medium.

The datum we measure in DOT is the photon current through the surface of the medium to be probed, that is
\begin{equation}\label{EQ:DOT Bdary Current}
	J_g(t, \bx) : = \left.\left(\gamma\dfrac{\partial u}{\partial n}\right)\right|_{(0, T] \times \partial\Omega},
\end{equation}
in a sufficiently long time interval $(0, T]$. Note that because of the Robin boundary condition in \eqref{EQ:Diff DOT Time}, the current datum $J_g$ is equivalent to the boundary photon density datum $u_{|(0, T] \times\partial\Omega}$. 

To simplify the presentation, we denote the map from boundary probing source $g$ to boundary current $J_g$ as $\Lambda^{{\rm DOT}}_{\gamma, \sigma}$. We can therefore write $J_g=\Lambda^{\text{DOT}}_{\gamma, \sigma} g$.  The objective of DOT is to reconstruct the coefficient $\gamma$ and $\sigma$ from the map %
\begin{equation}\label{EQ:DOT Data Time}
	\Lambda^{\text{DOT}}_{\gamma, \sigma}: g \mapsto J_g.
\end{equation}
In practice, one can only put sources at a finite number of locations on the boundary. Assume there are $N_s$ illuminating sources on the boundary, then the problem becomes to reconstruct the coefficients from $N_s$ pairs of data $\{g_s, J_{g_s}\}_{s=1}^{N_s}$.

\paragraph{The diffusion model in PAT.} The same diffusion equation holds for the short pulse of NIR photons sent into the medium in a PAT experiment. The derivation of the model can be found in~\cite{BaUh-IP10}. In a nutshell, the initial pressure field generated by the photoacoustic effect is determined by the total energy absorbed by the medium during the short pulse. Let $U(\bx) = \dint_0^{+\infty} u(t,\bx) dt$ and $G(\bx) = \dint_0^{+\infty} g(t,\bx) dt$. The diffusion for PAT is the equation for $U$ that we can obtain by integrating~\eqref{EQ:Diff DOT Time} over the time variable with the fact that $u(+\infty, \bx)=0$ on $\Omega$. We verify that $U$ solves the following stationary diffusion equation
\begin{equation}\label{EQ:Diff PAT}
	\begin{array}{rcll}
  	-\nabla\cdot \gamma \nabla U + \sigma U &=& 0 &\mbox{in}\ \ \Omega,\\
    U(\bx) +  \kappa\gamma \dfrac{\partial U}{\partial n} &=& G(\bx)& \mbox{on}\ \partial\Omega.
	\end{array}	
\end{equation}
\RED{
The initial pressure field inside the medium generated by the photoacoustic effect is proportional to the local energy absorbed by the medium. It can be written as
\begin{equation}\label{EQ:QPAT Data}
	H_g=\left.\Gamma(\bx) \sigma(\bx) U(\bx)\right|_{\bar\Omega},
\end{equation}
where $\sigma(\bx) U(\bx)$ is the local absorbed energy at $\bx\in\Omega$, and $\Gamma$ is a dimensionless thermodynamic constant that describes the photoacoustic efficiency of the underlying medium, measuring the conversion efficiency of heat energy to pressure. This initial pressure field then propagates as an ultrasound signal, often modeled by the acoustic wave equation, and is detected by receivers located on the surface of the medium. The measured ultrasound datum is then used to reconstruct $H$ in the first stage of PAT. In the rest of the work, we assume that the first step of PAT is done, and we are given $H$ as the internal datum we can use directly.}

We denote the map from the same boundary probing source $g$ as in the DOT problem \eqref{EQ:Diff DOT Time} to internal pressure distribution $H_g$ as $\Lambda_{\Gamma,\gamma, \sigma}^{{\rm QPAT}}$. Then we can write $H_g=\Lambda_{\Gamma, \gamma, \sigma}^{\text{QPAT}} g$.  The objective of QPAT is to reconstruct the coefficient $\Gamma$,$\gamma$ and $\sigma$ from the map %
\begin{equation}\label{EQ:QPAT Data Map}
	\Lambda^{\text{QPAT}}_{\Gamma,\gamma, \sigma}: g \mapsto H_g.
\end{equation}

\paragraph{The compatibility condition for the coupling.} The inverse problem we consider in this work is to couple the DOT inverse problem with the QPAT inverse problem. \RED{We aim at reconstructing information on $\Gamma$, $\gamma$, and $\sigma$ from data encoded in both~\eqref{EQ:DOT Data Time} and~\eqref{EQ:QPAT Data Map}.} To use both types of data, we need to ensure that they are compatible with each other since they are generated from the same source, and the light propagation processes are modeled by the same diffusion equation. For this reason, we integrate the diffusion equation~\eqref{EQ:Diff DOT Time} in the domain $\bbR_+\times\Omega$ and apply the divergence theorem to obtain 
\[
	\int_{\Omega}\dfrac{1}{c} \Big[u(t,\bx)\Big|_0^{+\infty} \Big] d\bx-\int_{\partial\Omega} \int_0^{+\infty} \bn\cdot \gamma\nabla u\, dS(\bx) dt+\int_\Omega\int_0^{+\infty} \sigma(\bx) u(t, \bx)\, dt d\bx=0.
\]
This, together with the initial condition of $u$ and the fact that $u(+\infty, \bx)=0$, leads to the compatibility condition:
\begin{equation}\label{EQ:Compat Cond}
	\int_\Omega \frac{H_g(\bx)}{\Gamma(\bx)}\, d\bx -\int_{\partial\Omega} \int_0^{+\infty} J_g(t,\bx) \, dt dS(\bx)=0.
\end{equation}
Since this compatibility condition originates from the physical model of photon propagation,

We finish this section with the following remarks.
\RED{
\begin{remark}[Nature of the coupling strategy]
    Photoacoustic tomography itself is a hybrid imaging that couples ultrasound imaging with diffuse optical imaging. PAT has also been coupled to many other imaging modalities, such as optical coherence tomography and elastography, to achieve additional advantages; see, for instance, ~\cite{ElMiSc-M2AS17,SiTh-JBO19,GaScWi-JMIV15,YaChWaSh-IEEE24} and references therein. The coupling of QPAT with DOT that we study here is slightly different in spirit as we are not using any physics to couple PAT with DOT. Instead, we are only utilizing the fact that the internal data reconstructed in PAT (i.e., $H$ in~\eqref{EQ:QPAT Data}) and the measured data in DOT are related to the same optical properties of the underlying medium. In some sense, our coupling is more of the nature of data fusion.
\end{remark}
\begin{remark}[On boundary conditions] 
We have assumed in the presentation so far that the illuminating sources used in QPAT and those used in DOT are the same. This is necessary for {\bf Stage I} of the algorithm we present in~\Cref{SEC:Alg} to work. However, this is not a super restrictive condition to assume. In practice, it is often the case that the illuminating optical source $g(t, \bx)$ in~\eqref{EQ:Diff DOT Time} is only a temporal modulation of a spatial pattern $\wt g(\bx)$ given in the laser source used, that is, $g(t, \bx)=\wt g(\bx) h(t)$. In the photoacoustic case, $h(t)=\delta(t-0_+)$~\cite{BaUh-IP10}. Therefore, $G(\bx)=\wt g(\bx)$ in~\eqref{EQ:Diff PAT}. In the stationary DOT, $h(t)=1$, and $g(t, \bx)=\wt g(\bx)$. Therefore, to require that the sources in DOT be the same as those in QPAT is essentially the same as requiring that we use the same spatial component $\wt g$ (controlled by the laser source used), not the time modulation mechanism, in both modalities, which is not hard to implement. 
\end{remark}
}

\section{Theoretical considerations on the coupling}
\label{SEC:Theory}

We review briefly some theoretical results on DOT and PAT and their implications on our coupled inverse problem. Recall that the coefficients to be reconstructed are $\Gamma(\bx)$, $\gamma(\bx)$ and $\sigma(\bx)$, and we assume throughout the rest of the paper that $\gamma(\bx)\in \cC^2(\bar\Omega)$ for simplicity.

In the frequency domain, the diffusion model in DOT ~\eqref{EQ:Diff DOT Time}, via the Fourier transformation in time, can be written as 
\begin{equation}\label{EQ:Diff DOT Freq}
	\begin{array}{rcll}
  	i\dfrac{\omega}{c}\wh u(\omega, \bx)-\nabla\cdot \gamma \nabla\wh  u +\sigma\wh u &=& 0 &\mbox{in}\ \ \Xi\times\Omega,\\
    \wh u(\omega, \bx) +  \kappa\gamma \dfrac{\partial \wh u}{\partial n} &=& \wh g(\omega, \bx)& \mbox{on}\ \ \Xi\times\partial\Omega,
	\end{array}	
\end{equation}
where $\Xi=[\omega_{\rm min}, \omega_{\rm max}]$ is the interval of modulation frequency that is accessible.

For any pair of $\gamma(\bx)$ and $u(t,\bx)$, we introduce the Liouville transform $v(t, \bx)=\sqrt{\gamma} u(t, \bx)$. Then the above DOT problem can be written into the standard Schr\"{o}dinger form:
\begin{equation}\label{EQ:Diff DOT Freq v}
\begin{array}{rcll}
	\Delta \wh v(\omega, \bx) - (\dfrac{\Delta\sqrt \gamma}{\sqrt \gamma} + \dfrac{\sigma}{\gamma}+ \dfrac{i\omega}{c\gamma})\wh v& = & 0 & \mbox{in}\ \Xi\times\Omega, \\
\wh v+\kappa \gamma^{3/2}\bn\cdot\nabla\dfrac{\wh v}{\sqrt{\gamma}} &= & \sqrt{\gamma} \,\wh g(\omega, \bx) &\mbox{on} \ \Xi\times\partial\Omega,
\end{array}
\end{equation}
with the boundary current data in the form
\[
 \wh J_g(\omega, \bx)  = \left.\left(\bn\cdot \gamma \nabla \dfrac{\wh v}{\sqrt{\gamma}}\right)\right|_{\Xi \times \partial\Omega}.
\]
The PAT diffusion problem~\eqref{EQ:Diff PAT}, under the similar Liouville transform $V(\bx)=\sqrt{\gamma}U(\bx)$, takes the form:
\begin{displaymath}
\begin{array}{rcll}
\Delta V - (\dfrac{\Delta\sqrt \gamma}{\sqrt \gamma} + \dfrac{\sigma}{\gamma})V & = & 0 &\mbox{in}\ \ \Omega, \\
V+\kappa \gamma^{3/2}\bn\cdot\nabla\dfrac{V}{\sqrt{\gamma}} &= & \sqrt{\gamma}G(\bx) &\mbox{on} \  \partial\Omega,
\end{array}
\end{displaymath}
with the internal data $H$ of the form
\[
	H_g(\bx)=\left.\dfrac{\Gamma\sigma}{\sqrt{\gamma}} V\right|_{\bar\Omega}.
\]
To simplify the notations, we further define
\RED{
\begin{equation}\label{EQ:Potentials}
	q:=\dfrac{\Delta\sqrt \gamma}{\sqrt \gamma} + \dfrac{\sigma}{\gamma} \qquad \mbox{and}\qquad \pi:=\dfrac{\Gamma\sigma}{\sqrt{\gamma}}.
\end{equation}
}
The results of Sylvester and Uhlmann~\cite{SyUh-AM87} in dimension three and Nachman~\cite{Nachman-AM96} in dimension two state that the full Dirichlet-to-Neumann map, that is, the map \RED{$\wh\Lambda_{\gamma,\sigma}(\omega):\wh g(\omega,\cdot)
        \mapsto \wh J_g$}
with $\kappa=0$, uniquely determines the complex potential
\begin{equation}\label{EQ:Potential Complex}
	\dfrac{\Delta \sqrt{\gamma}}{\sqrt{\gamma}}+ \dfrac{\sigma}{\gamma}+\dfrac{i\omega}{c \gamma}=q(\bx)+\dfrac{i\omega}{c \gamma}.
\end{equation}
Arridge and Lionheart~\cite{ArLi-OL98}, based on previous results in~\cite{SyUh-AM87}, observed that the real part of this potential gives $q(\bx)$ while the imaginary part of the potential gives $\dfrac{i\omega}{c \gamma}$. Therefore, if $\omega\neq 0$, $\gamma$ and $q$ (and therefore $\sigma$) are uniquely reconstructed from the DOT data. This means that when $\omega\neq 0$, one can uniquely reconstruct $\gamma$ and $\sigma$ from the Dirichlet-to-Neumann data. The stability of the reconstruction, however, is only logarithmic; see, for instance, the discussions in~\cite{Mandache-IP01,Isakov-Book06} and references therein. When time-domain measurements are available, one can uniquely reconstruct $(\gamma, \sigma)$ with Lipschitz stability if the initial condition in~\eqref{EQ:Diff DOT Time} is known (which is hard to do in practice) and non-zero~\cite{ImYa-IP98}. We refer to~\cite{Uhlmann-BMS14,BeFrVe-SIAM21,Foschiatti-JMAA24,Alessandrini-AA88,Mandache-IP01,AlVe-AAM05,Rondi-AAM06} and references therein for more detailed discussions on similar inverse problems.

Meanwhile, it has been shown in~\cite{BaUh-IP10,BaRe-IP11}, when $\kappa=0$, the data encoded in the map $\Lambda^{\rm QPAT}_{\Gamma, \gamma, \sigma}: g \mapsto H_g$ are enough to uniquely and stably, in appropriate norms, reconstruct the potential $q$ and $\pi$ when $\gamma_{|\partial\Omega}$ is known. 
When the boundary value of $\gamma$, $\gamma_{|\partial\Omega}$ is known, one can reconstruct any two of the three coefficients from two data sets. When $\gamma_{|\partial\Omega}$ is not known, it is shown in~\cite{ReGaZh-SIAM13} that $\gamma$ can not be uniquely reconstructed from given internal data. 

\paragraph{Uniqueness for Coupled Inverse Problem.}  Following the above discussion, we conclude that the coupled inverse problem of reconstruction $(\Gamma, \gamma, \sigma)$ from $\Lambda_{\gamma, \sigma}^{\rm DOT}$ and $\Lambda_{\Gamma, \gamma, \sigma}^{\rm QPAT}$ has the following properties:

(i) When $\kappa=0$, that is, with Dirichlet boundary conditions, one can uniquely reconstruct $\Gamma$, $\gamma$, and $\sigma$ simultaneously from the map $\Lambda_{\gamma, \sigma}^{\rm DOT}$ and one internal data $H_g$. This is because the data $\Lambda_{\gamma, \sigma}^{\rm DOT}$ allow the unique reconstruction of $\gamma$ and $\sigma$. Once they are reconstructed, the corresponding solution $U$ is also reconstructed. Therefore, the additional internal data $H$ will allow the unique reconstruction of $\Gamma$. 

(ii) When $\kappa$ is nonzero, the DOT and PAT data from two well-chosen illuminations $g_1$ and $g_2$, $\{(g_j, J_{g_j}), H_{g_j}\}_{j=1}^2$, are sufficient to uniquely determine $(\Gamma, \gamma, \sigma)$ as a generalization of the result in \cite{ReGaZh-SIAM13}. 

Both (i) and (ii) provide the uniqueness of reconstructing three coefficients with data at a single optical wavelength, which are different from the results in~\cite{BaRe-IP12,MaRe-CMS14} where illuminations at multiple optical wavelengths are used, together with additional \emph{a priori} information, to achieve similar uniqueness of reconstructions.

\section{Three-stage reconstruction of \texorpdfstring{$(\Gamma, \gamma, \sigma)$}{}}
\label{SEC:Alg}

The analysis in Section~\ref{SEC:Theory} guarantees that the coupling allows a unique reconstruction of $\Gamma$, $\gamma$ and $\sigma$ without using multi-spectral data, which motivates us to develop methods to reconstruct the three parameters simultaneously. In this section, we propose a three-stage image reconstruction method for the coupled DOT-PAT system. In the first stage, we reconstruct the value of $\gamma$ on the boundary of the domain. In the second stage, we reconstruct $\gamma$ and $\sigma$ in the entire domain. In the third stage, we reconstruct the  Gr\"uneisen coefficient $\Gamma$. We assume that we have data collected from $N_s$ sources for all modulation frequencies. The three-stage method is summarized as Algorithm~\ref{ALG:Three-Stage}.

\paragraph{Stage I: reconstructing $\gamma_{|\partial\Omega}$.} To reconstruct $\gamma$ on the boundary, we use two pairs of data $\{ (g_j, J_j), H_j \}_{j=1}^2$, and observe that
\[
	\gamma U_2^2\bn\cdot\nabla\dfrac{H_1}{H_2}=\gamma U_2^2\bn\cdot\nabla\dfrac{U_1}{U_2}=U_2\gamma \bn\cdot\nabla U_1- U_1\gamma \bn\cdot\nabla U_2.
\]
Using the boundary conditions \eqref{EQ:DOT Bdary Current} and in \eqref{EQ:Diff PAT}, we obtain the relationships on the boundary by integrating these conditions over time: $ \bn\cdot\nabla U_i =\int_0^{\infty}J_idt $ and $  U_i = G_i =\int_0^{\infty}g_idt $.

Therefore, we obtain the analytic formula for the reconstruction of $\gamma$ on the boundary
\begin{equation}\label{EQ:Gamma Bdary Rec}
	\gamma 
 = \dfrac{U_2\bn\cdot\nabla U_1-U_1\bn\cdot\nabla U_2}{U_2^2\bn\cdot\nabla\dfrac{H_1}{H_2}} 
 =\dfrac{\int_0^{\infty}g_2\,dt\int_0^{\infty}J_1\,dt-\int_0^{\infty}g_1\,dt\int_0^{\infty}J_2\,dt}{(\int_0^{\infty}g_2\,dt)^2\, \bn\cdot\nabla\dfrac{H_1}{H_2}} \; 
\end{equation}
provided that $\bn\cdot\nabla\dfrac{H_1}{H_2}\neq 0$ on $\partial\Omega$.

\paragraph{Stage II: reconstructing $(\gamma, \sigma)$.} We first solve the following minimization problem:
\begin{equation}\label{couple}
    (\gamma, \sigma)=\argmin_{\gamma,\sigma} \Phi(\gamma,\sigma),
\end{equation}
where the objective function $\Phi$ has three different components: one from the misfit of the internal PAT observations $H_i$ to the data $H_i^*$, one from the misfit of the boundary DOT data, and the third from regularization. More precisely, it is written as
\begin{equation}\label{EQ:OBJ}
    \Phi(\gamma,\sigma)=\Phi_{\rm PAT}(\gamma,\sigma)+\Phi_{\rm DOT}(\gamma,\sigma)+  \edited{\cR(\gamma, \sigma)},
\end{equation}
where 
\begin{equation}\label{EQ:Phi PAT}
\Phi_{\rm PAT}(\gamma, \sigma) = \dfrac{1}{2}\sum_{i=2}^{N_s}\sum_{j<i} \int_{\Omega} ( H_i^*H_j-H_i H_j^*)^2 d\bx\,,
\end{equation}
\begin{equation}\label{phi-dot}
\Phi_{\rm DOT}(\gamma, \sigma) =
\dfrac{1}{2}\sum_{j=1}^{N_s}\sum_{i=1}^{N_w} \int_{\partial\Omega}\left|\gamma \dfrac{\partial u_j^{\omega_i}}{\partial n}-J_j^{\omega_i}\right|^2 dS(\bx) \,,
\end{equation}
\RED{with $N_\omega$ denotes the total number of frequencies in $\Xi$} and the regularization functional $\cR(\gamma, \sigma)$ is taken to be of the Tikhonov type
\begin{equation}\label{reg}
\edited{	\cR(\gamma, \sigma) = \frac{1}{2}\Big(\beta_{\gamma}\|\nabla \gamma\|_{L^2(\Omega)}^2 +\beta_{\sigma}\|\nabla \sigma\|_{L^2(\Omega)}^2\Big)\,.}
\end{equation}
The reason why we only reconstruct $\gamma$ and $\sigma$ from minimizing $\Phi$ is that although the function itself depends on the unknown coefficient $\Gamma$, its minimizers do not. This is the key observation that leads to our three-stage method. We provide more discussions on this issue in~\Cref{SEC:Sensitivity}.

\paragraph{Stage III: reconstructing $\Gamma$.} Once the coefficients $\gamma$ and $\sigma$ are reconstructed, we reconstruct the Gr\"uneisen coefficient $\Gamma$ as
\begin{equation}\label{stage3}
	\Gamma(\bx) = \dfrac{\sum_{j=1}^{N_s}H_j^*(\bx)}{\sigma(\bx) \sum_{j=1}^{N_s} U_j(\bx)}.
\end{equation}
Due to the fact that optical signals decay fast away from the source location, signals measured on detectors far from the source are very weak. This average-then-reconstruct method is, in general, better than the reconstruct-then-average approach $\Gamma=\frac{1}{N_s+1}\sum_{j=1}^{N_s}\frac{H_j(\bx)}{\sigma(\bx) U_j(\bx)}$ since $\sum_{j=1}^{N_s} U_j(\bx)$ is less likely to be zero than each of the individual $U_j$.

Combining the three stages discussed above, we obtain the following \Cref{ALG:Three-Stage} for simultaneous reconstruction of $\sigma, \gamma, \Gamma$. 
\begin{algorithm}
\caption{Three-stage DOT-PAT Algorithm}
\label{ALG:Three-Stage}
\begin{algorithmic}[1]
\Require{Data: $\{g_j, J_j, H_j^*\}_{j=1}^{N_s}$ }
\State Stage I: \RED{Analytic recovery of $\gamma_{|\partial\Omega}$ using~\eqref{EQ:Gamma Bdary Rec}}.
\State Initialize guesses for $(\sigma_0$, $\gamma_0, \Gamma_0)$.
\State Stage II: Reconstruct $(\gamma,\sigma)$ from solving the optimization problem \eqref{couple}. \label{alg:opt}
\State Stage III: Recover $\Gamma$ by evaluating the average formula \eqref{stage3}.
\State Update $\sigma$ by optimizing $\Phi_{\text{PAT}}$ under constraints \eqref{EQ:Diff PAT}  with fixed $\gamma$ and $\Gamma_0$ in \eqref{EQ:QPAT Data}. \label{alg:update-sigma}
\Ensure{reconstruction of $\gamma, \sigma,\Gamma$}
\end{algorithmic}
\end{algorithm}

The extra step in \Cref{alg:update-sigma} of ~\Cref{ALG:Three-Stage} is there because we observe that $\sigma$ is only mildly updated when solving \eqref{alg:opt}, thus we need an extra step to fully update and recover $\sigma$. We refer to the numerical experiment sections to further discuss the reason and effect of this step.

In solving the optimization problems in ~\Cref{alg:update-sigma,alg:opt} of~\Cref{ALG:Three-Stage}, 
we use the BFGS quasi-Newton method~\cite{KlHi-IP03}, which requires gradient computations. The gradients of $\Phi$ with respect to $\gamma$ and $\sigma$ can be derived from the standard adjoint method~\cite{HiPiUlUl-Book08}. We omit the derivation process as they can be found adapted from calculations in existing calculations as those in~\cite{ReGaZh-SIAM13}. We summarize the results in the following lemma.
\begin{lemma}\label{thm}
Assume sufficient regularity, $\Phi$ is Fr\'echet differentiable with respect to $\gamma$ and $\sigma$. The derivatives are given as, assuming that $\gamma$ is given on the boundary, as
\begin{equation}
	\Phi'(\gamma, \sigma)[\delta \gamma]=\Phi_{\rm PAT}'(\gamma, \sigma)[\delta \gamma]+\Phi_{\rm DOT}'(\gamma, \sigma)[\delta \gamma]+ \edited{ \cR'(\gamma, \sigma)[\delta\gamma]},
\end{equation}
\begin{equation}
	\Phi'(\gamma, \sigma)[\delta \sigma]=\Phi_{\rm PAT}'(\gamma, \sigma)[\delta \sigma]+\Phi_{\rm DOT}'(\gamma, \sigma)[\delta \sigma]+\edited{\cR'(\gamma, \sigma)[\delta\sigma]},
\end{equation}
where the three parts of the contributions are respectively given as follows. 

(i) The PAT part:
\begin{equation}\label{eq:PAT-grad-gamma}
	\Phi_{\rm PAT}'(\gamma, \sigma)[\delta \gamma]=\sum_{i=2}^{N_s}\sum_{j<i}\int_{\Omega} \delta\gamma \big[\nabla W_i\cdot \nabla U_j - \nabla \wt W_j \cdot \nabla U_i\big] d\bx- \int_{\partial\Omega} \delta\gamma \big[W_i\dfrac{\partial U_j}{\partial n}-\wt W_j\dfrac{\partial U_i}{\partial n}\big]dS(\bx),
\end{equation}
\begin{equation}\label{eq:PAT-grad-sigma}
	\Phi_{\rm PAT}'(\gamma, \sigma)[\delta \sigma]=\sum_{i=2}^{N_s}\sum_{j<i}\int_{\Omega} \delta\sigma Z_{ij}\Gamma\big[H_i^* U_j-U_iH_j^*\big] d\bx + \int_{\Omega} \delta\sigma\big[W_i U_j - \wt W_j U_i\big] d\bx,
\end{equation}
where we have analytic reconstruction of $\gamma_{|\partial\Omega}$ from proposed stage I and so $\delta\gamma_{|\partial\Omega}=0$, $W_{i}$ and $\wt W_{j}$ are respectively solutions to 
\begin{equation}
-\nabla \cdot \gamma \nabla W_{i} + \sigma W_{i}=-Z_{ij} H_i^*\Gamma\sigma\ \  \mbox{in}\ \Omega, \qquad 
W_{i}+\kappa\gamma\dfrac{\partial W_{i}}{\partial n}=0\ \ \mbox{on}\ \partial\Omega.
\end{equation}
\begin{equation}
-\nabla\cdot\gamma\nabla\wt W_{j} + \sigma \wt W_{j}=-Z_{ij} H_j^*\Gamma\sigma\ \  \mbox{in}\ \Omega, \qquad 
\wt W_{j} + \kappa\gamma\dfrac{\partial \wt W_{j}}{\partial n} =0\ \ \mbox{on}\ \partial\Omega,
\end{equation}
with $Z_{ij}=(H_i^*H_j-H_iH_j^*)/H_i^*H_j^*$.

(ii) The DOT part:
\begin{equation}
\begin{aligned}
	\Phi_{\rm DOT}'(\gamma, \sigma)[\delta \gamma]&= {\rm Re}\Big(\sum_{j=1}^{N_s}\sum_{i=1}^{N_w} \int_{\Omega} \delta\gamma  \nabla \overline{w_j}^i\cdot \nabla u_j^i\,d\bx\Big)\,,
 \end{aligned}
\end{equation}
\begin{equation}
	\Phi_{\rm DOT}'(\gamma, \sigma)[\delta \sigma]={\rm Re}\Big(\sum_{j=1}^{N_s} \sum_{i=1}^{N_w} \int_{\Omega} \delta\sigma \overline w_j^i u_j^i\, d\bx\Big)\,,
\end{equation}
where $\{\overline{w_j}^i\}_{i=0,j=0}^{N_w,N_s}$ are complex conjugates of $\{w_j^i\}_{i=0,j=0}^{N_w,N_s}$, with total number of frequencies $N_w$. $w_j^i$ are solutions to 
\begin{equation}
\dfrac{-i\omega_i}{c}w_j^i-\nabla\cdot\gamma\nabla w^i_j + \sigma w^i_j=0\ \  \mbox{in}\ \Omega, \qquad 
w^i_j + \kappa\gamma\dfrac{\partial w^i_j}{\partial n} =z_j^i\ \ \mbox{on}\ \partial\Omega,
\end{equation}
with $z_j^i=(\gamma \dfrac{\partial u_j^{\omega_i}}{\partial n}-J_j^{\omega_i})/|J_j^{\omega_i}|$. 
\\[1ex]
(iii) The regularization part
\begin{equation}
	\edited{\cR'(\gamma, \sigma)[\delta\gamma, \delta\sigma]=\left[\beta_{\gamma}\int_{\Omega}\delta\gamma \Delta \gamma\,d\bx, \beta_{\sigma}\int_{\Omega}\delta\sigma \Delta \sigma \,d\bx\right].}
\end{equation}
\end{lemma}

To save computational cost in practical uses, we consider an alternative objective function for the PAT part:
\begin{equation}\label{alternative}
\Phi_{\text{PAT}}(\gamma, \sigma) = \dfrac{1}{2}\sum_{j=2}^{N_s}\int_{\Omega} (H_1^*H_j-H_j H_1^*)^2 d\bx,
\end{equation}
which fixes one of the reference data in $(H_i^*H_j-H_iH_j^*)$ to be $H_1$ and $H_1^*$, rather than going through all the data pairs. It reduces the double summation in the original objective function into a single summation while also factoring out $\Gamma$ as Equation\eqref{joint}.
The gradients can be calculated similarly as~\eqref{eq:PAT-grad-gamma} and~\eqref{eq:PAT-grad-sigma}. The only difference is changing the double summations to the corresponding single summations. 

\RED{
\begin{remark}
It is clear by now that {\bf Stage I} of our algorithm requires that the illuminating sources used in QPAT and those used in DOT are the same. This is the only place where such a requirement is needed. Therefore, this requirement of using the same sources in DOT and QPAT can be removed if the boundary value of $\gamma$ does not need to be reconstructed, a situation that is often assumed in the QPAT literature~\cite{BaUh-IP10,BaRe-IP11}.
\end{remark}
}
\section{Numerical experiments}
\label{SEC:Num}
We now present some numerical simulation results to demonstrate the performance of the proposed three-stage reconstruction algorithm.

The data used are synthetic data that are created by solving the forward problems \eqref{EQ:Diff DOT Freq} and \eqref{EQ:Diff PAT} on a regular mesh with a two-dimensional rectangular geometry and coefficients that are set to mimic the real data range. The domain is taken to be $[0,2]\times[0,2]$, covered by a finite element mesh of $6561$ nodes and $12800$ triangle elements. The results are interpolated on a $81 \times 81$ grid in MATLAB. In the DOT reconstruction stage, we use a frequency to the medium-speed rate, $\frac{\omega}{c}$, within $10^{-7}\times[0,2,4,5,7]$ to mimic real modulation frequencies. We use the alternative objective function \eqref{alternative} %
for the PAT part to save computation. 
To ensure different terms in the objective have compatible magnitudes, we rescale the objective functional with measured data as follows: %
\begin{equation}\label{EQ:Rescaled Obj}
\begin{array}{rcl}
\Phi_{\rm PAT}(\gamma, \sigma) &=& \dfrac{1}{2}\dsum_{j=2}^{N_s}\dsum_{j<i} \int_{\Omega} \left(\dfrac{H_i^*H_j-H_i H_j^*}{H_i^*H_j^*} \right)^2 d\bx\,,\\[0.1ex]
\Phi_{\rm DOT}(\gamma, \sigma) &=&
\dfrac{1}{2}\dsum_{j=1}^{N_s}\dsum_{i=1}^{N_w} \int_{\partial\Omega}\left|(\gamma \dfrac{\partial u_j^{\omega_i}}{\partial n}-J_j^{\omega_i})/J_j^{\omega_i}\right|^2 dS(\bx) \,.
\end{array}
\end{equation}
This reweighting technique is also because that the optical signals decay fast away from the source location, and signals measured on detectors far from the source are very weak and therefore play little role in the unweighted objective function. This rescaling cancels the effect of the fast decaying of magnitudes for data far away from the source, and ensures they contribute with similar magnitudes to the objective compared with the close-to-source part of the data. 

\RED{The regularization coefficients \edited{$\beta_{\sigma}$ and $\beta_{\gamma}$} are used differently for $\sigma$ and $\gamma$, and their values are determined in a trial-and-error manner as we do not have a theory for the optimal parameter selection. The main consideration in the parameter selection process is the ratio between the magnitude of the regularization term and that of the data mismatch term. For the fairness of comparisons, we pre-determined the values of the parameters and kept them for all the experiments. Specifically, the parameter for the absorption is $\beta_{\sigma}=10^{-6}$, and the parameter for the diffusion coefficient is $\beta_\gamma=10^{-5}$ in all the simulations.
}

\RED{The data are collected from $N_s=36$ illuminating sources evenly spaced on the boundary. For each source, $320$ detectors are used in the DOT measurement. Detectors are also evenly spaced. The spatial components of the sources are the restrictions Gaussian functions on the segment of the boundary they are located in. More precisely,
\begin{equation*}
g_s(t, \bx)=\Big(1+5\, \exp\big(-\frac{1}{0.02}|\bx-\bx_s|^2\big)\Big)\delta(t)\, ,
\end{equation*}
where $\bx_s$ is the center of source $s$. %
}

In the following discussions, we will denote the data obtained from the direct solutions of DOT and PAT equations under MATLAB \texttt{assempde} as the \textit{noise-free} data though they might inherit noise from the initial guess for coefficients and the MATLAB interpolation algorithm. We reserve the term \textit{noisy data} for those data with additional additive noise by multiplying the original data with $(1+0.05\texttt{random})$ where \texttt{random} is a uniformly distributed random variable in $[-1,1]$. We fixed the noise level to be $5\%$ for later experiments. We also refer to all Dirichlet boundary conditions as the rescaled (nondimensionalized) extrapolation length $\kappa=0$ and the Robin boundary condition as $\kappa=0.2$.
 
The later optimization problems are all solved with the BFGS quasi-Newton method and updated with the adjoint-derived gradients in \Cref{SEC:Alg} with noise-free observational data. \Cref{ALG:Three-Stage} is implemented within the MATLAB Optimization Toolbox. \RED{The first-order optimality tolerance is set to be $10^{-9}$ with a max iteration to be $200$.} 
Considering the nature of the optical tomography problem is unstable and the fact that the initial point should include some information about the true coefficients, the initial guesses for all later experiments are from sending the true coefficients through the MATLAB function \texttt{imgaussfilt} which filters the truth with a 2-D Gaussian smoothing kernel under a standard deviation of $5$.
\begin{figure}[!htb]
\settoheight{\tempdima}{\includegraphics[width=.23\linewidth]{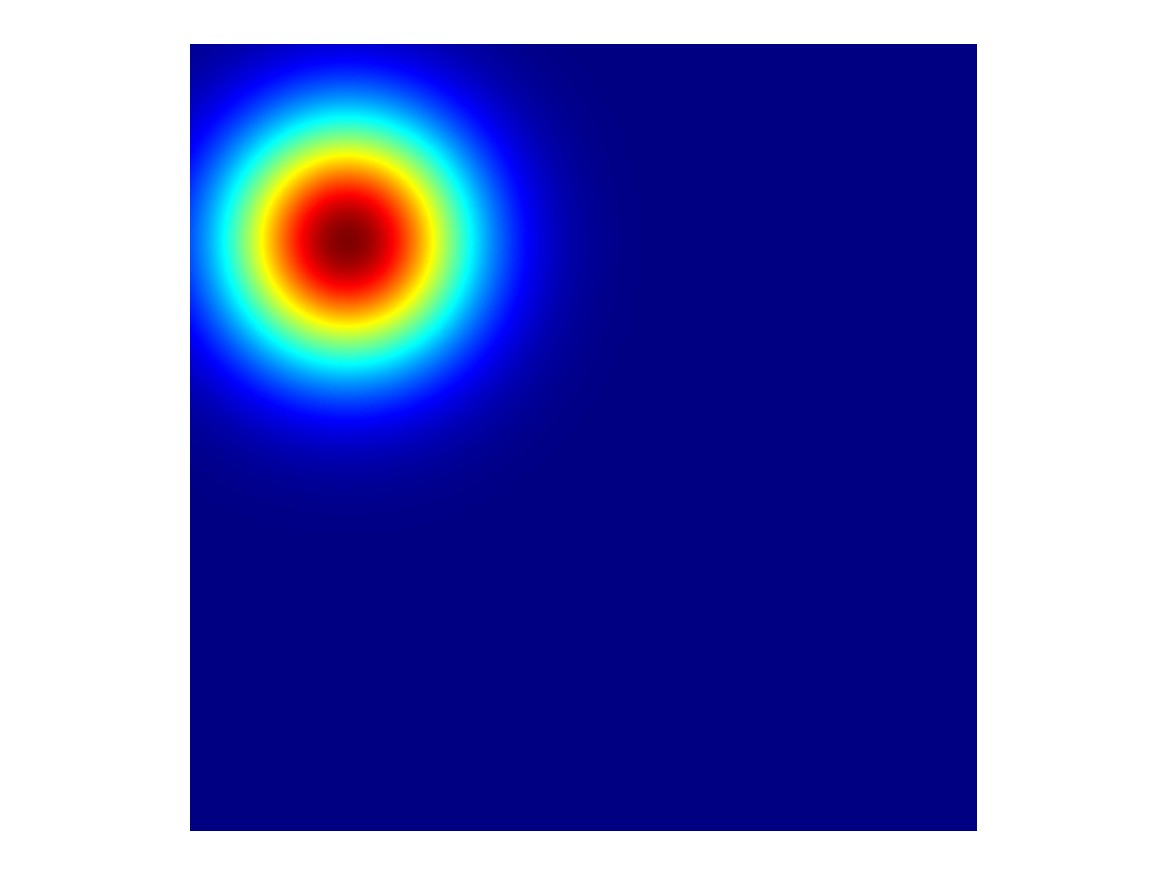}}%
\centering
\begin{tabular}{@{\hspace{-.3ex}}c@{\hspace{-0.1ex}}c@{\hspace{-3.5ex}}c@{\hspace{-2.7ex}}c@{\hspace{-1ex}}c@{\hspace{-2.7ex}}c@{\hspace{-0.7ex}}}
\rowname{$\sigma$}&
\includegraphics[width=.23\linewidth]{FinalFigures/couple/sigmat.jpg}&
\includegraphics[width=.23\linewidth]{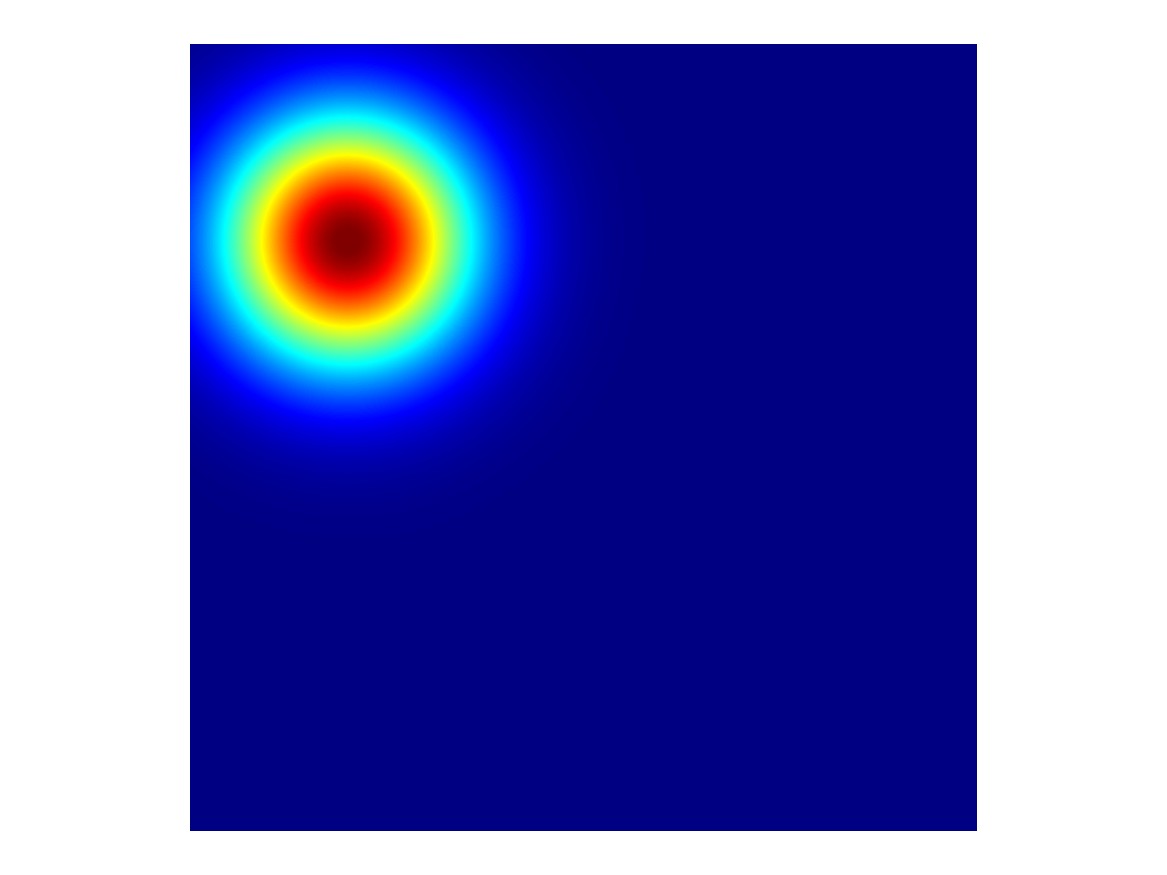}&
\includegraphics[width=.23\linewidth]{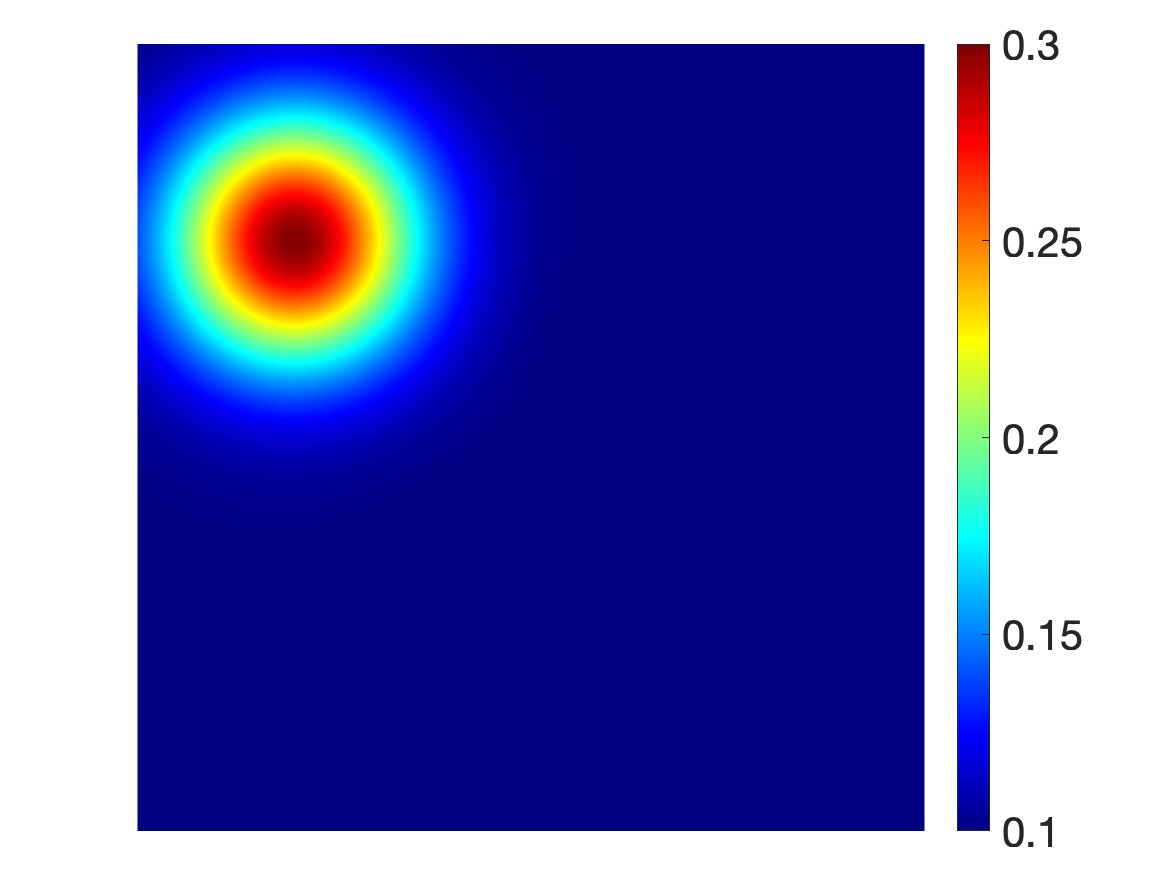}&
\includegraphics[width=.23\linewidth]{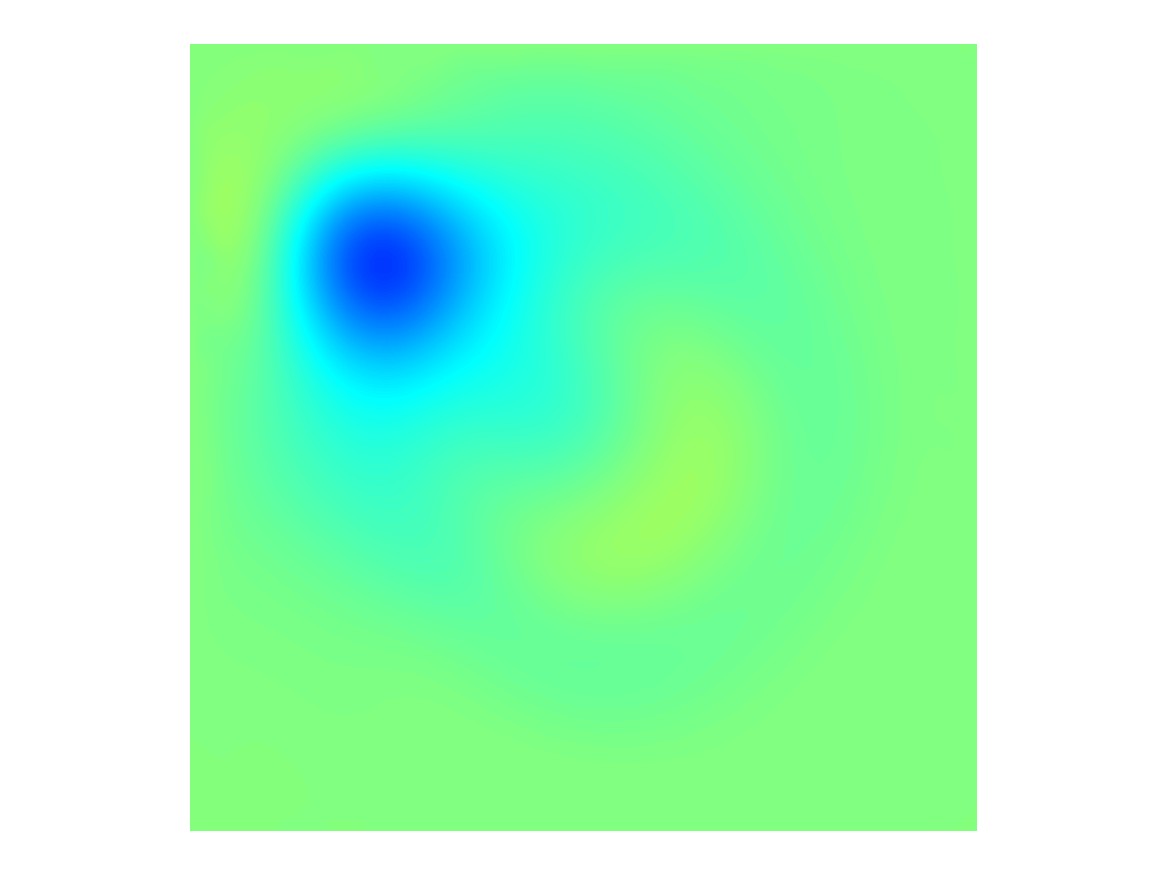}&
\includegraphics[width=.23\linewidth]{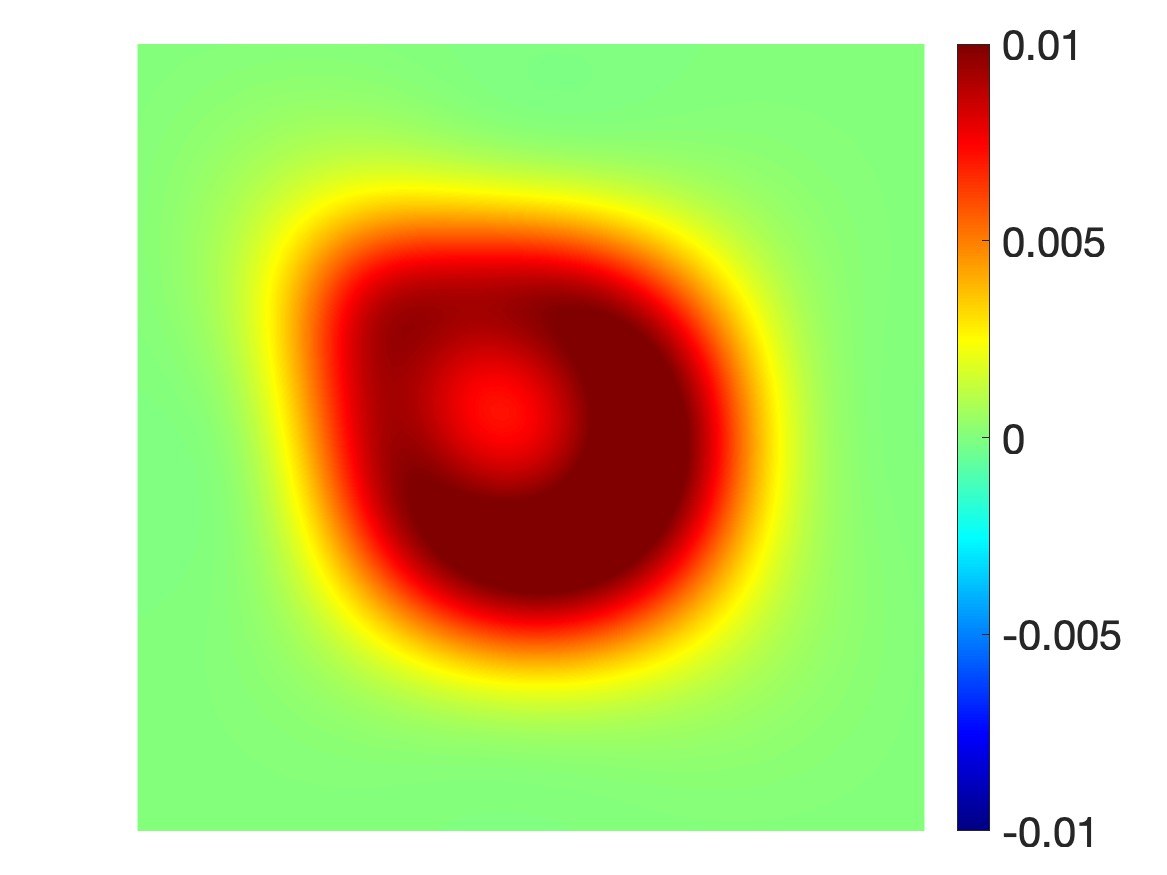}\\
 \rowname{$\gamma$}& 
\includegraphics[width=.23\linewidth]{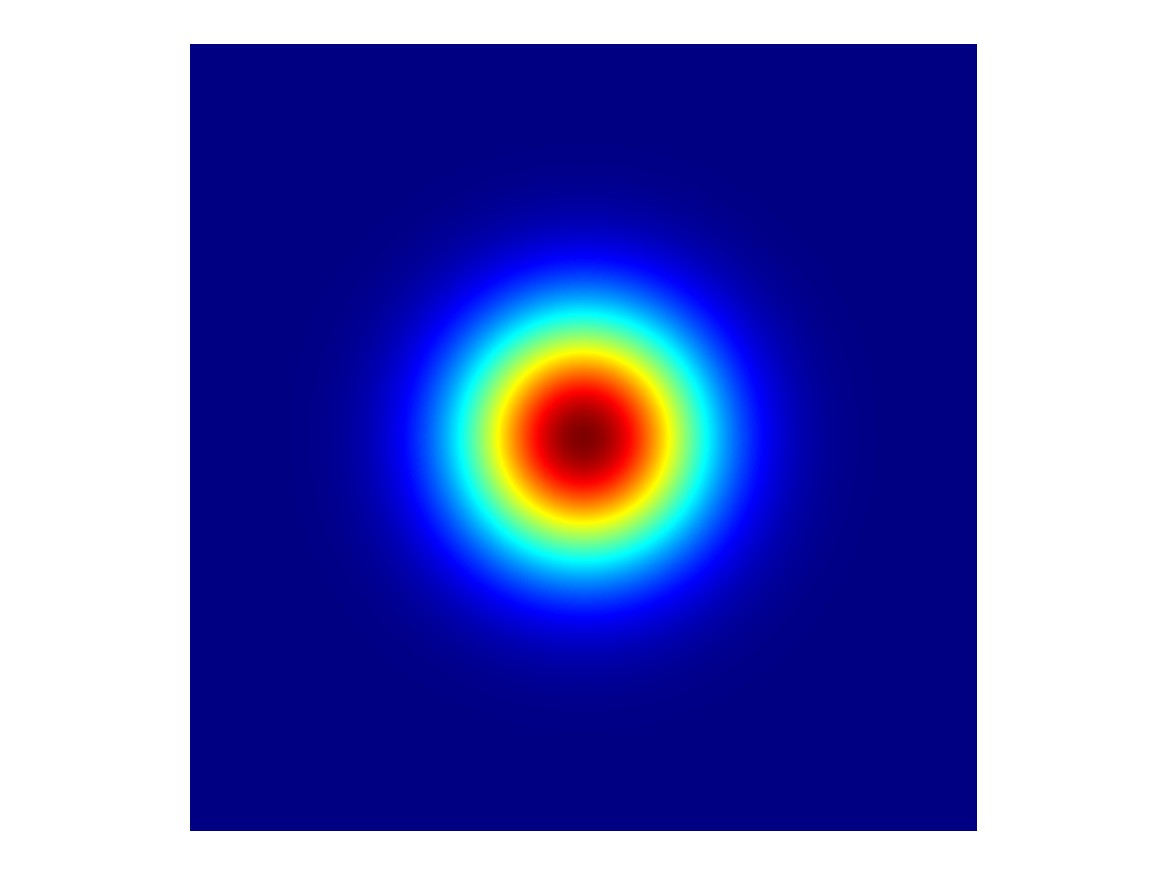}&
\includegraphics[width=.23\linewidth]{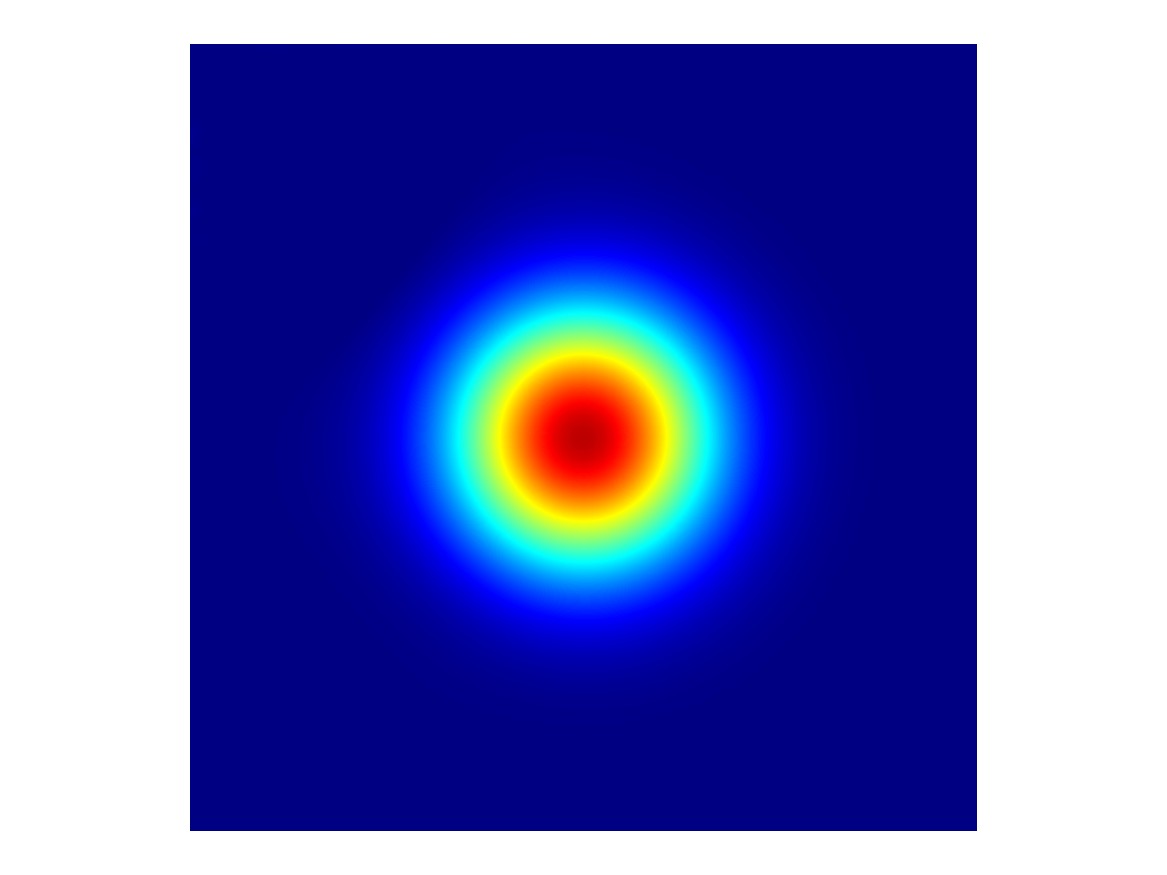}&
\includegraphics[width=.23\linewidth]{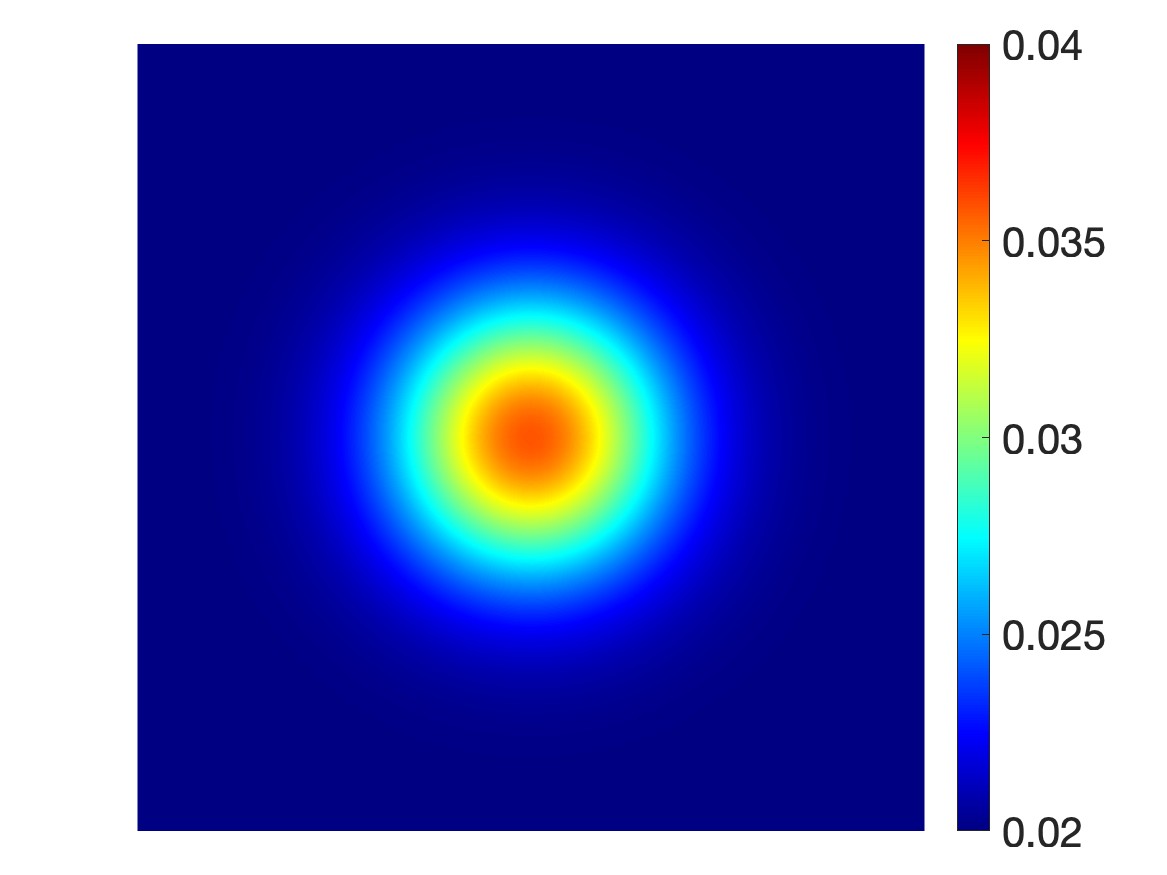}&
\includegraphics[width=.23\linewidth]{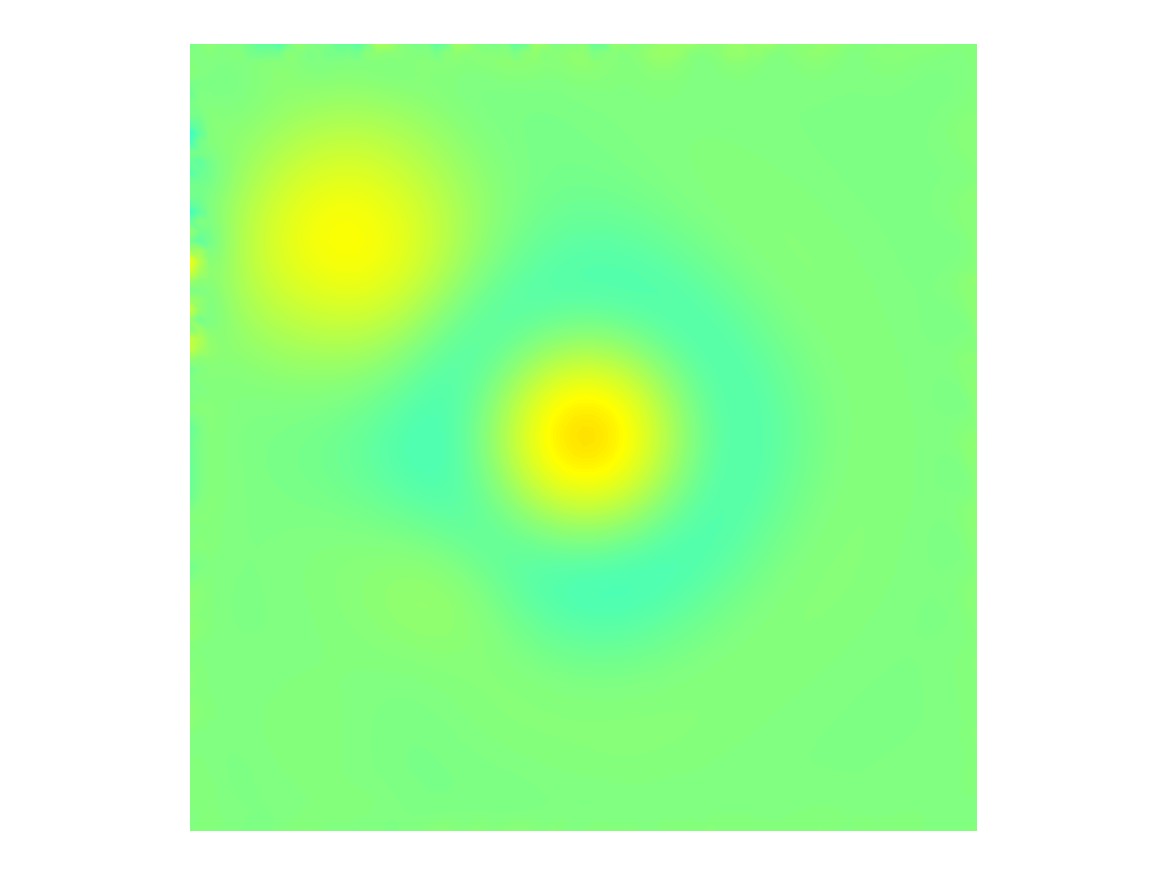}&
\includegraphics[width=.23\linewidth]{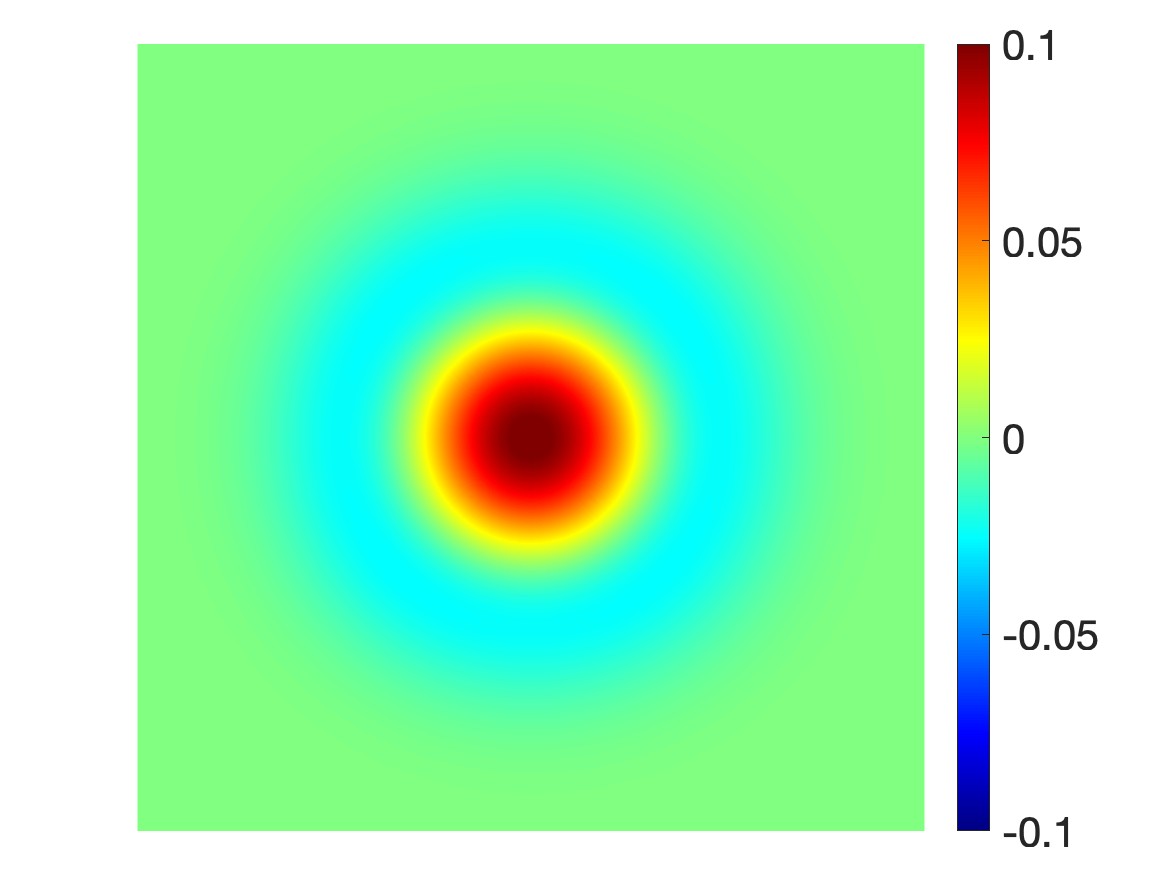}\\
& Truth & Coupling & QPAT & Coupling Error & QPAT Error
\end{tabular}
\caption{\RED{Comparison of the reconstructions using QPAT+DOT coupling data (second column) with those using only QPAT data (third column), when $\Gamma$ is known.}
}%
\label{FIG:Exp I}
\end{figure}

\RED{\paragraph{Experiment I.} We start with numerical simulations where we assume that $\Gamma$ is known already and we are only reconstructing $\sigma$ and $\gamma$. We compare the reconstruction results using both QPAT and DOT data (using our Algorithm \ref{ALG:Three-Stage} but skipping Stage III given that $\Gamma$ is known) with those where only QPAT data are used. The reconstruction in the latter case (i.e. when only QPAT data are used) is achieved with the following minimization problem
\begin{equation}
    (\sigma,\gamma)=\argmin_{\sigma,\gamma} \sum_{s=1}^{Ns} \int_\Omega \left(\dfrac{H_s-H^*_s}{H_s^*}\right)^2 d\bx + \edited{ \cR(\gamma, \sigma)},
\end{equation}
with the same regularization term as that in~\eqref{EQ:OBJ}. The true coefficients are given by
\begin{equation}
    \sigma(\bx) = 0.1 + 0.3 e^{-\frac{|\bx-\bx_1|^2}{2\times 0.05^2}}, \qquad \gamma(\bx) = 0.02 + 0.04 e^{-\frac{|\bx-\bx_2|^2}{2\times 0.05^2}}, 
\end{equation}
where $\bx_1=(0.4,1.5)$ and $\bx_2=(1.0, 1.0)$.
} 

\RED{The reconstructions of $\sigma$ (top row) and $\gamma$ (bottom row) from both the coupling data (QPAT and DOT data) and the QPAT only data are presented in~\Cref{FIG:Exp I}. Shown are the true coefficients (first column), the reconstructions with QPAT+DOT data (second column), the reconstructions with only QPAT data (third column), and the respective errors (defined as the difference between the true and the reconstructed coefficients) of the reconstructions (fourth and fifth columns). While the improvement is minimal, we do observe that the results from the coupled data reconstructions are more accurate (and this is consistently observed in other numerical experiments we did that are not presented here due to the limitation of space), as seen easily from the error plots. A careful observation of the error plots shows that both the magnitude and the spatial distribution of the errors are quite different for the two sets of reconstructions presented. It seems that the error of the coupling reconstruction is large where $\sigma$ has large values, while the error of the QPAT-only reconstruction tends to be large at locations where $\gamma$ has large values.}

\RED{A closer look at the reconstructions in~\Cref{FIG:Exp I} shows a clear sign of cross-talk effect between $\sigma$ and $\gamma$ in both the coupled and the QPAT-only reconstructions. Due to the large difference between the magnitudes of the coefficients, it is not easy to see the cross-talk in the reconstructed images. However, the error plots in the right two columns show that large errors of $\sigma$ reconstructions happen at the places where the $\gamma$ coefficient is large, and vice versa. While such cross-talk effects between $\sigma$ and $\gamma$ observed here are less severe compared to those observed in the DOT literature, it is still a clear sign that the changes in $\sigma$ and $\gamma$ can sometimes create very similar changes in QPAT and DOT data.
}

\begin{figure}[!htb]
\settoheight{\tempdima}{\includegraphics[width=.23\linewidth]{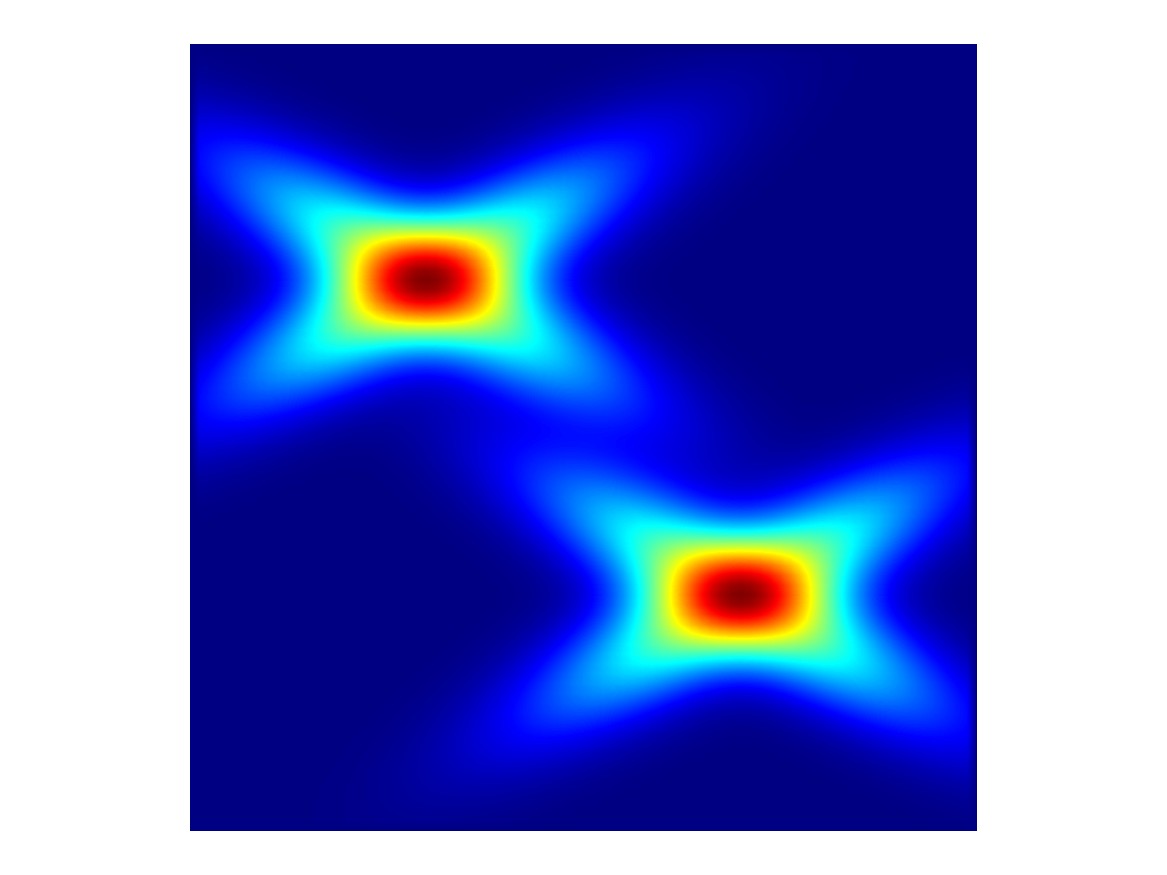}}%
\centering
\begin{tabular}{@{\hspace{-.3ex}}c@{\hspace{-0ex}}c@{\hspace{-3.5ex}}c@{\hspace{-2.7ex}}c@{\hspace{-1ex}}c@{\hspace{-2.7ex}}c@{}}
\rowname{$\sigma$}&
\includegraphics[width=.23\linewidth]{FinalFigures/GaussianMixture/sigmat.jpg}&
\includegraphics[width=.23\linewidth]{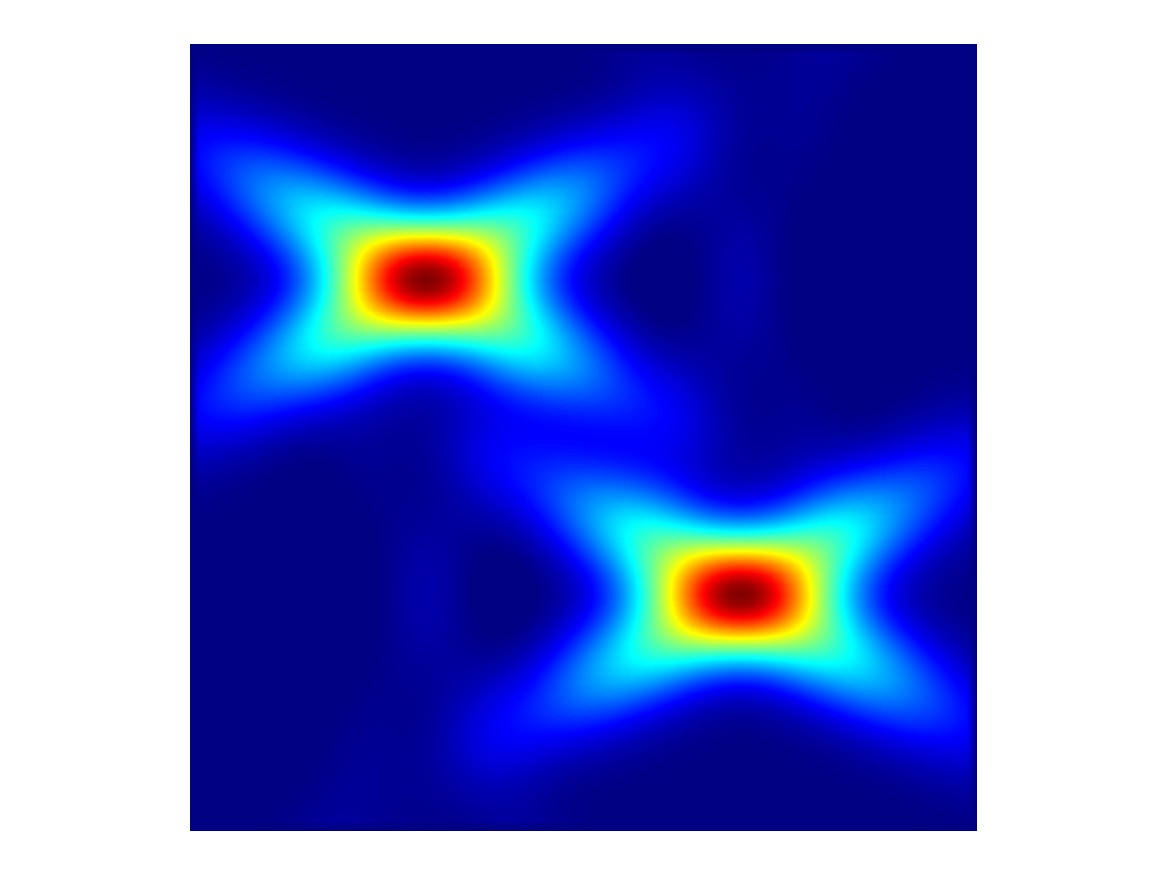}&
\includegraphics[width=.23\linewidth]{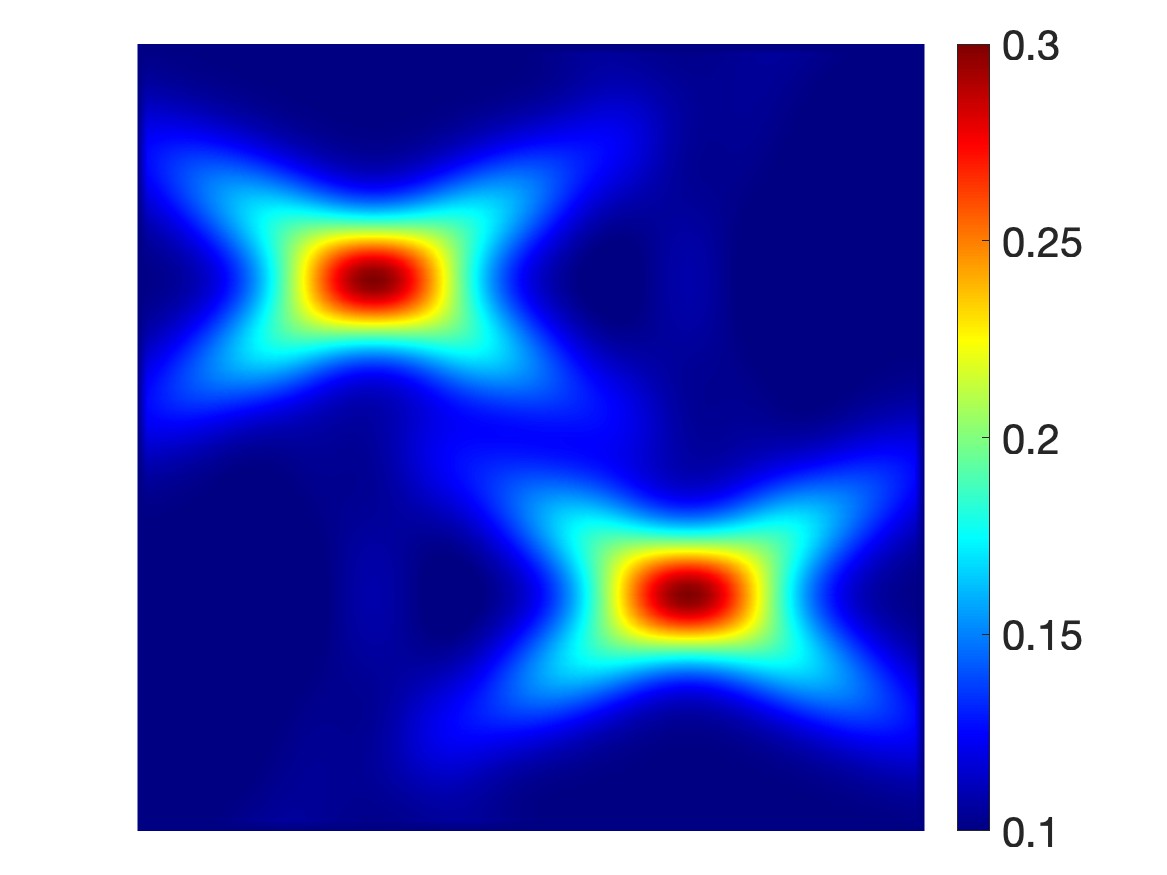}&
\includegraphics[width=.23\linewidth]{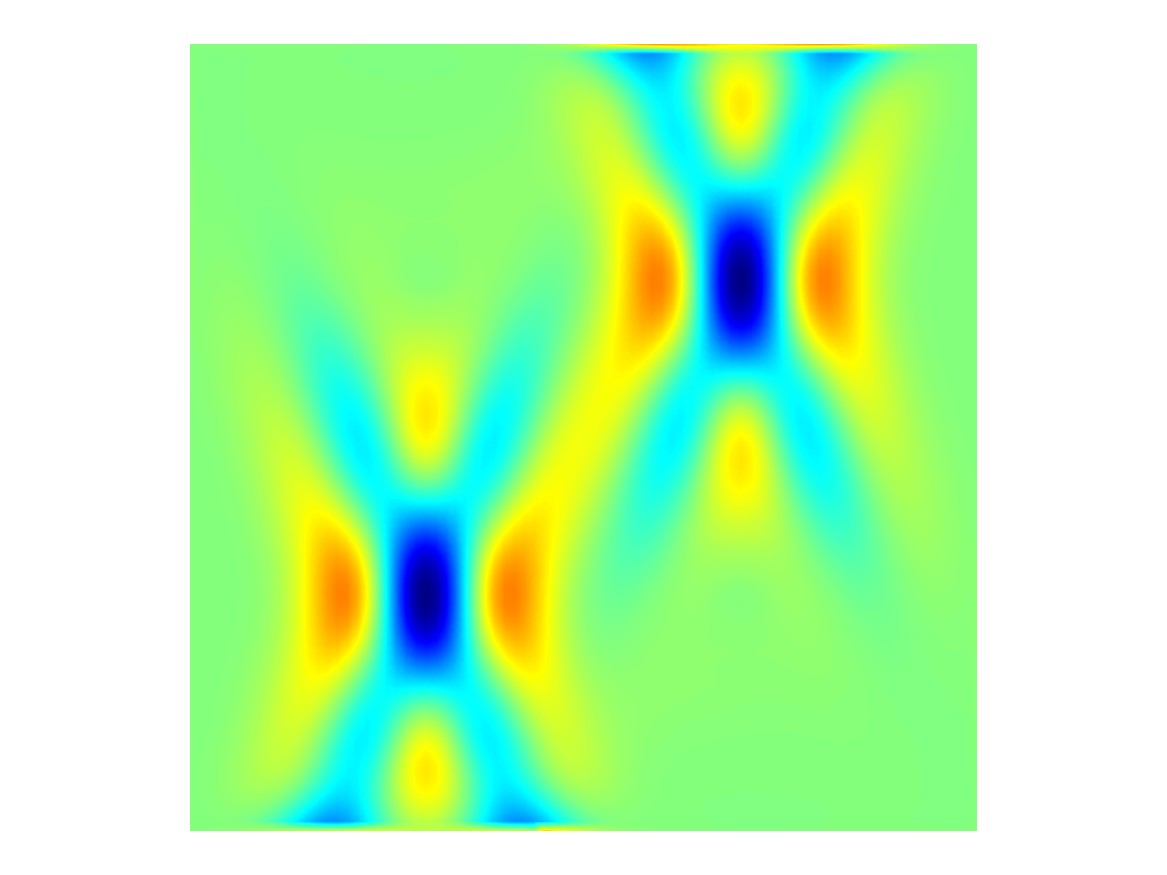}&
\includegraphics[width=.23\linewidth]{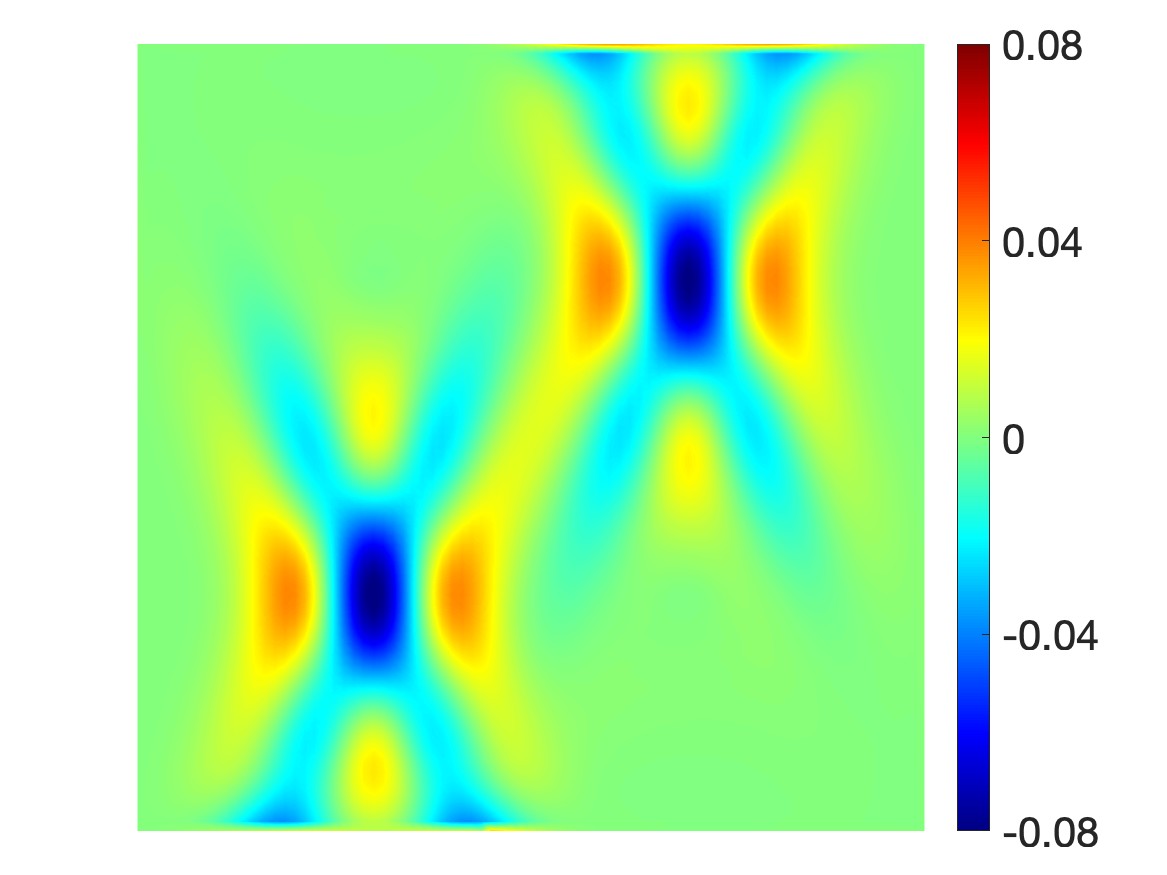}\\
\rowname{$\gamma$}&
\includegraphics[width=.23\linewidth]{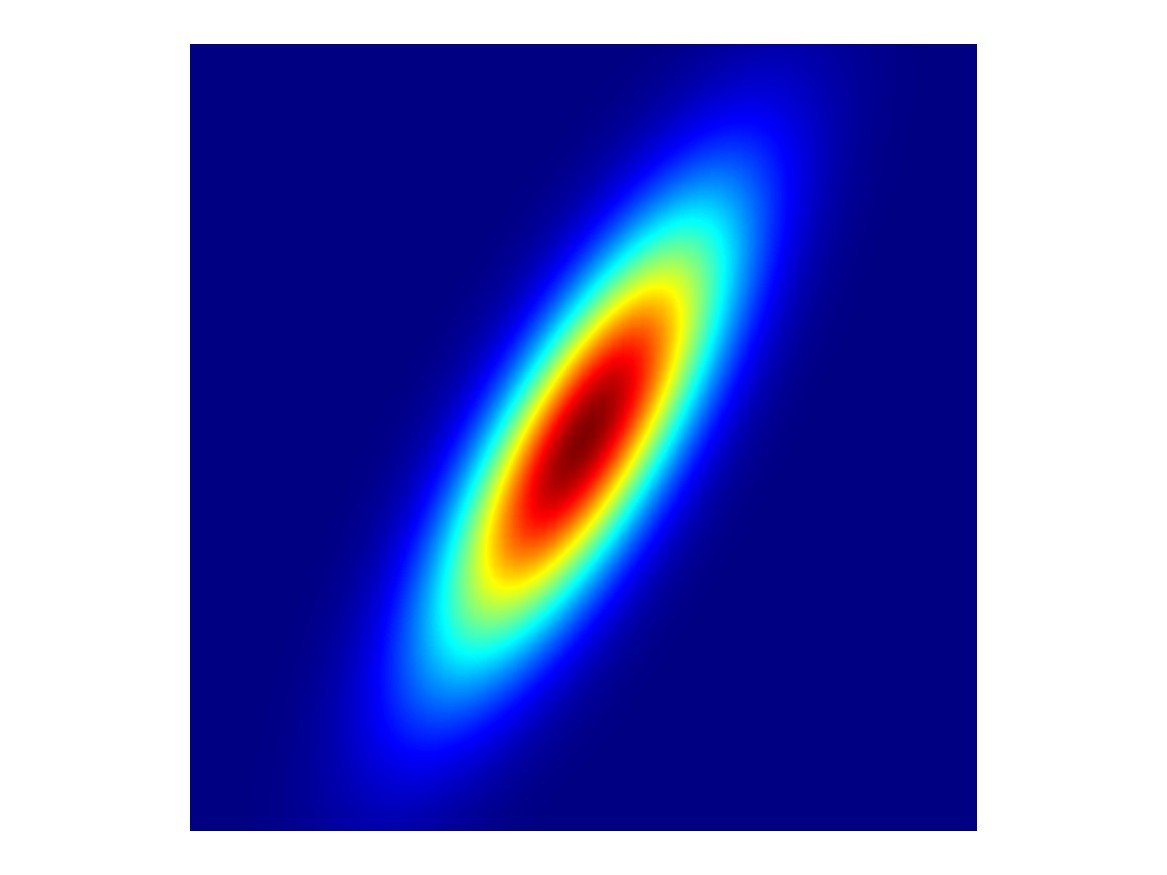}&
\includegraphics[width=.23\linewidth]{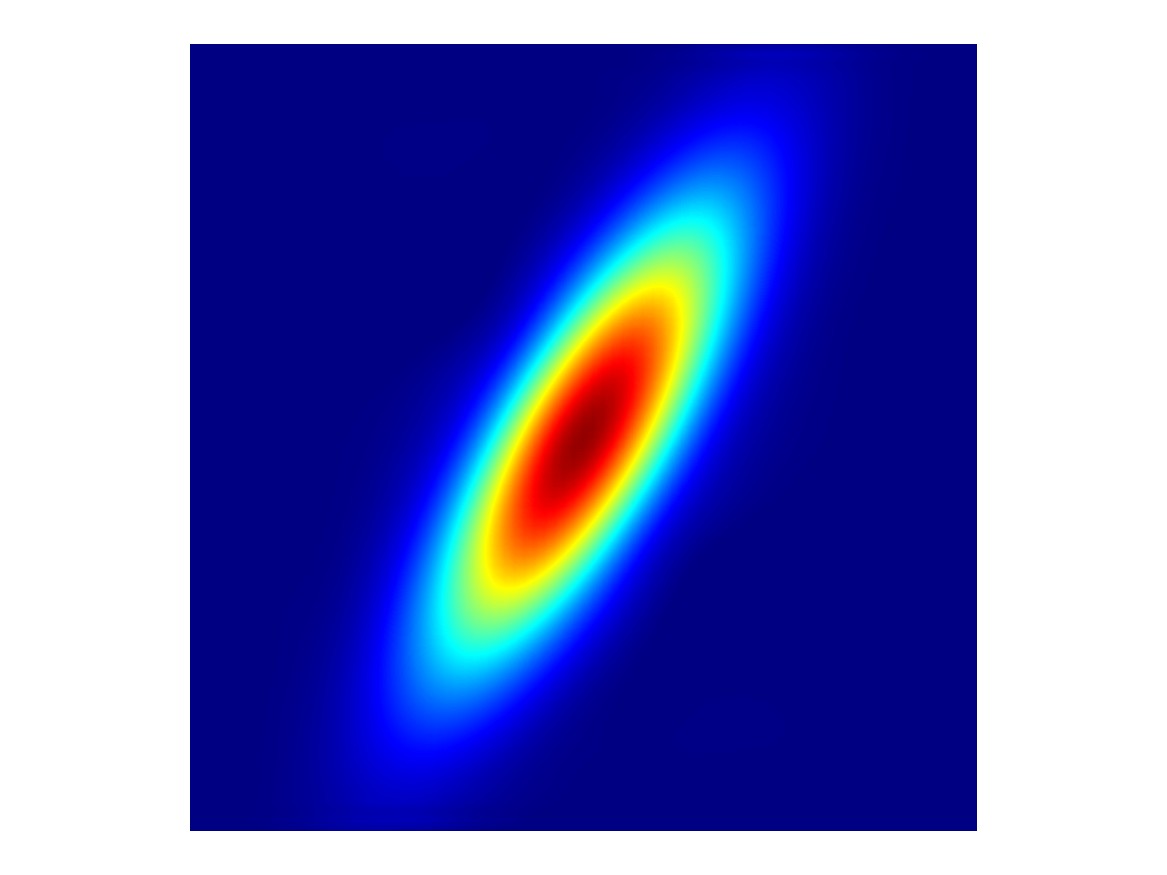}&
\includegraphics[width=.23\linewidth]{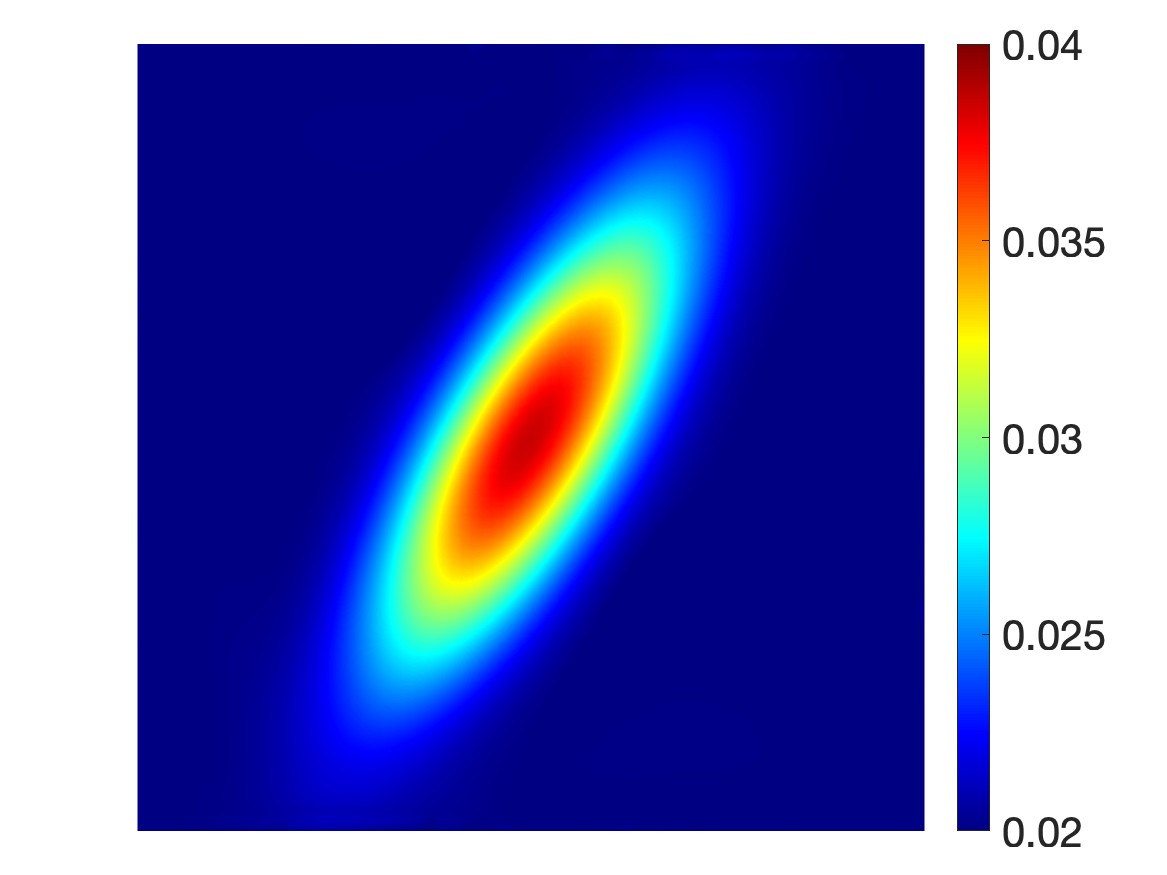}& 
\includegraphics[width=.23\linewidth]{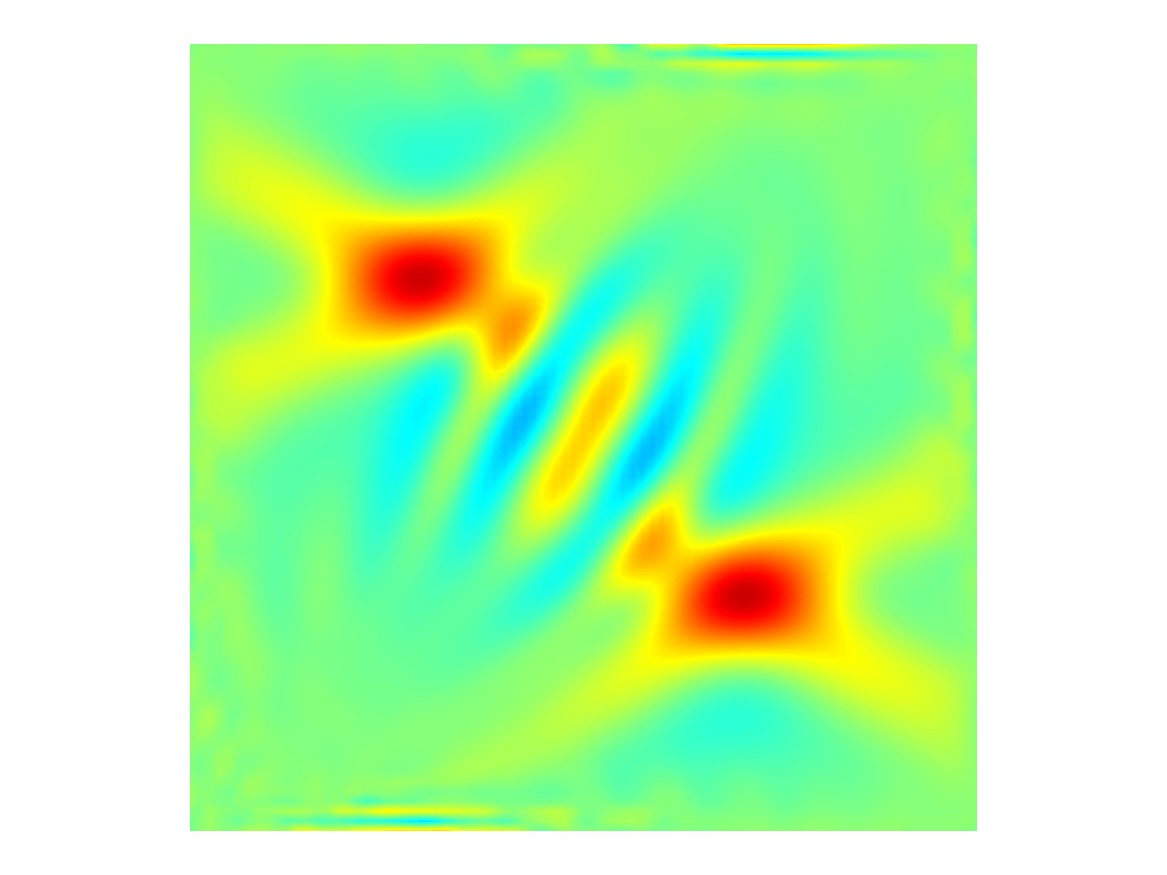}&
\includegraphics[width=.23\linewidth]{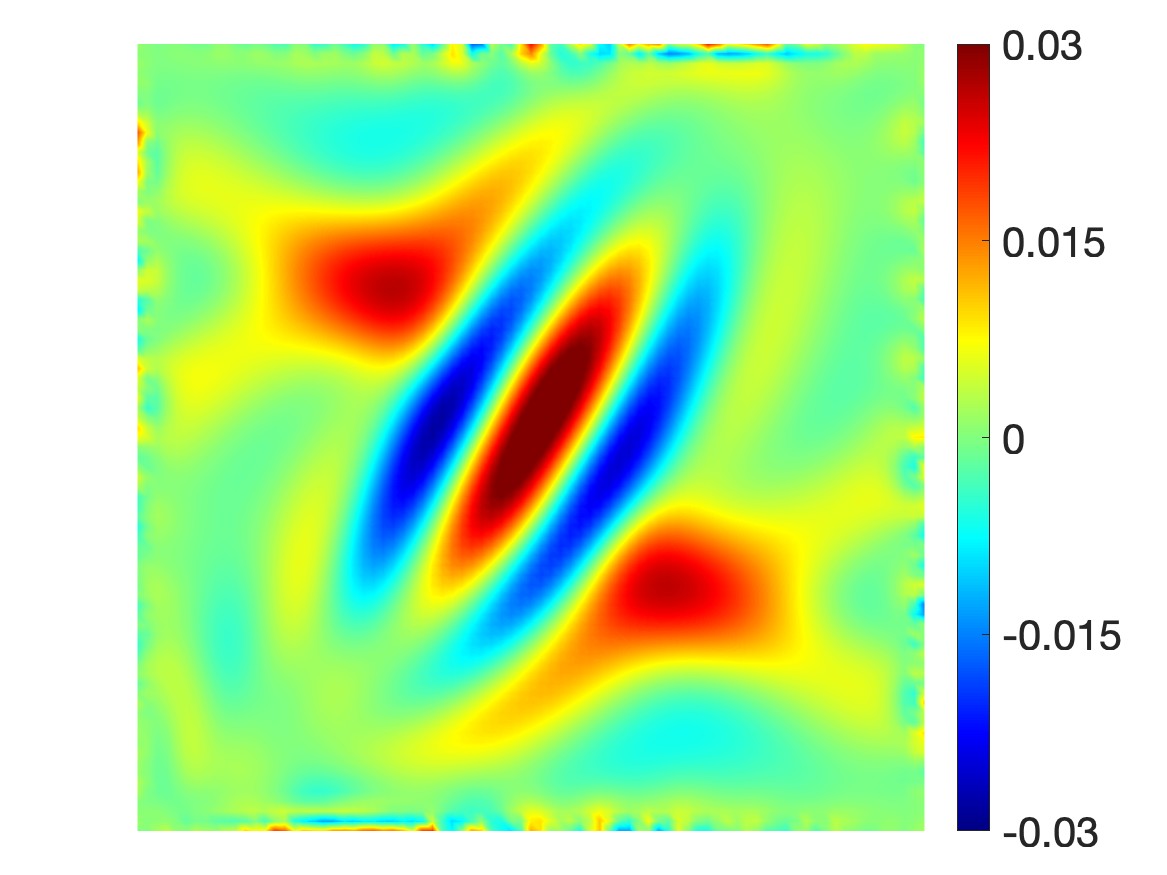}\\
\rowname{$\Gamma$}&
\includegraphics[width=.23\linewidth]{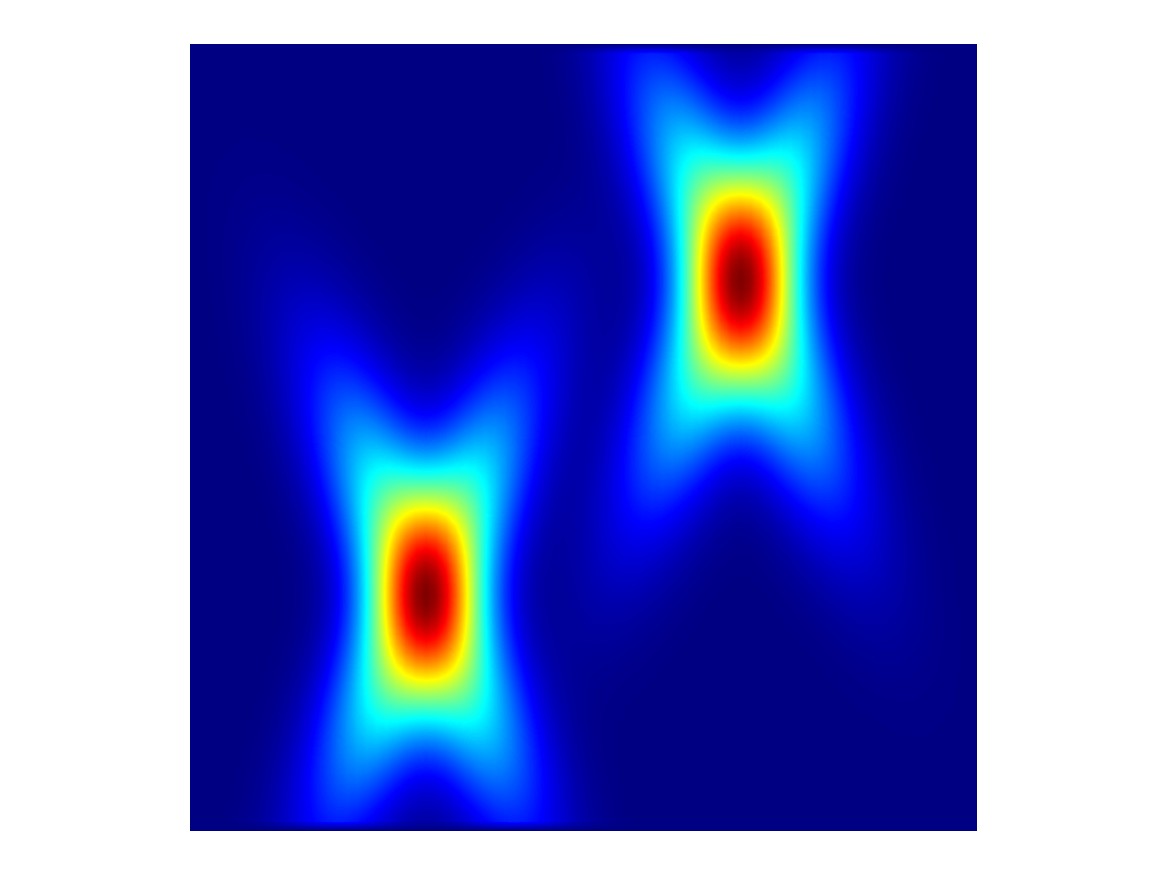}&
\includegraphics[width=.23\linewidth]{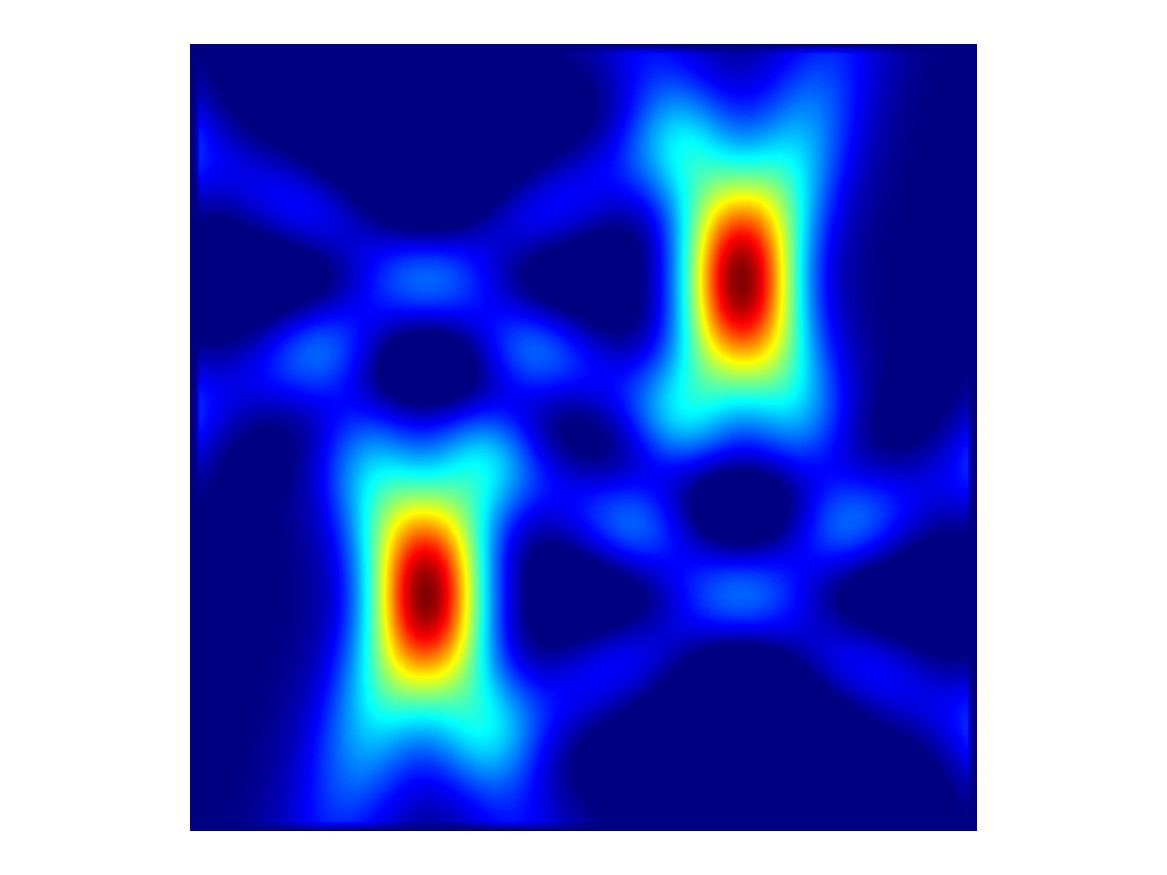}&
\includegraphics[width=.23\linewidth]{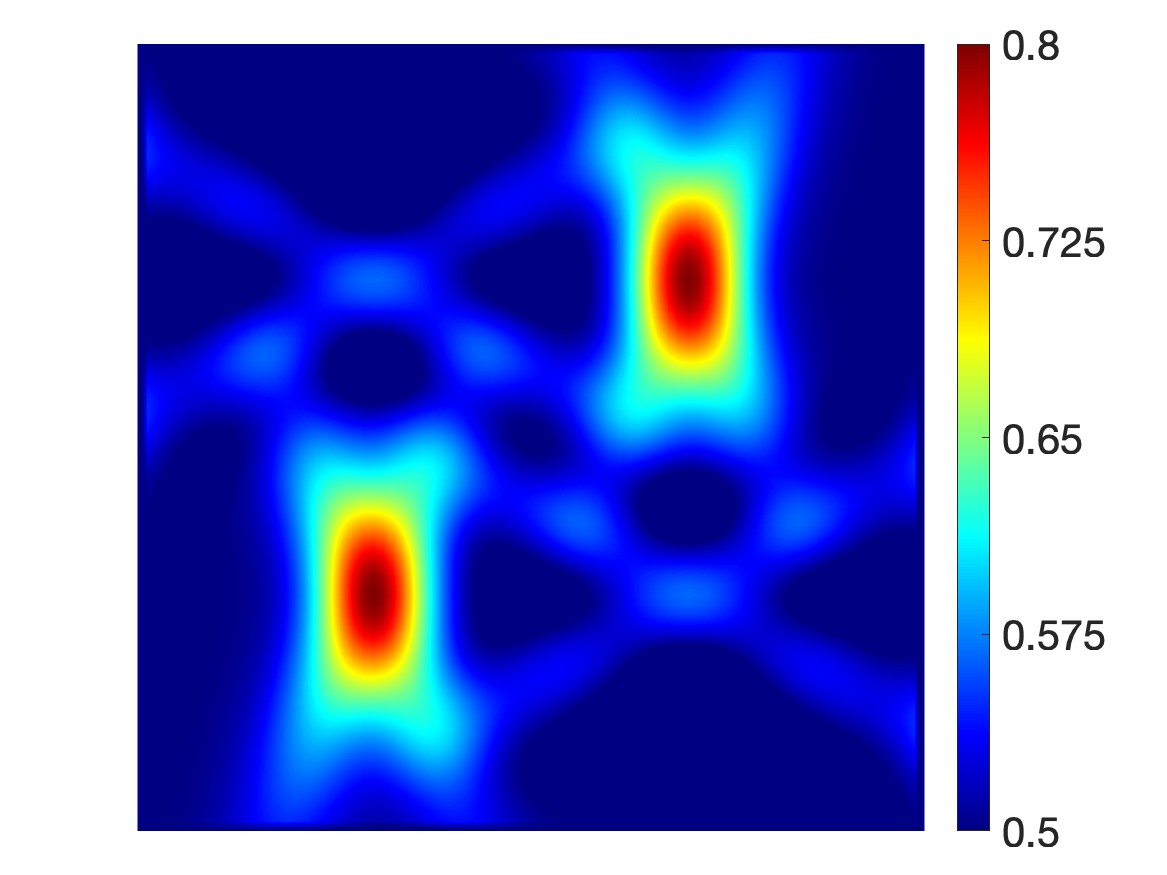}&
\includegraphics[width=.23\linewidth]{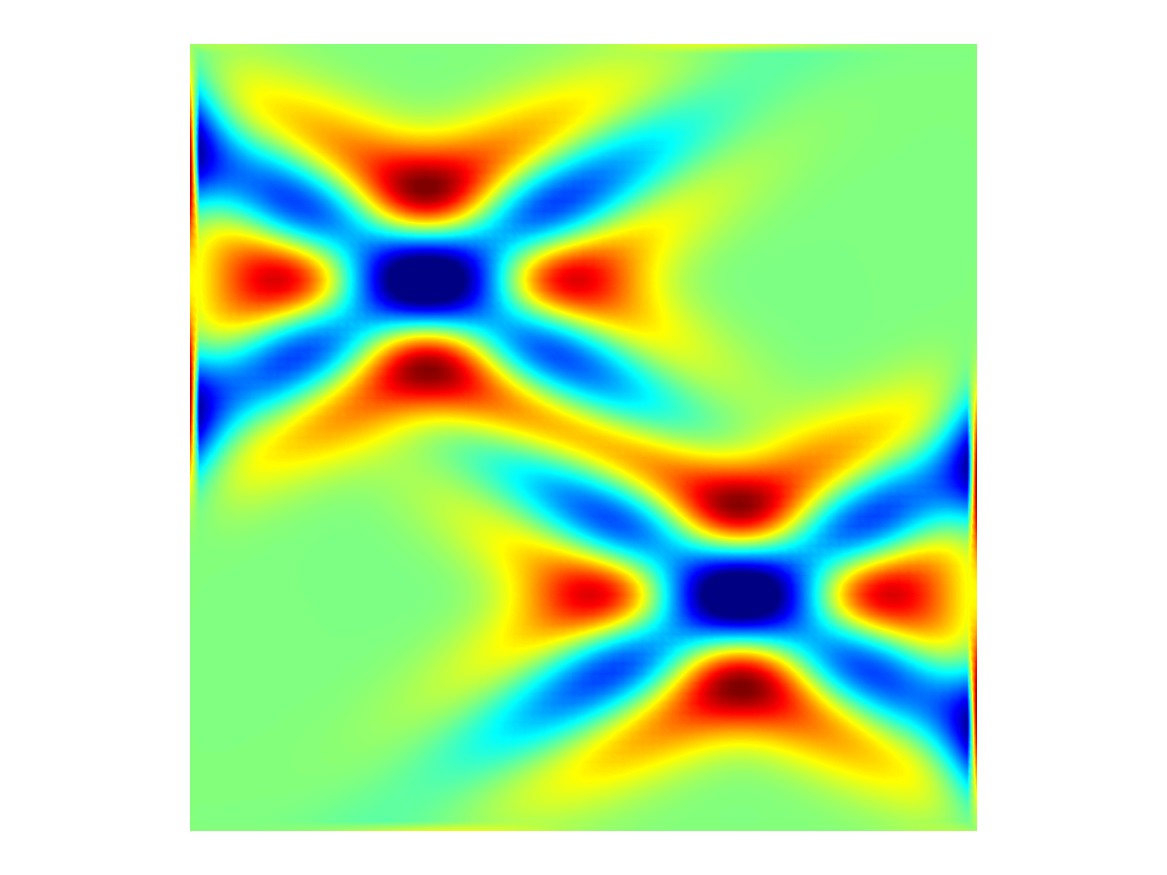}&
\includegraphics[width=.23\linewidth]{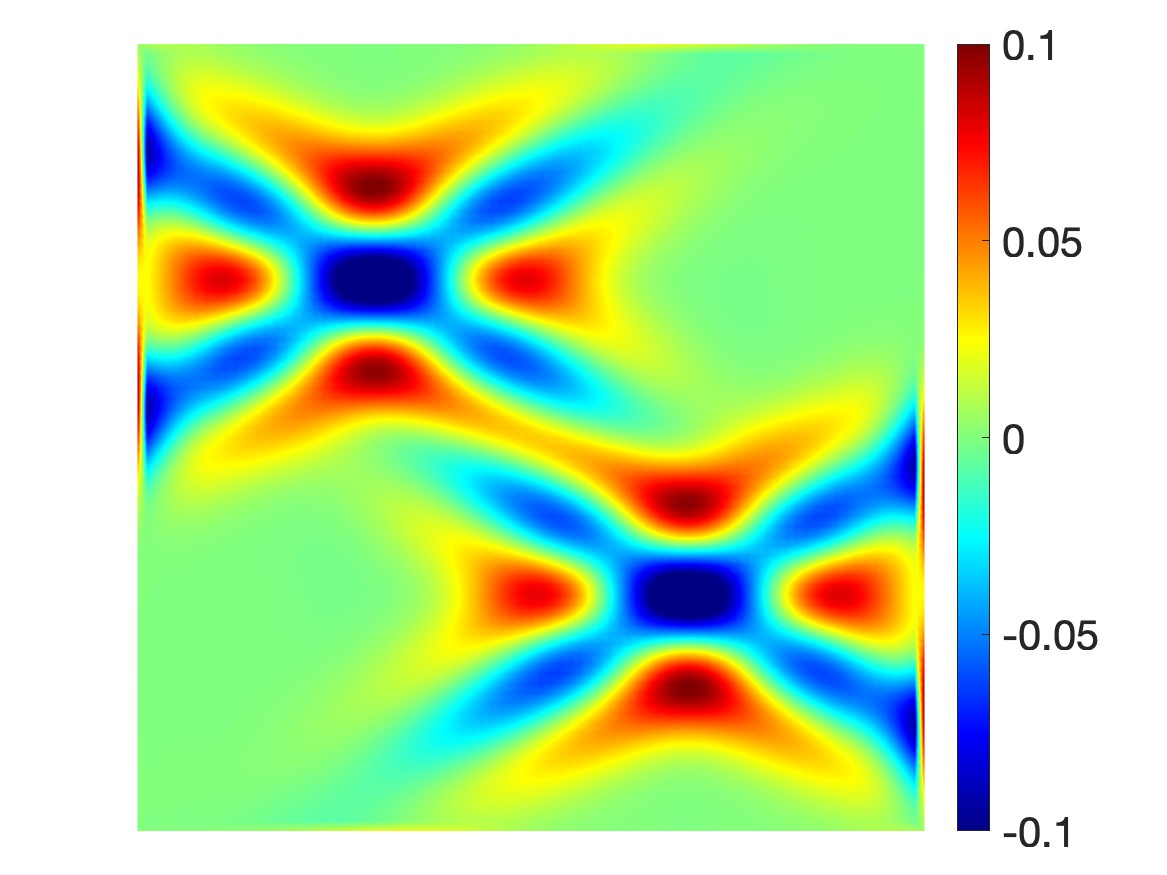}\\[0.1ex]
&Truth & Dirichlet & Robin &  Dirichlet Error&  Robin Error%
\end{tabular}
\caption{True and reconstructed absorption (top row), diffusion (middle row), and Gr\"uneisen (bottom row) coefficients in Experiment II. Shown from left to right are the true coefficients (first column), reconstructions with Dirichlet and Robin boundary conditions (middle two columns), and the relative error of the reconstructions (right two columns).
}%
\label{FIG:e2}
\end{figure}
\paragraph{Experiment II.}  We now simulate the case of reconstructing all three coefficients $\sigma$, $\gamma$, and $\Gamma$. We start with a case where the coefficients are smooth functions consisting of superpositions of Gaussian functions. To be precise, the true coefficients are 
\begin{equation}\label{EQ:oval}
\begin{aligned} 
\sigma(\bx) &=&& \exp \left(-\frac{1}{2}|R_{\pm\pi/3}\Sigma_\sigma(\bx-\bx_0^1)|^2\right)
 +\exp \left(-\frac{1}{2}|R_{\pm\pi/3}\Sigma_\sigma(\bx-\bx_0^2)|^2\right),\\
\gamma(\bx) &=&& 0.04\exp\left(-\frac{1}{2}|R_{\pi/6}\Sigma_\gamma(\bx-\bx_1)|^2\right),\\
\Gamma(\bx)&=&& \exp\left(-\frac{1}{2}|R_{\pm\pi/8}\Sigma_\Gamma(\bx-\bx_2^1)|^2\right)  +\exp\left(-\frac{1}{2}|R_{\pm\pi/8}\Sigma_\Gamma(\bx-\bx_2^2)|^2\right),\\
\end{aligned}
\end{equation}
where the centered points for these Gaussians are $\bx_0^1=(0.5, 1.5), \bx_0^2=(1.5, 0.5), \bx_1=(1.0, 1.0), \RED{\bx_2^1=(0.5, 0.5), \bx_2^2=(1.5, 1.5)}$.
$R_\theta=\begin{bmatrix}
    \cos\theta & -\sin\theta\\
    \sin\theta & \cos\theta
\end{bmatrix}$ is the rotation matrix, $\Sigma_\sigma=\Sigma_\Gamma=\begin{bmatrix}
    0.25 & 0 \\ 0 & 4
\end{bmatrix}$ and $\Sigma_\gamma=\begin{bmatrix}
    0.01 & 0 \\ 0 & 0.1
\end{bmatrix}$ are the variance matrix for $\sigma,\gamma,\Gamma$ accordingly as $\Sigma=\begin{bmatrix}
    \frac{1}{\sigma_x^2} & 0 \\ 0 & \frac{1}{\sigma_y^2}
\end{bmatrix}$.

The true coefficients, the reconstructions under two different boundary conditions for the diffusion model, as well as the relative error in the reconstructions, are shown in~\Cref{FIG:e2}. The reconstructions are performed with noise-free synthetic data. Simple visual observation shows that the quality of reconstructions for such smooth coefficients is very high. In particular, while the qualify of reconstructing $(\sigma,\Gamma)$ is similar in Dirichlet and Robin boundary conditions, that for $\gamma$ is better under the Dirichlet condition for the diffusion model. The relative errors for $\sigma$ and $\gamma$ are relatively small with values at most $3\%-6\%$, while the maximum of the relative errors of $\Gamma$ reaches $10\%$. 

\begin{figure}[!htb]
\settoheight{\tempdima}{\includegraphics[width=.25\linewidth]{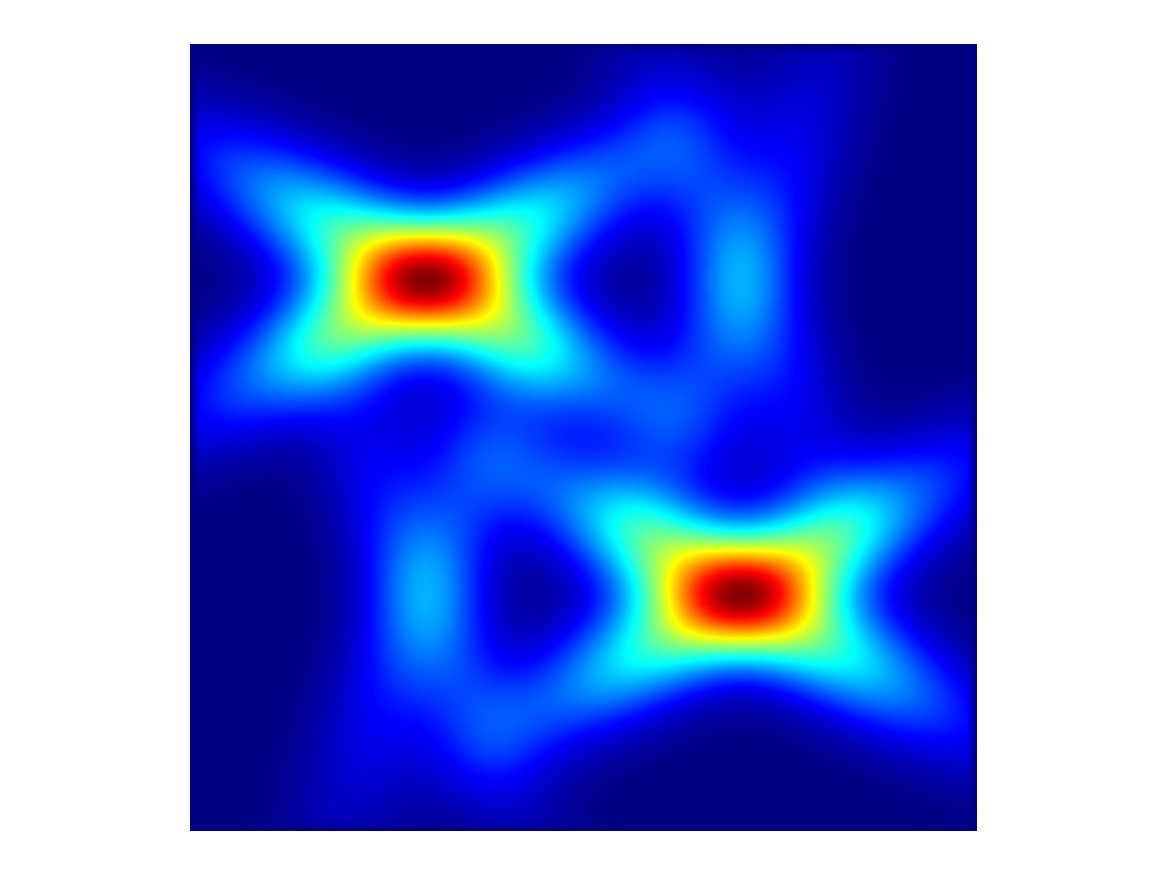}}%
\centering
\setlength{\tabcolsep}{1pt} %
\renewcommand{\arraystretch}{0.8} %
\begin{tabular}{@{}c@{}c@{}c@{}c@{}}
\rownamerev{Dirichlet}&
\includegraphics[width=.25\linewidth]{FinalFigures/GaussianMixture/crossTrue.jpg}&
\includegraphics[width=.25\linewidth]{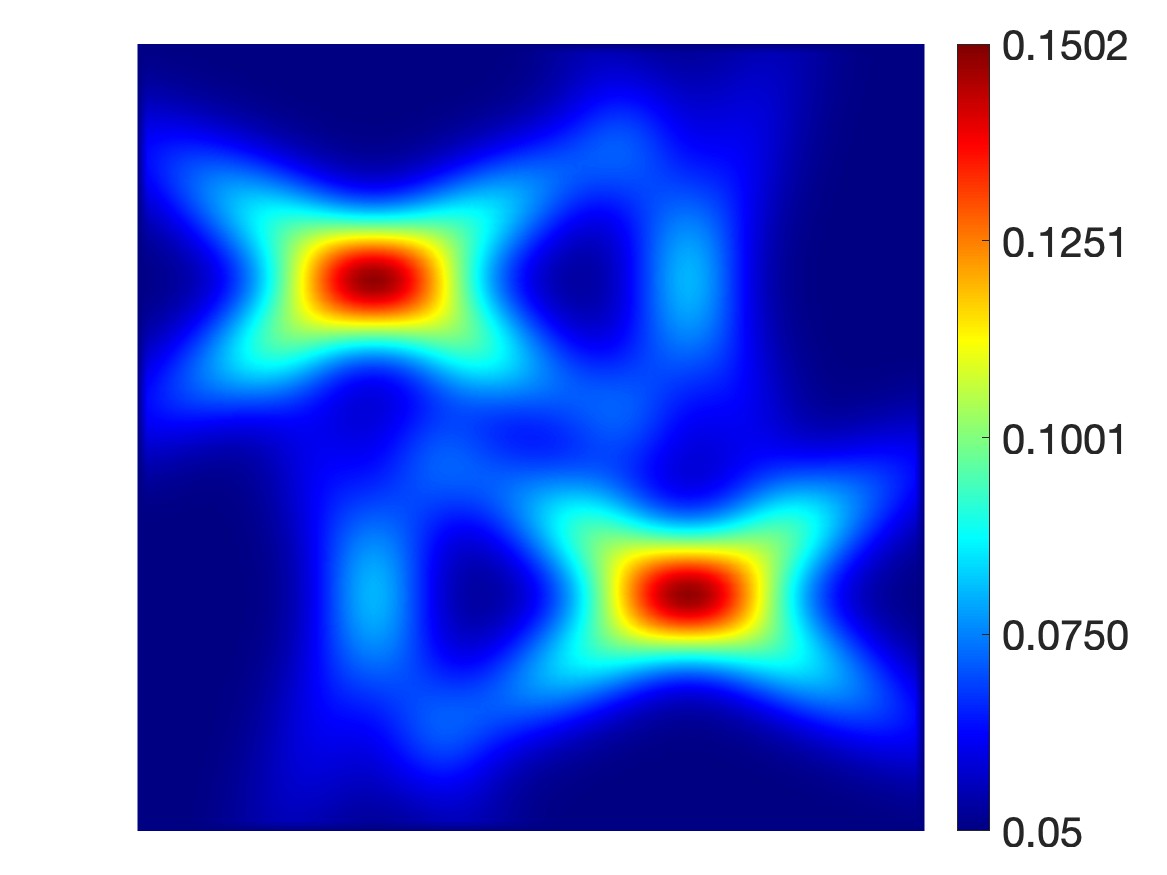}&
\includegraphics[width=.25\linewidth]{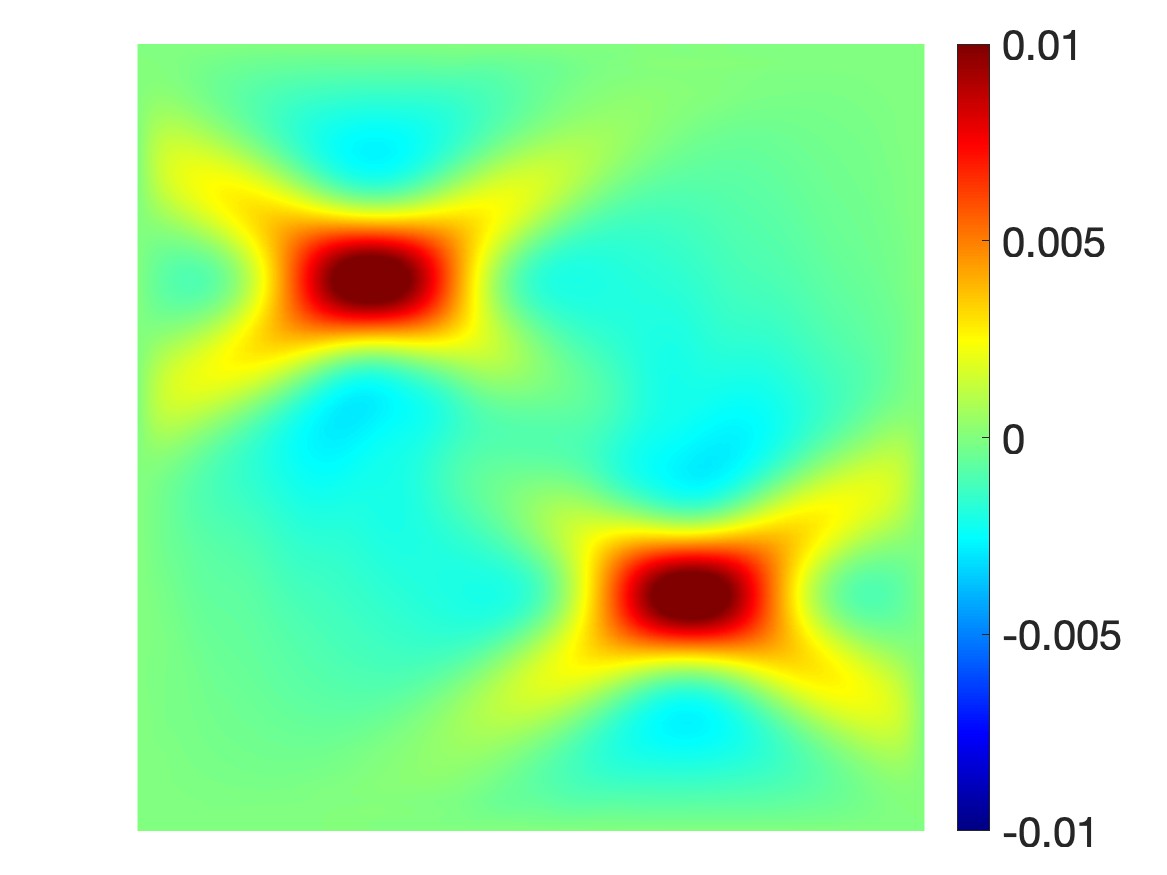}\\[-1ex]
\rownamerev{Robin} & 
\includegraphics[width=.25\linewidth]{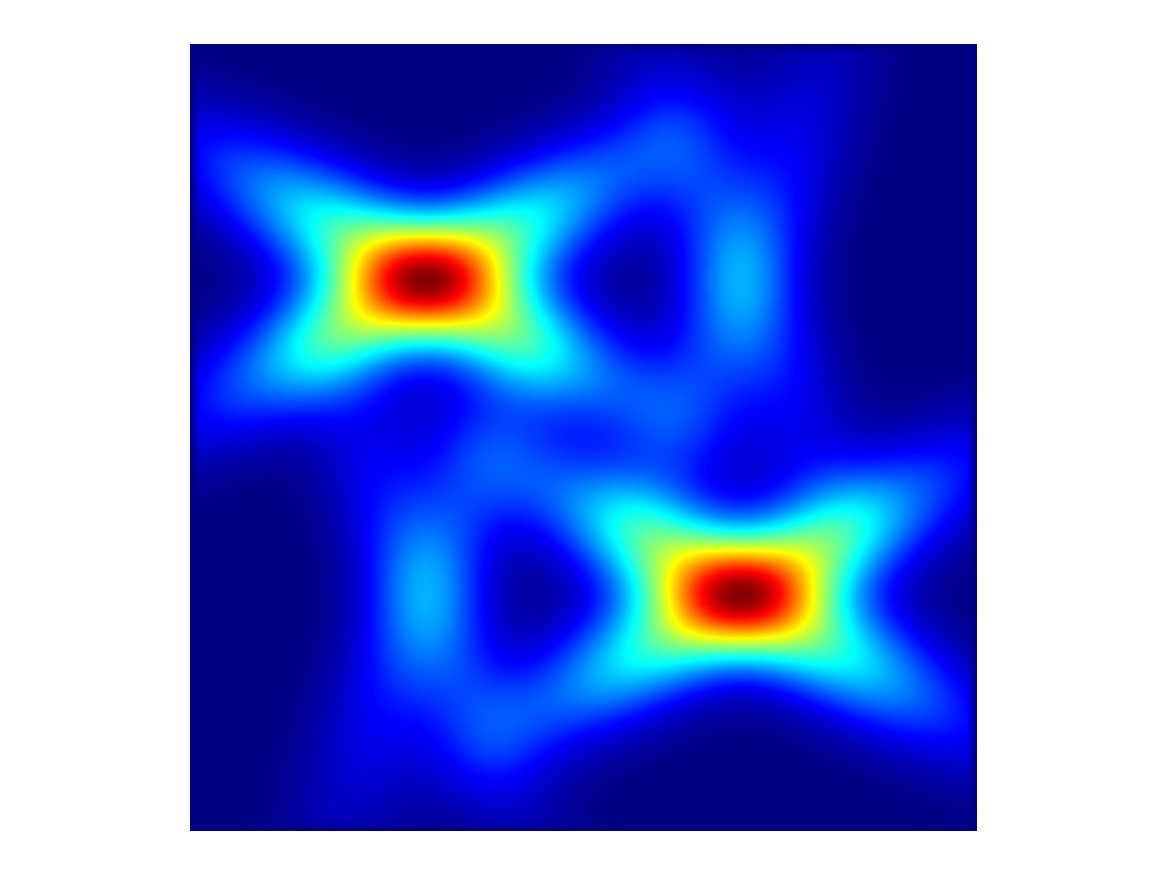}&
\includegraphics[width=.25\linewidth]{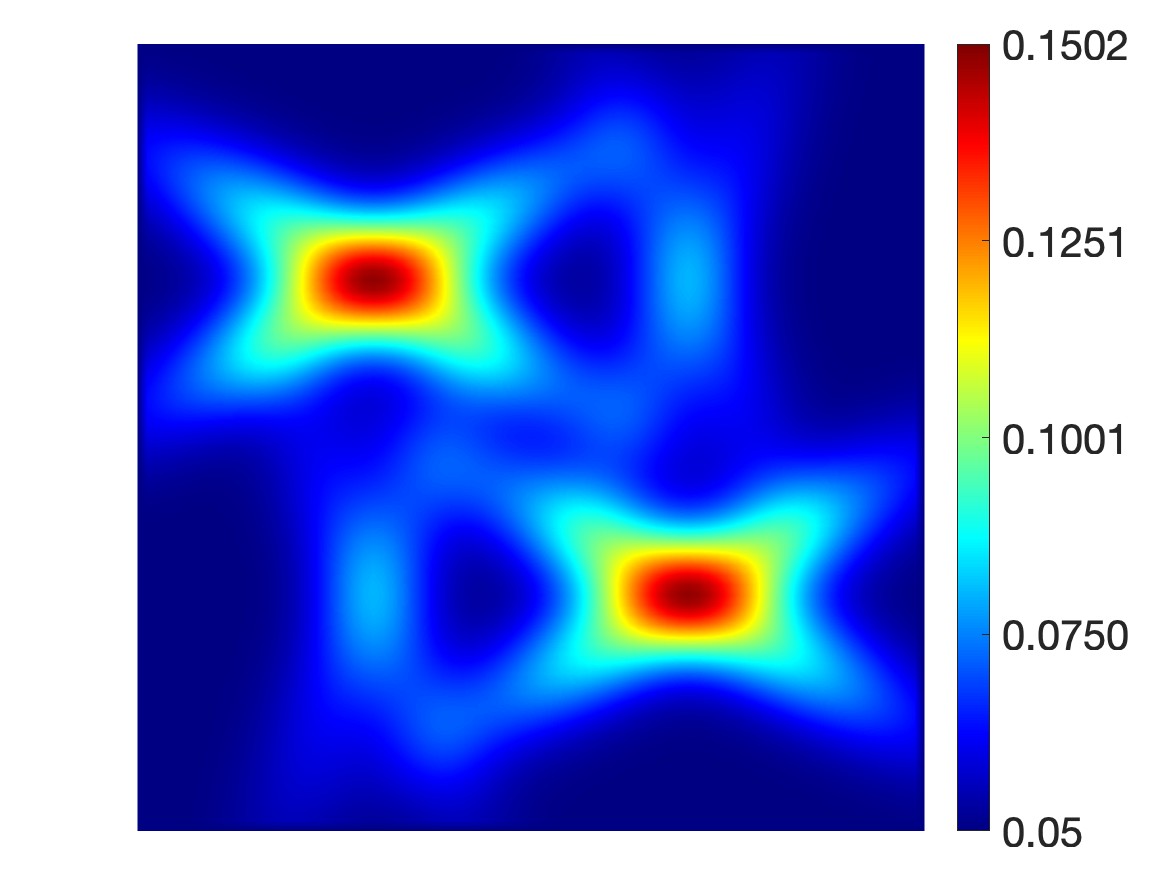}&
\includegraphics[width=.25\linewidth]{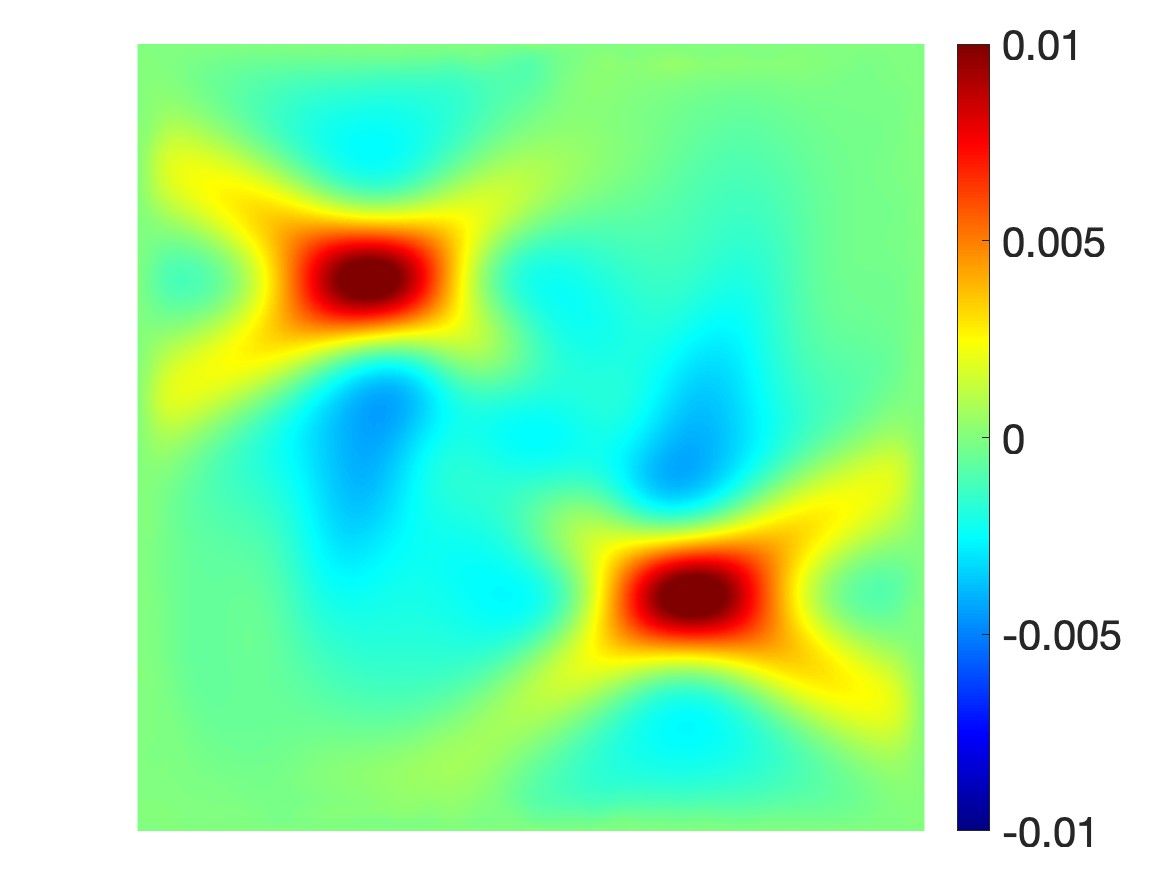}\\& True $\Gamma \sigma$ & Reconstructed $\Gamma\sigma$ & Relative error\\[-1ex]
\end{tabular}
\caption{Reconstruction of the product $\mu:=\Gamma\sigma$ in Experiment II with Dirichlet (top) and Robin (bottom) boundary conditions. Shown are the true $\Gamma\sigma$ (first column), the reconstructions (second column), and relative errors (third column) in the reconstructions.
}
\label{FIG:e2-cross}
\end{figure}
\RED{We again observe clearly cross-talk effect in the reconstruction, mainly between $\sigma$ and $\Gamma$ this time. The error plots in the right two columns of~\Cref{FIG:e2} show that large errors of $\sigma$ reconstructions happen at the places where the $\Gamma$ coefficient is large and vice versa. To further highlight the cross-talk effect, we show in~\Cref{FIG:e2-cross} the reconstruction and the relative error of the product of the absorption coefficient and the Gr\"uneisen coefficient: $\Gamma\sigma$. It is clear that $\Gamma \sigma$ is reconstructed fairly well, with a relative error around $1\%$, better than the individual reconstructions of either $\gamma$ or $\sigma$ as shown in~\Cref{FIG:e2}.}

\begin{figure}[htb!]
\settoheight{\tempdima}{\includegraphics[width=.3\linewidth]{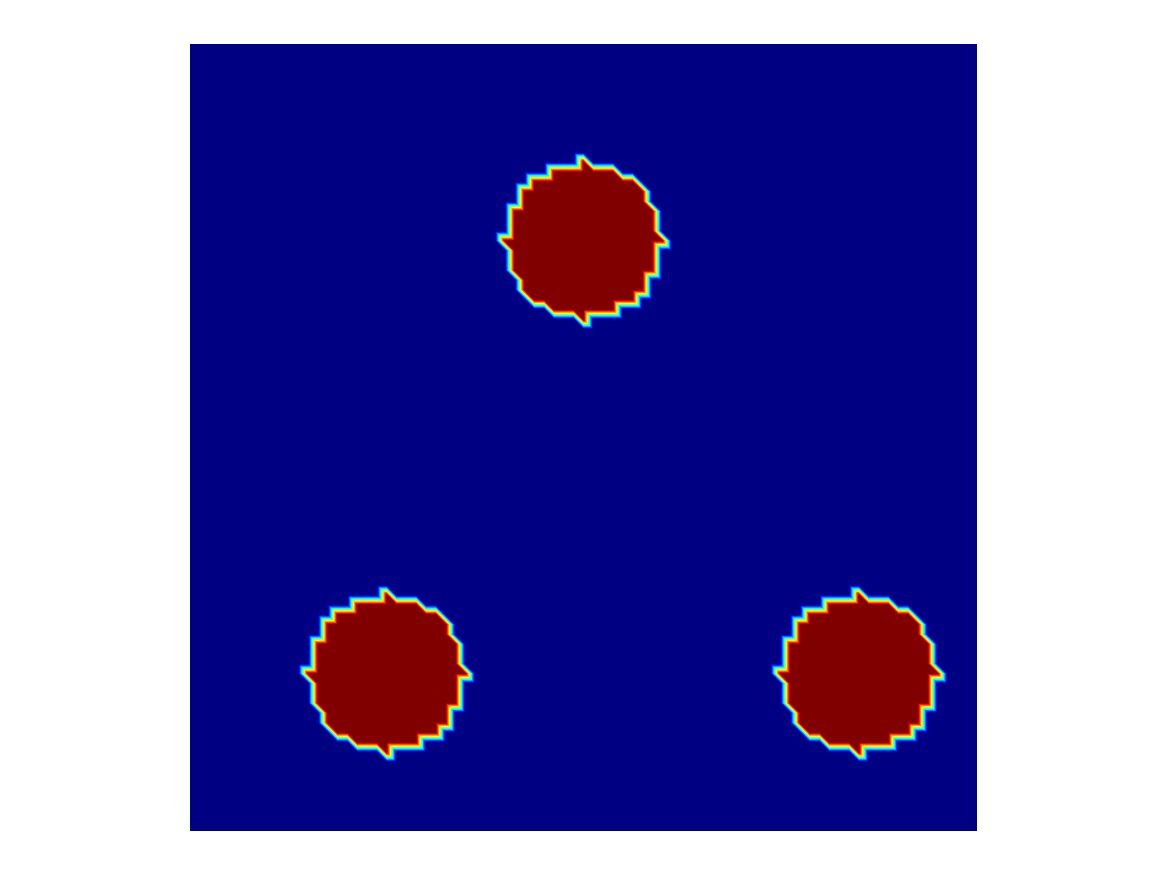}}%
\centering
\begin{tabular}
{@{\hspace{-.3ex}}c@{\hspace{-0ex}}c@{\hspace{-3.5ex}}c@{\hspace{-1ex}}c@{}}
\rowname{$\sigma$}&
\includegraphics[width=.3\linewidth]{FinalFigures/Piecewise/sigmat.jpg}&
\includegraphics[width=.3\linewidth]{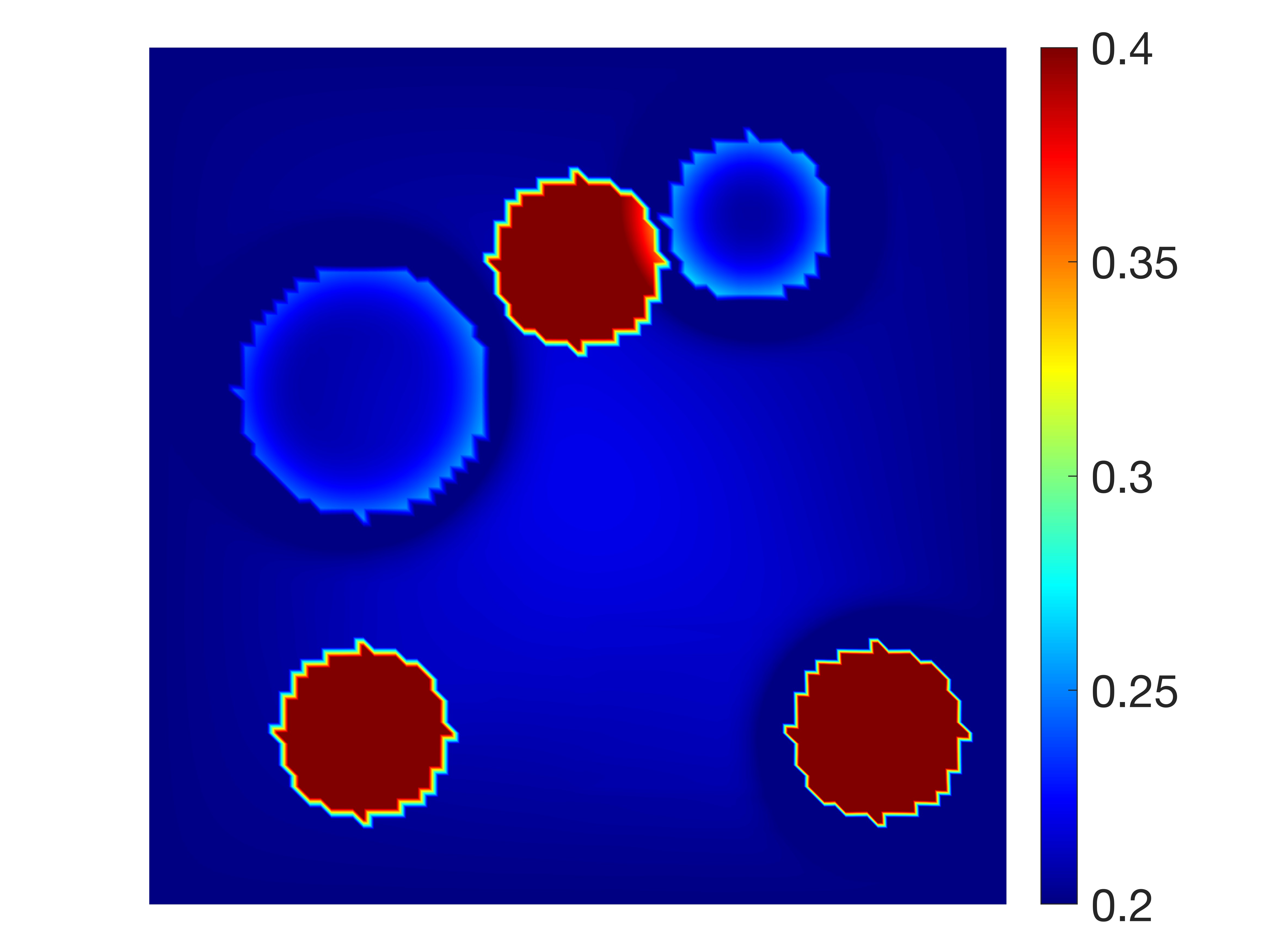}&
\includegraphics[width=.3\linewidth]{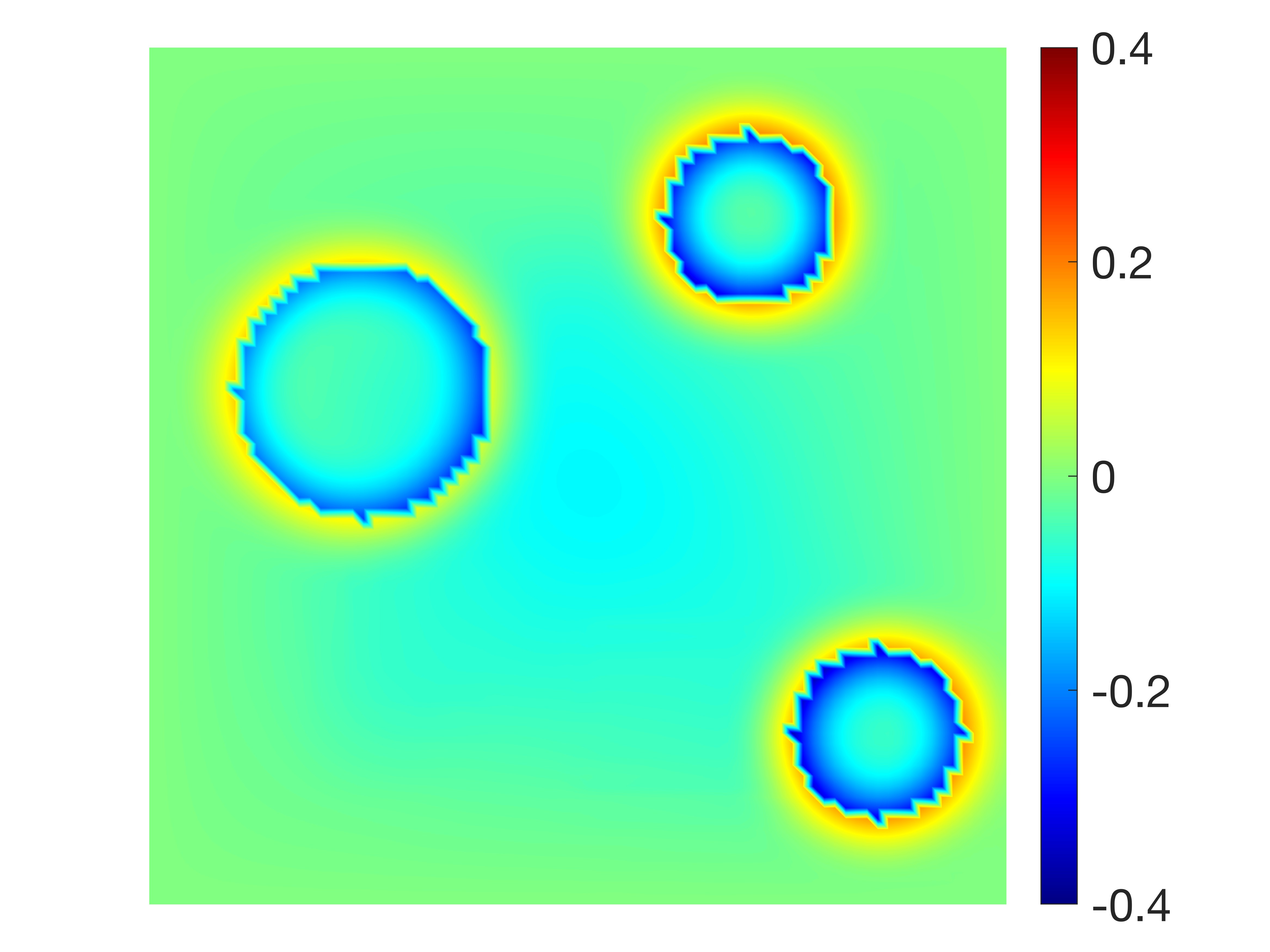}\\
\rowname{$\gamma$}&
\includegraphics[width=.3\linewidth]{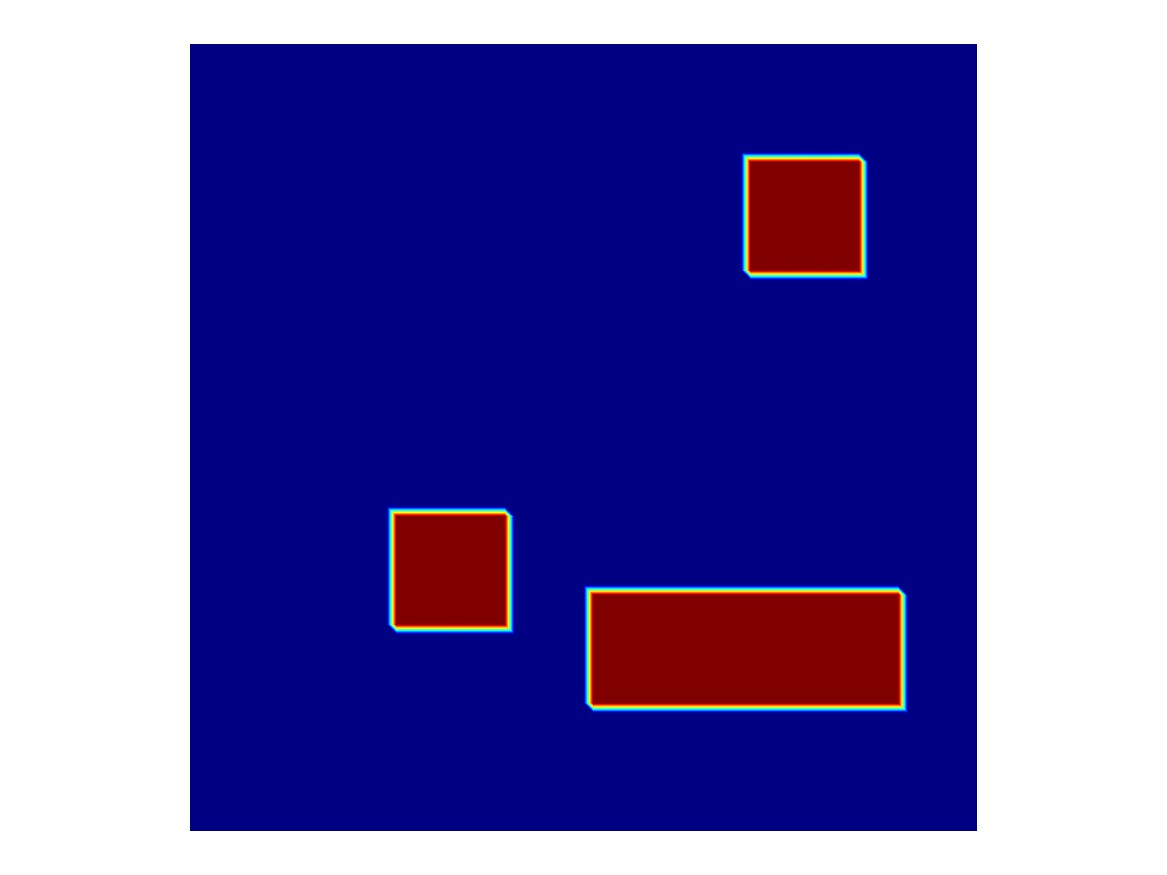}&
\includegraphics[width=.3\linewidth]{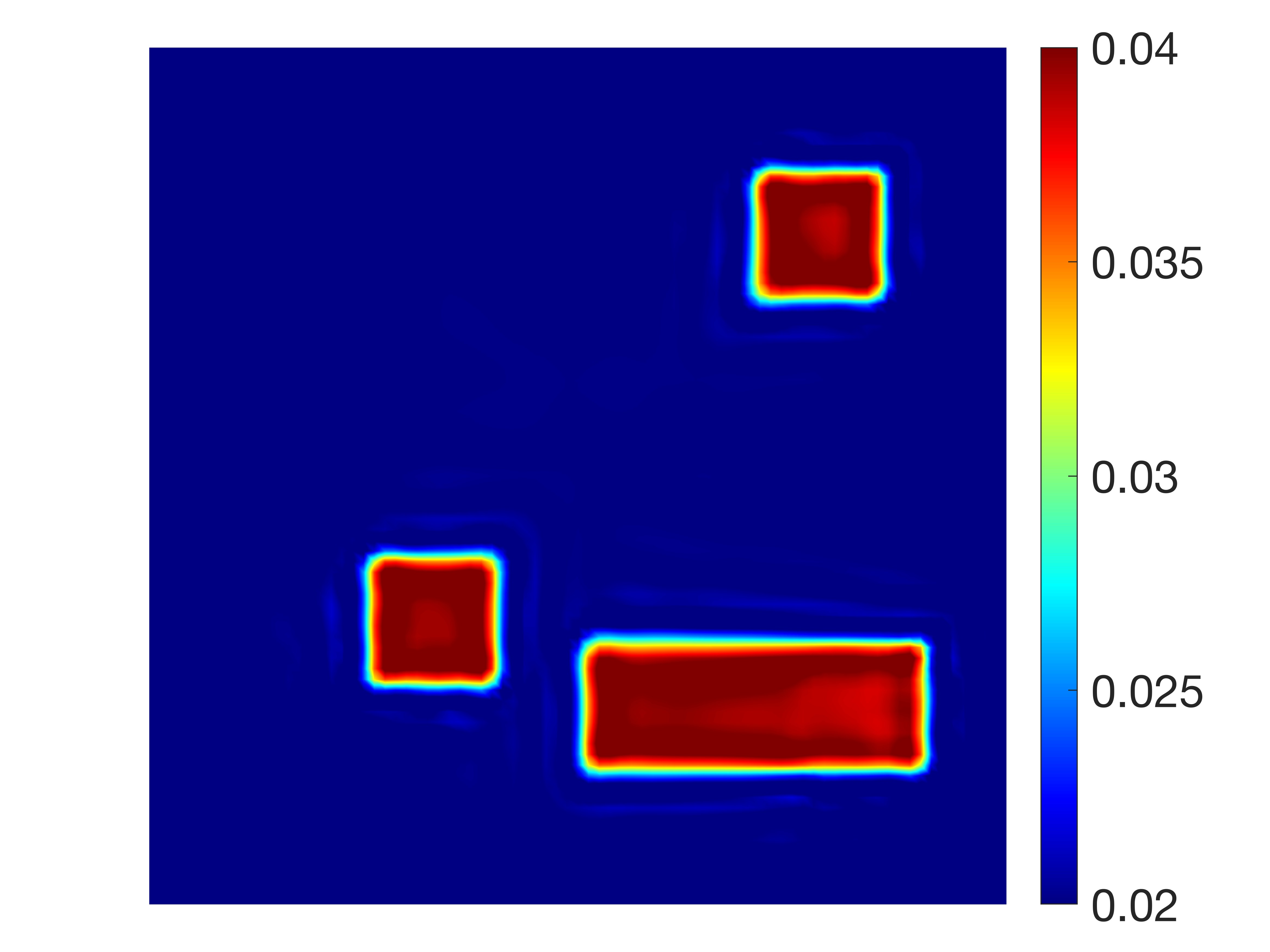}&
\includegraphics[width=.3\linewidth]{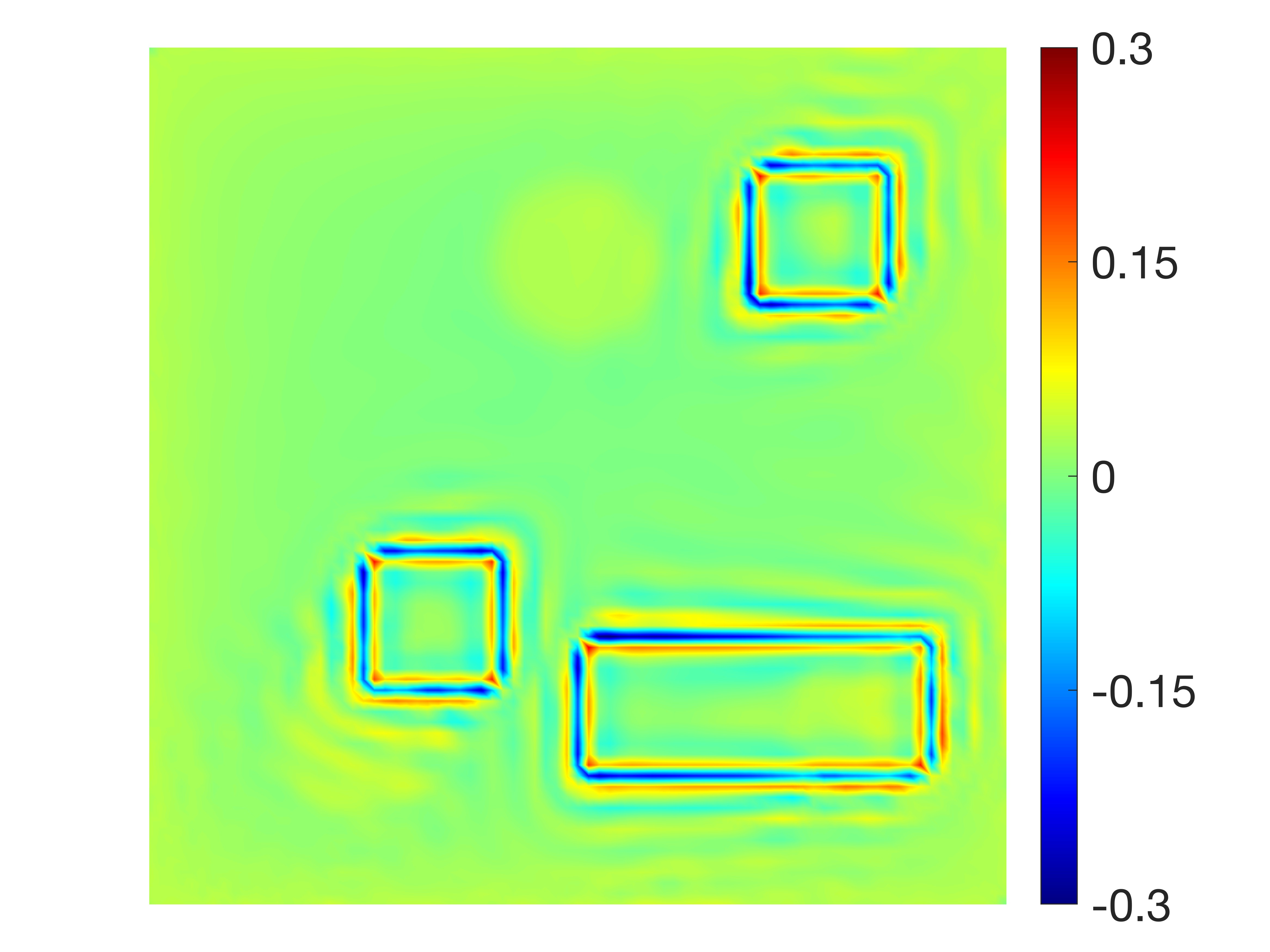}\\
\rowname{$\Gamma$} &
\includegraphics[width=.3\linewidth]{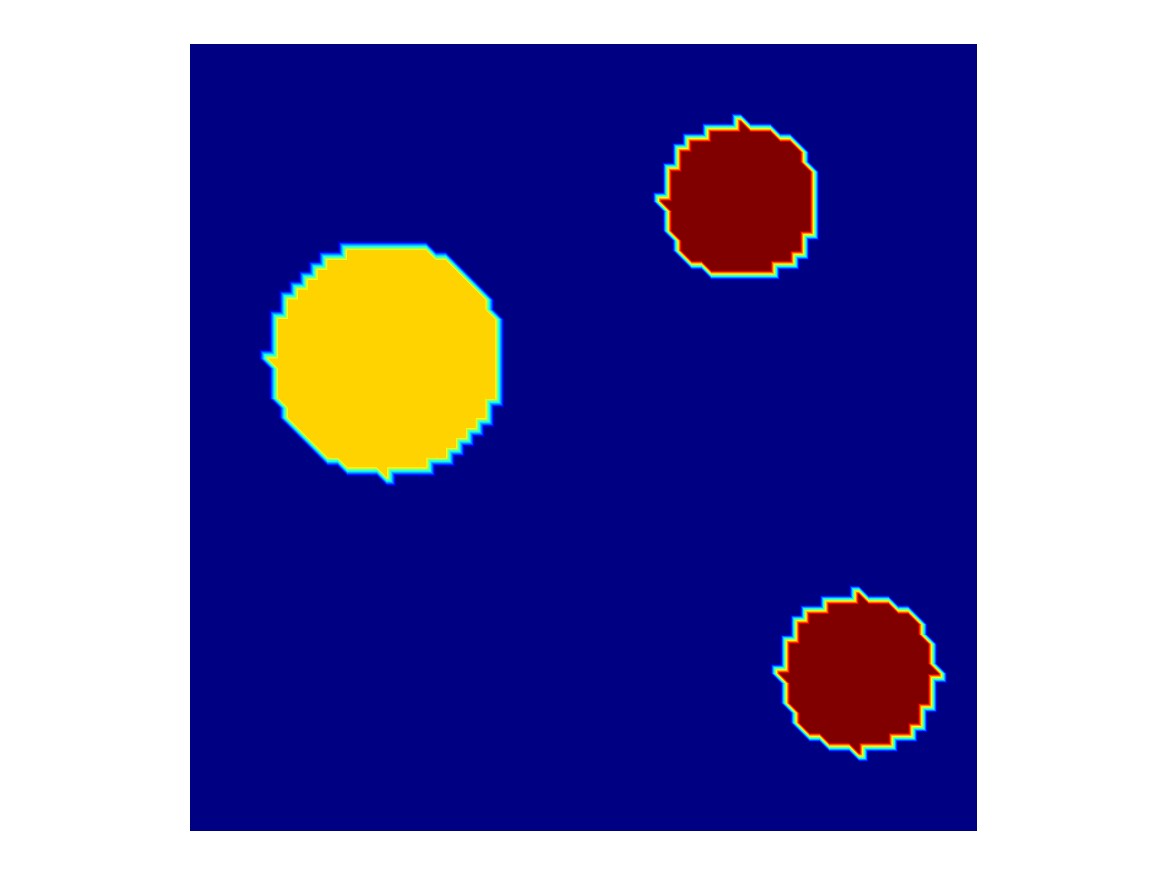}&
\includegraphics[width=.3\linewidth]{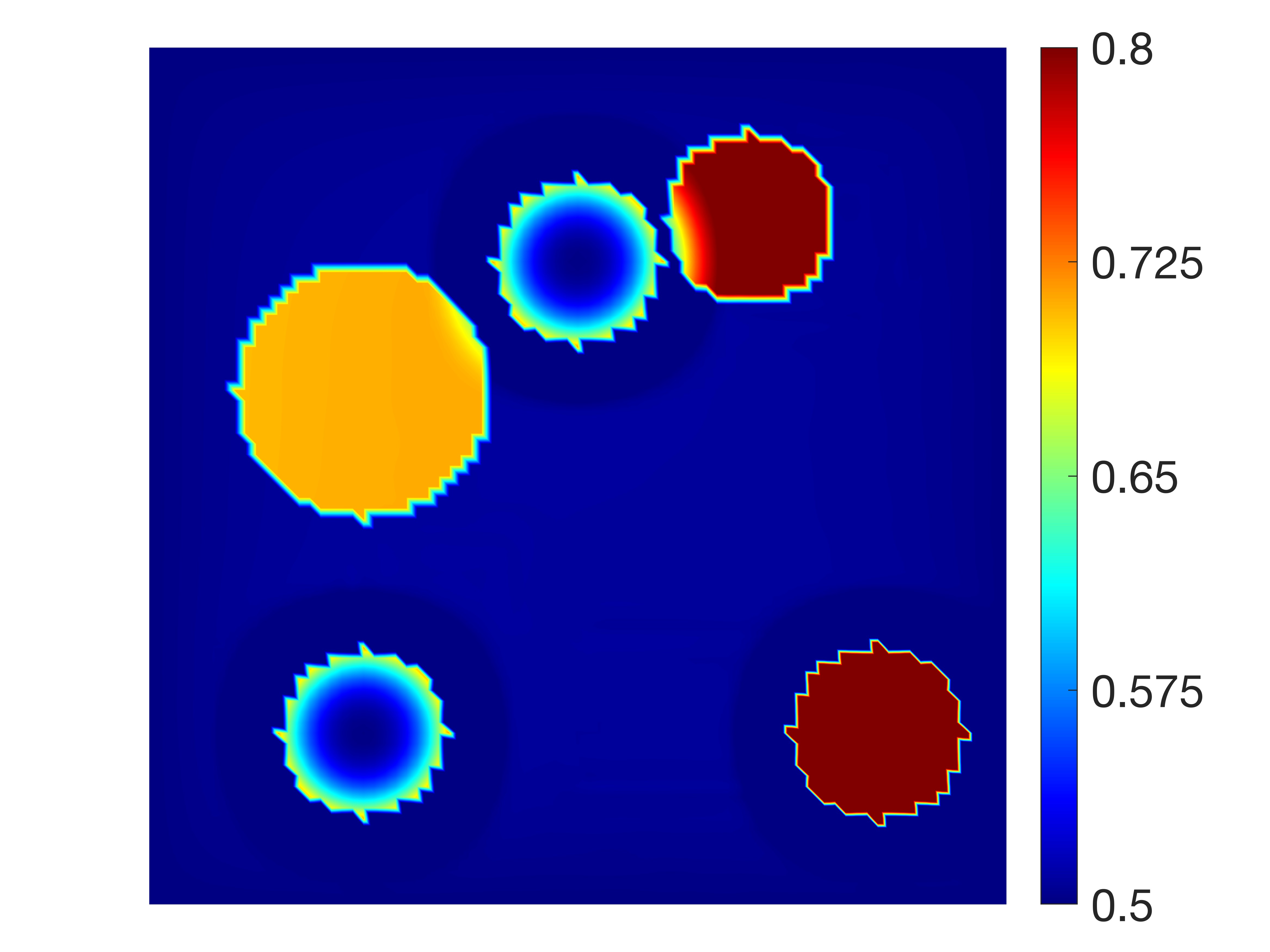}&
\includegraphics[width=.3\linewidth]{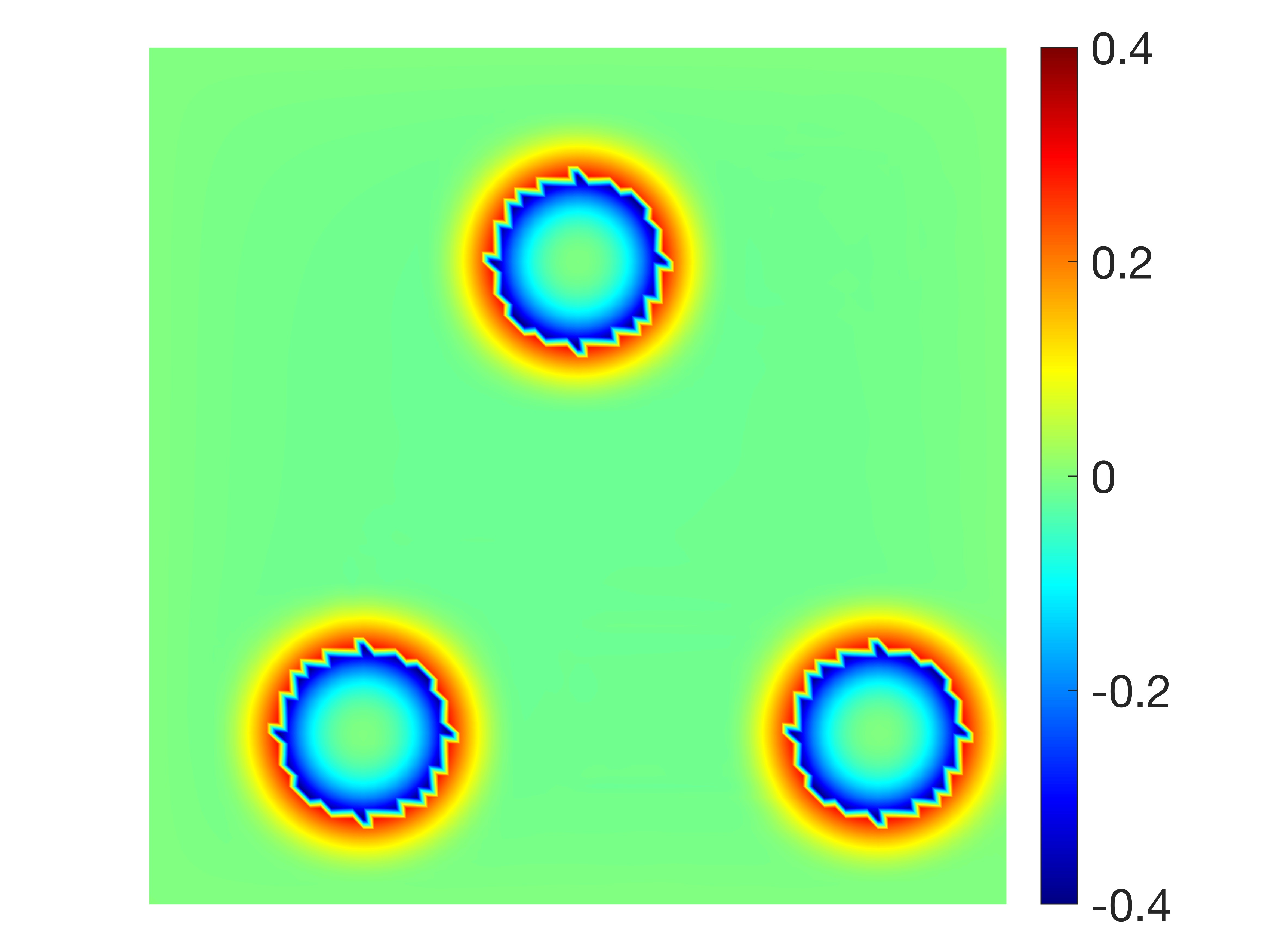}\\[0.1ex]
&True coefficients & Reconstructions &  Relative error%
\end{tabular}
\caption{True coefficients (left column), reconstructions from noise-free data (middle column), and the relative error (right column) in Experiment III.}%
\label{FIG:e3 noise-free}
\end{figure}
\paragraph{Experiment III.} We now reproduce the simulations in the previous numerical experiment with a case where the coefficients have discontinuities. The true coefficients are
\begin{equation}
    \begin{aligned}
    \sigma(\bx) &=&&0.2+0.2 \chi_{B_{0.2}((1,1.5))} +0.2 \chi_{B_{0.2}((1.7,0.4))}+0.2 \chi_{B_{0.2}((0.5,0.4))}\\
    \gamma(\bx) &=&&0.02 + 0.02\chi_{{[0.5,0.8]\times[0.5,0.8]}} + 0.02 \chi_{{[1.4,1.7]\times[1.4,1.7]}} + 0.02 \chi_{{[1,1.8]\times[0.3,0.6]}}\\
    \Gamma(\bx) &=&&0.5+0.2 \chi_{B_{0.3}((0.5,1.2))} +0.3 \chi_{B_{0.2}((1.4,1.6)))}+0.3 \chi_{B_{0.2}((1.7,0.4))}
    \end{aligned}
\end{equation}
where $\chi_A$ is the indicator function for set $A$, and $B_r(\bx)$ is a ball centered at $\bx$ with radius $r$. 

\begin{figure}[htb!]
\settoheight{\tempdima}{\includegraphics[width=.3\linewidth]{FinalFigures/Piecewise/sigmat.jpg}}%
\centering
\begin{tabular}
{@{\hspace{-.3ex}}c@{\hspace{-0ex}}c@{\hspace{-3.5ex}}c@{\hspace{-1ex}}c@{\hspace{-1ex}}c@{\hspace{-2.7ex}}c@{}}
\rowname{$\sigma$}&
\includegraphics[width=.3\linewidth]{FinalFigures/Piecewise/sigmat.jpg}&
\includegraphics[width=.3\linewidth]{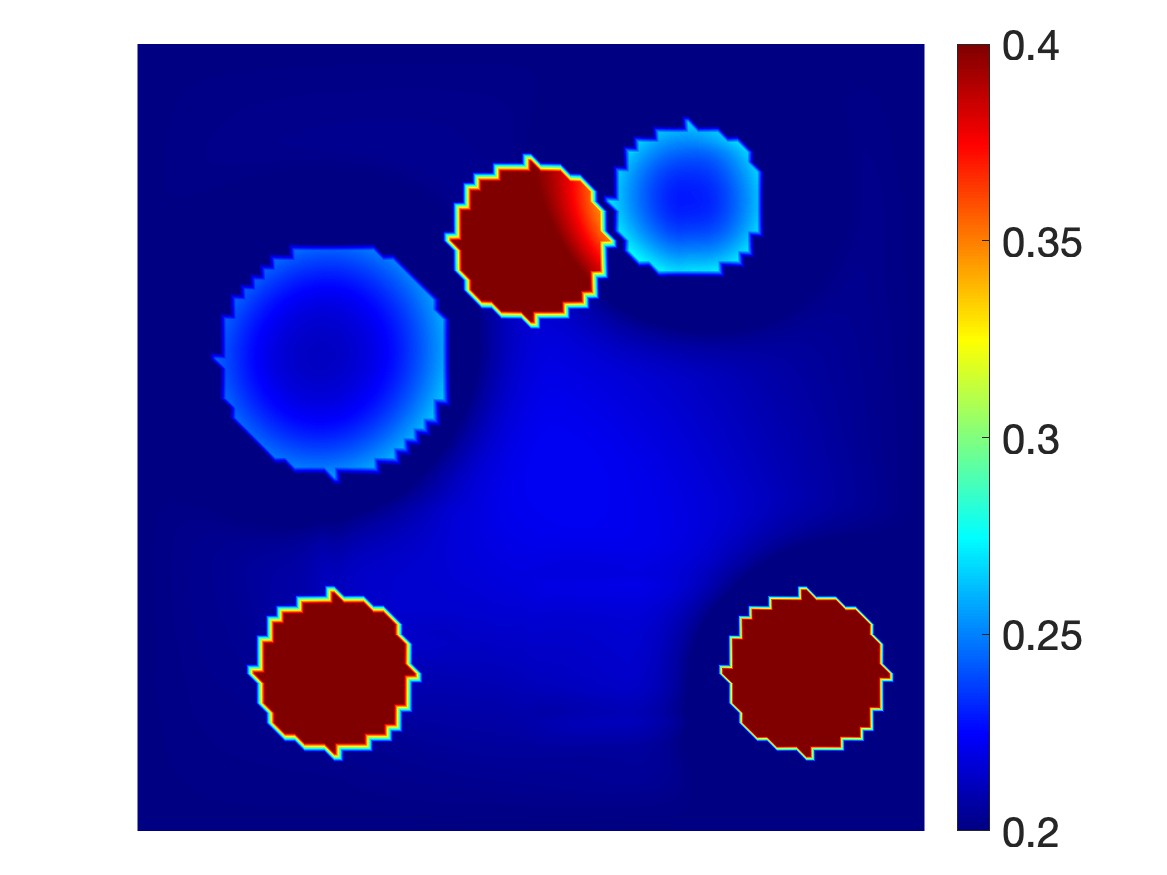}&
\includegraphics[width=.3\linewidth]{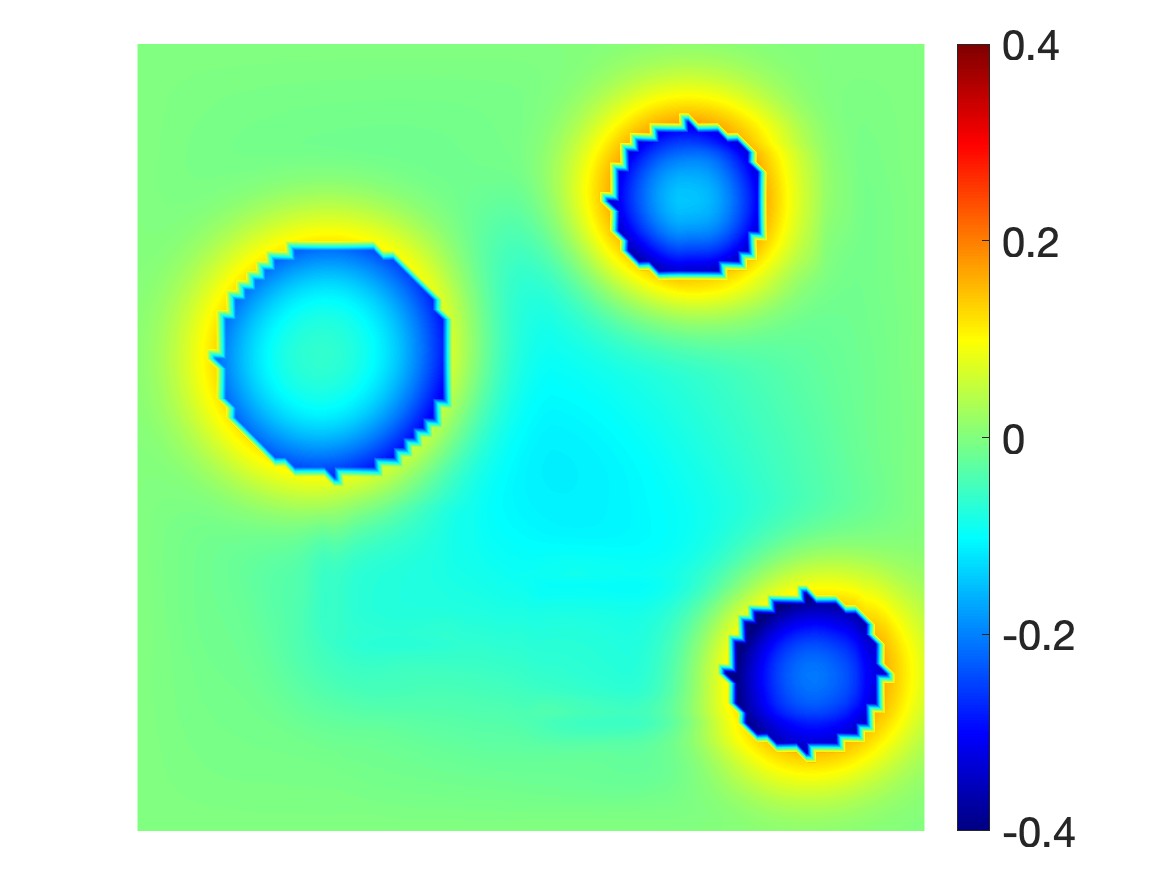}\\
\rowname{$\gamma$}&
\includegraphics[width=.3\linewidth]{FinalFigures/Piecewise/Difft.jpg}&
\includegraphics[width=.3\linewidth]{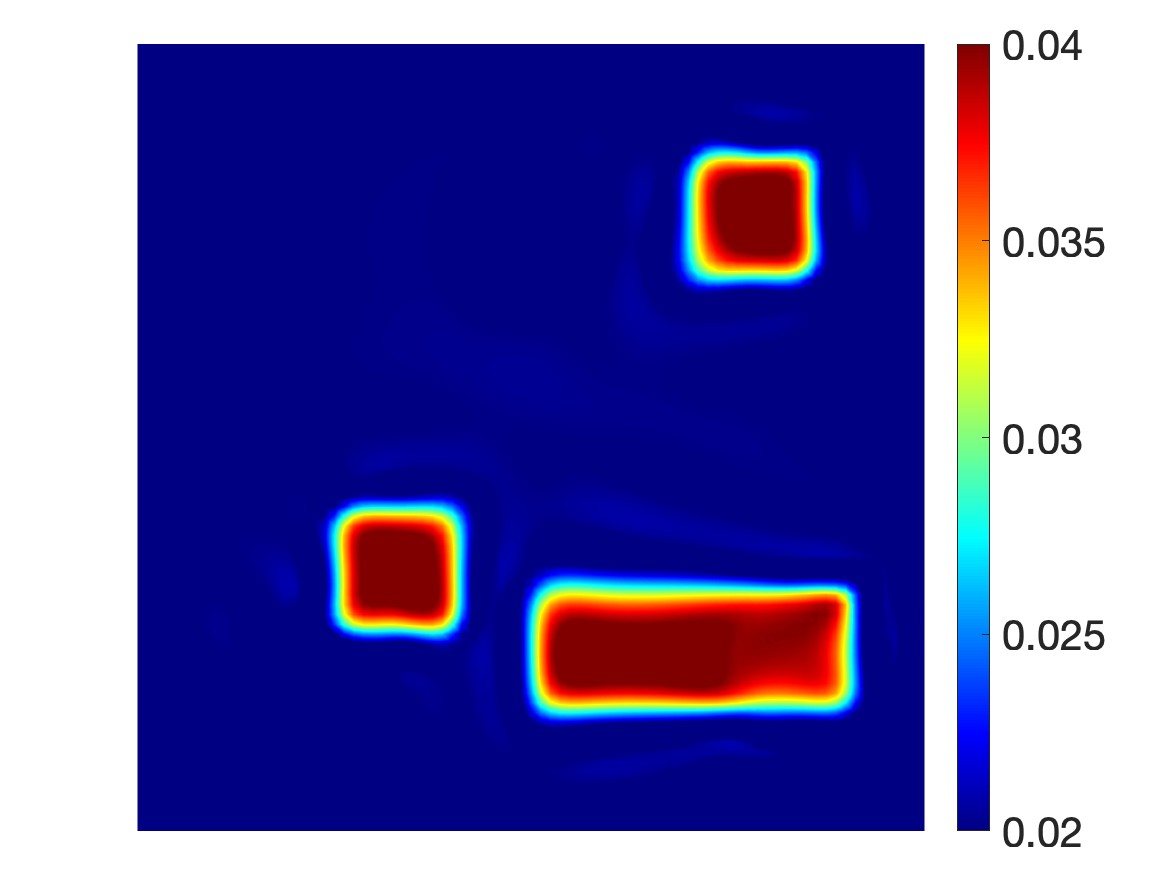}&
\includegraphics[width=.3\linewidth]{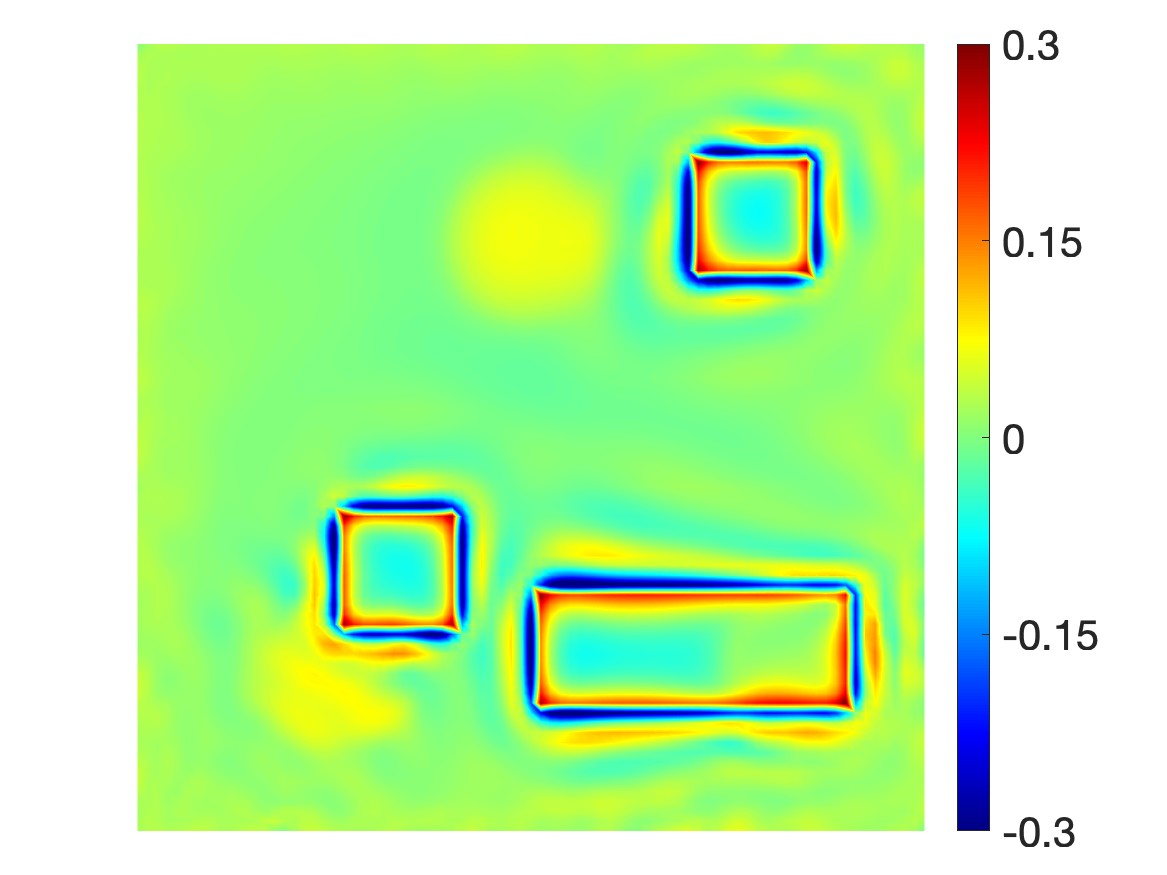}\\
\rowname{$\Gamma$} &
\includegraphics[width=.3\linewidth]{FinalFigures/Piecewise/Gammat.jpg}&
\includegraphics[width=.3\linewidth]{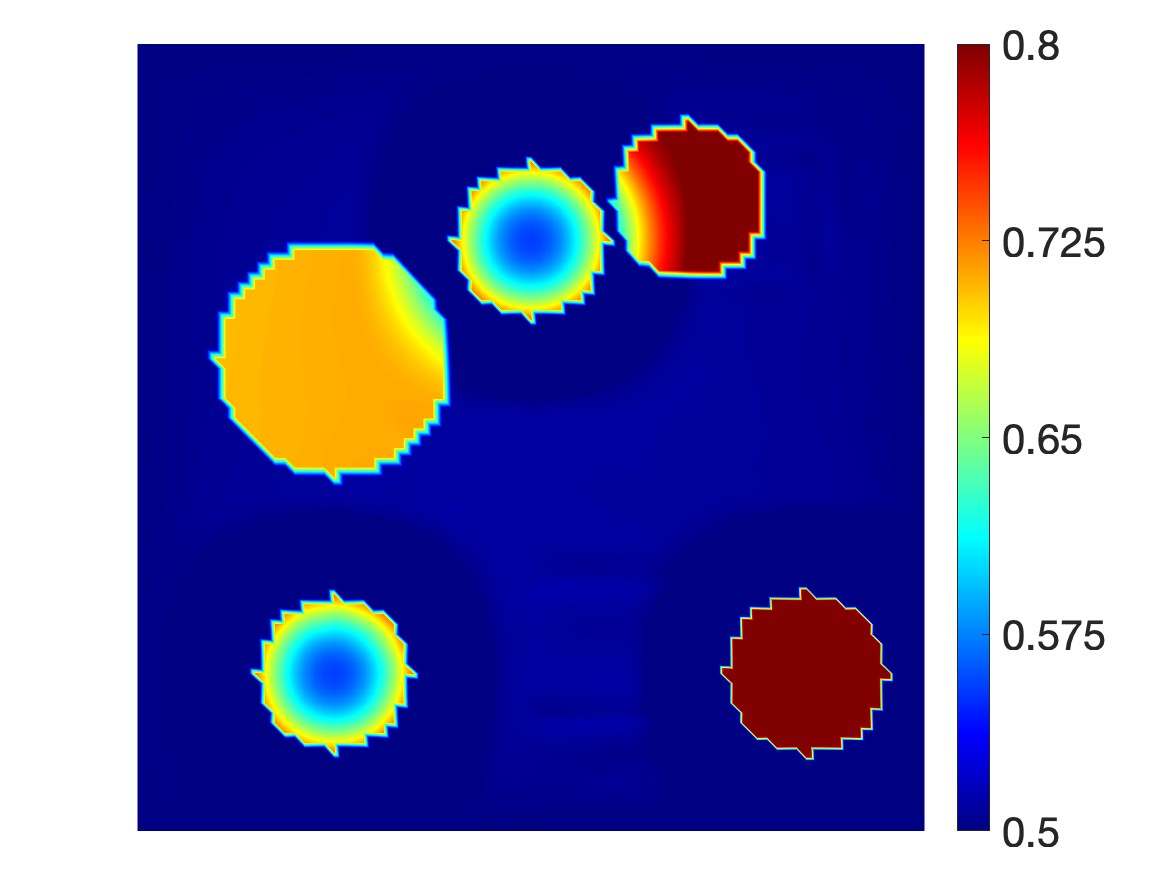}&
\includegraphics[width=.3\linewidth]{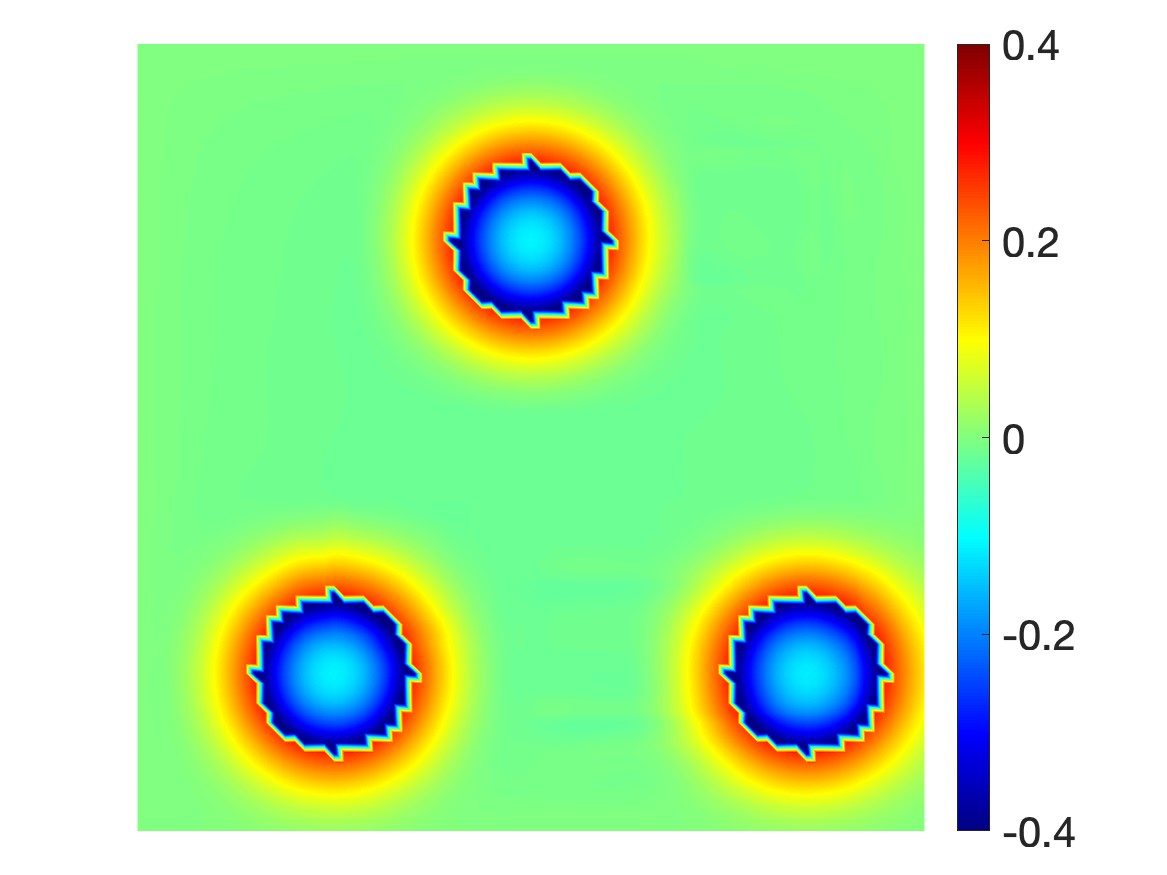}\\[0.1ex]
&True coefficients & Reconstructions &  Relative error
\end{tabular}
\caption{True coefficients (left), reconstructions from noisy data (middle column), and the relative error (right column) in Experiment III.}%
\label{FIG:e3 noisy}
\end{figure}

Reconstructions are performed with both noisy and noise-free data under the Dirichlet boundary condition. Reconstruction results with noise-free data and the relative errors are shown in~\Cref{FIG:e3 noise-free} while the same simulations with noisy data containing $5\%$ multiplicative noise are shown in~\Cref{FIG:e3 noisy}. We observe that our~\Cref{ALG:Three-Stage}
successfully recovers the basic shapes of these coefficients for both the noise-free and noisy data cases. As expected, the results with noisy data have slightly worse reconstructions in that the boundaries of these shapes are less clear than the ones with noise-free data.
 The relative errors of these constructions for coefficients with discontinuous values are larger than those in Experiment II, where smooth coefficients were pursued. This probably happens due to the Tikhonov regularization~\eqref{reg} we used in the objective function, which penalizes the nonsmooth coefficients more than the smooth ones. 
 
\begin{figure}[!htb]
\settoheight{\tempdima}{\includegraphics[width=.25\linewidth]{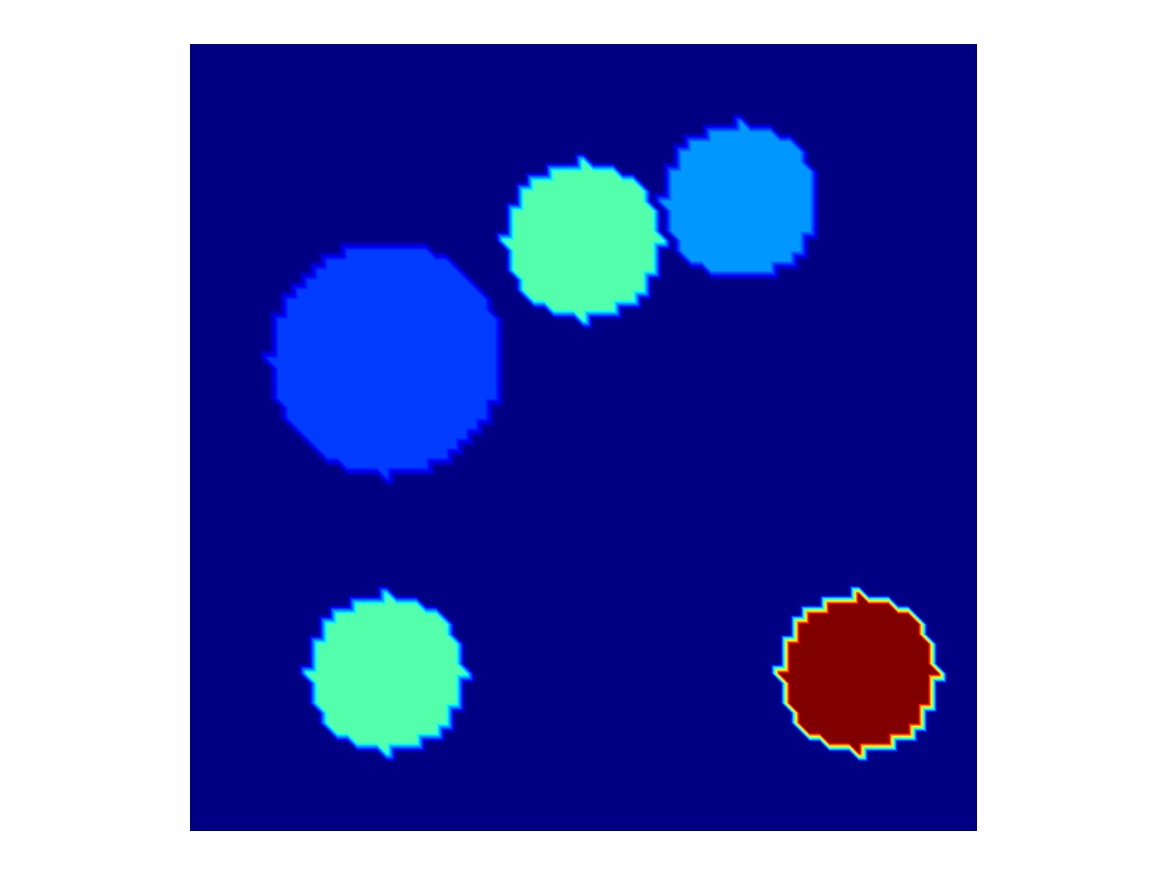}}%
\centering\begin{tabular}{@{}c@{}c@{}c@{}c@{}}
\rownamerev{Noise-free} &
\includegraphics[width=.25\linewidth]{FinalFigures/Piecewise/crossTrue.jpg}&
\includegraphics[width=.25\linewidth]{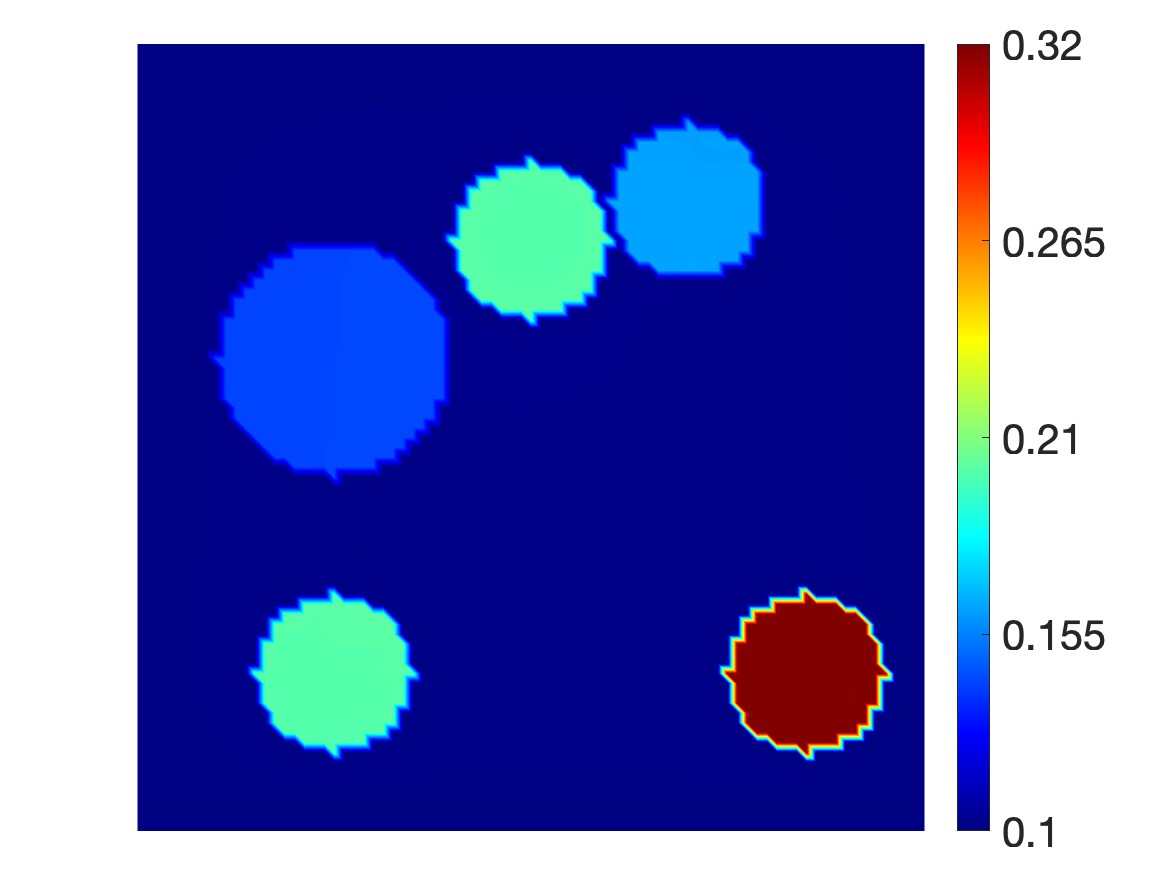}&
\includegraphics[width=.25\linewidth]{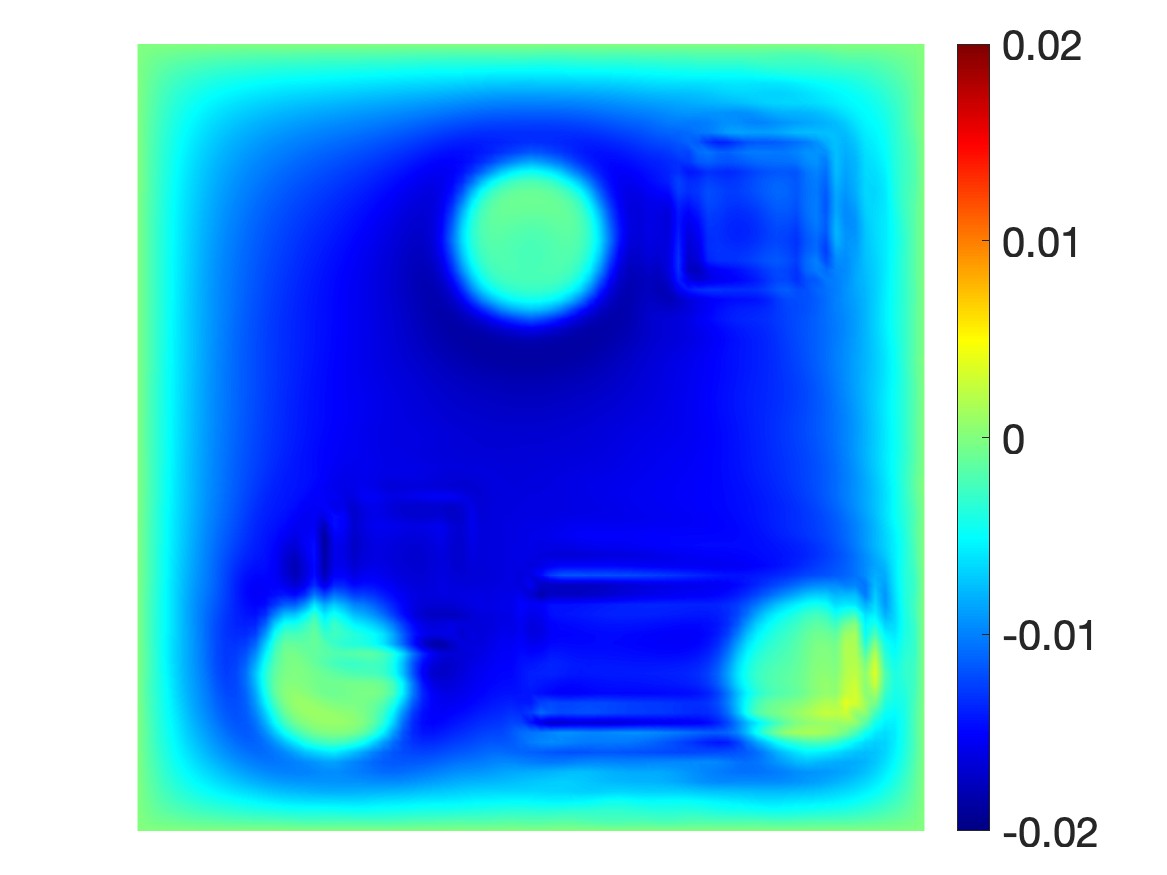}\\
\rownamerev{Noisy} &
\includegraphics[width=.25\linewidth]{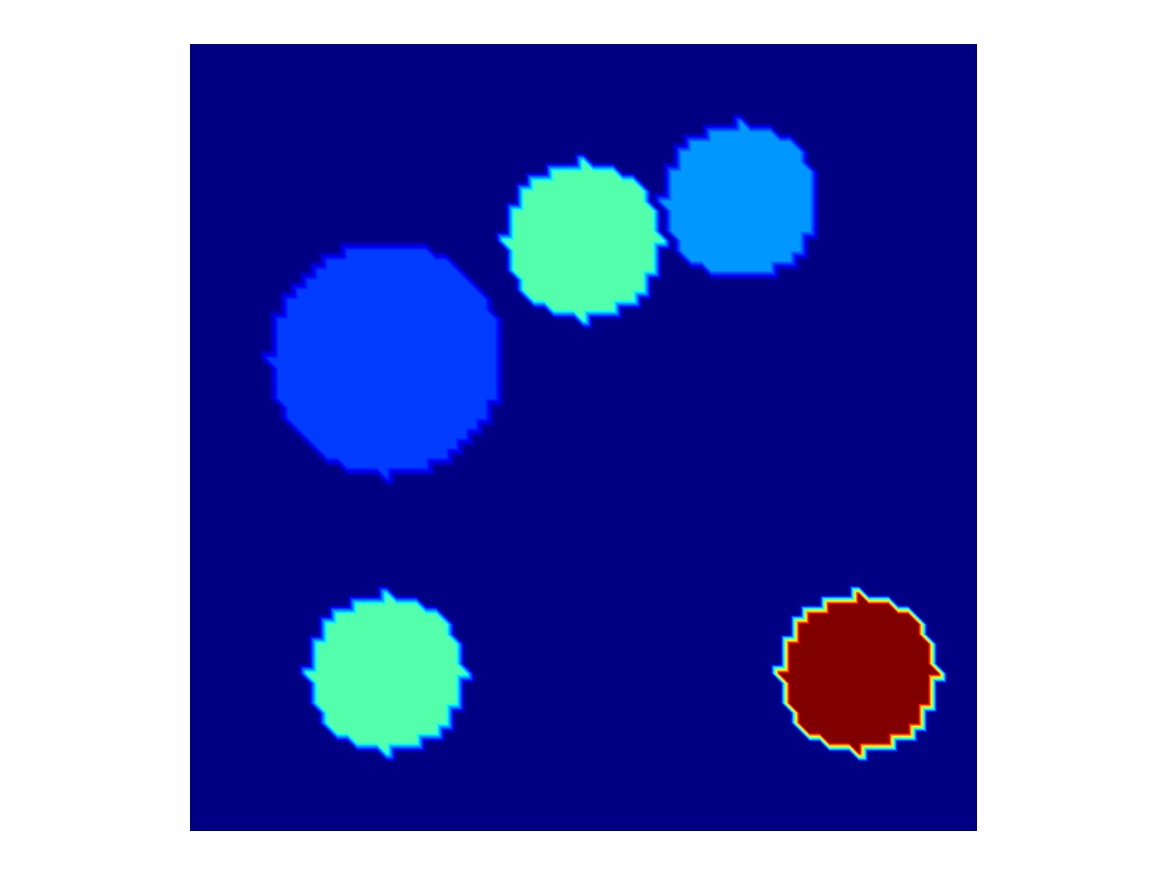}&
\includegraphics[width=.25\linewidth]{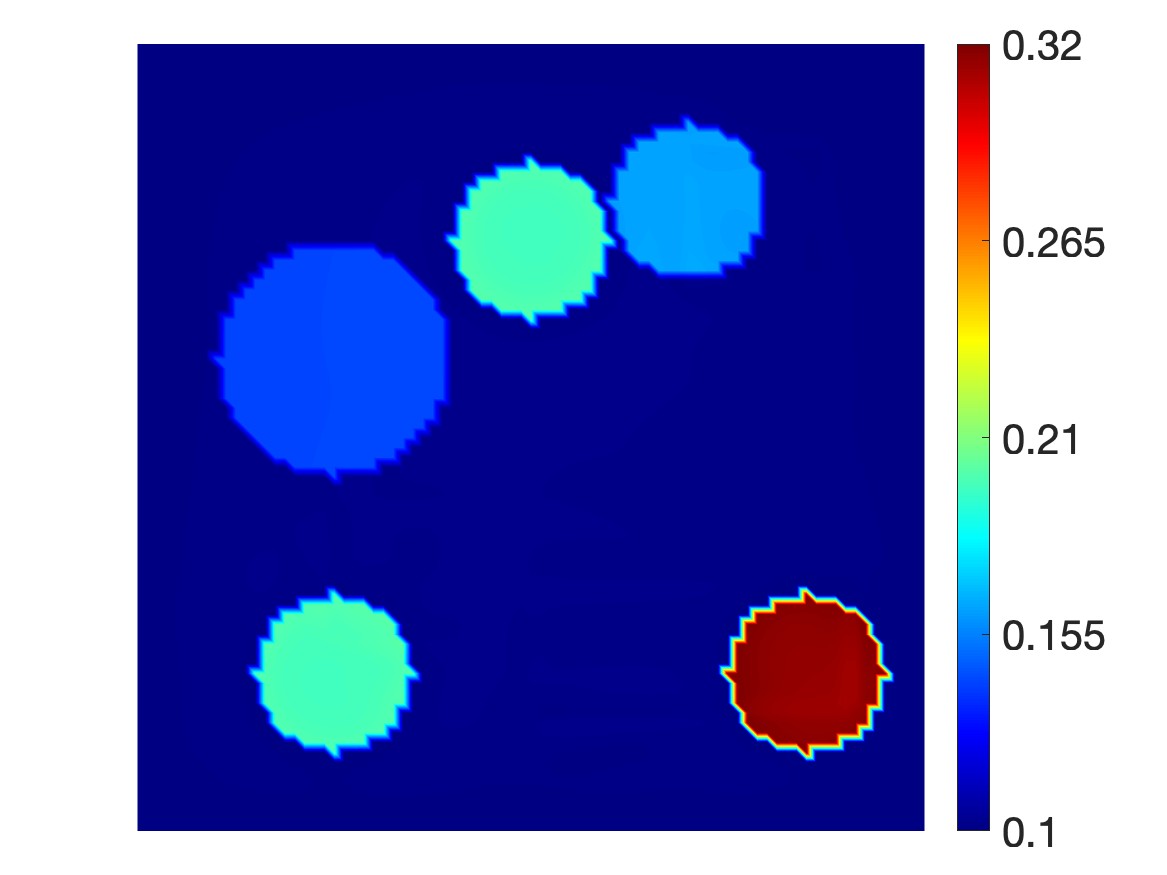}&
\includegraphics[width=.25\linewidth]{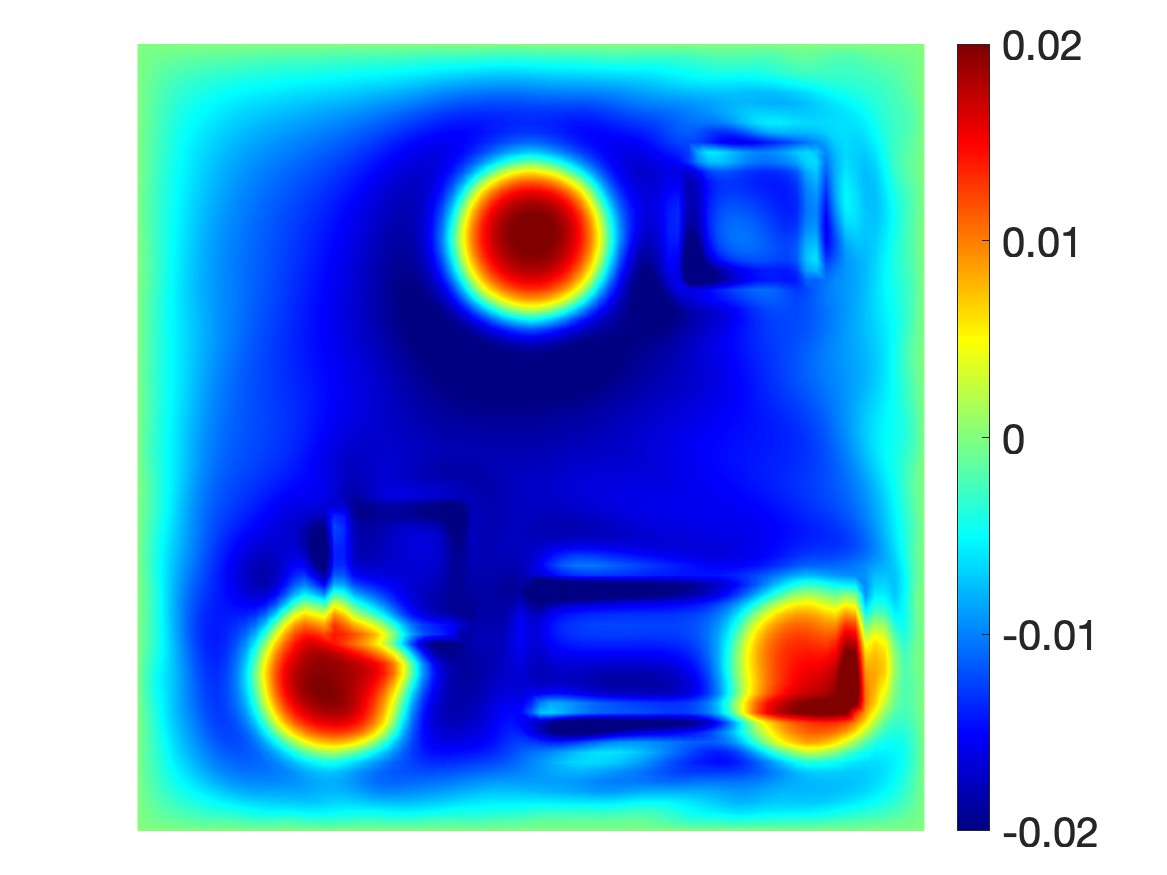}\\
& True $\Gamma\sigma$ & Reconstructed $\Gamma\sigma$ & Relative error
\end{tabular}
\caption{Reconstruction of the product $\Gamma\sigma$ in Experiment III with noise-free (top row) and noisy (bottom row) data.}
\label{FIG:e3-cross}
\end{figure}

Cross-talk between different coefficients is again clearly observed. The pattern of $\Gamma$ appears in the reconstruction of $\sigma$, while patterns of $\sigma$ appear in the reconstruction of $\Gamma$. These are shown in both the reconstruction results and their corresponding relative errors in~\Cref{FIG:e3 noise-free} and~\Cref{FIG:e3 noisy}. To further study the crosstalk between $\Gamma$ and $\sigma$, we show in~\Cref{FIG:e3-cross} the reconstruction of the product $\Gamma\sigma$. Once again, the reconstruction of $\Gamma\sigma$ is relatively better, with approximately $2\%$ relative error at most. 

\RED{Cross-talks between different coefficients have been observed in many multi-coefficient inverse problems. For instance, the cross-talk between the absorption and the diffusion coefficients we observed in Experiment I is well-known in the DOT literature~\cite{Arridge-IP99,ReBaHi-SIAM06}. Cross-talks can happen even when the uniqueness of the reconstruction is mathematically guaranteed. In such a case, they often occur due to the fact that changes in different coefficients can result in similar (but not identical) changes in the data. In the case of the cross-talk between $\sigma$ and $\gamma$, it is easy to see from the Liouville transform we discussed in~\Cref{SEC:Theory} that the potential the data in DOT is determined by $q+i\frac{\omega}{c\gamma}$. Under the same boundary condition, any combination of $\sigma$ and $\gamma$ that gives the same $q+i\frac{\omega}{c\gamma}$ will lead to the same data. It is clear that one can give many combinations of $\sigma$ and $\gamma$ that lead to similar $q+i\frac{\omega}{c\gamma}$ values (which then lead to similar data). The reconstruction algorithm is then confused about the differences between the different combinations when starting from initial guesses that are far from the true coefficients. Of course, if one is lucky enough to find the true (unique) global minimum of the objective functional, this cross-talk will disappear. This could be done by starting the algorithm from an initial guess that is extremely close to the true coefficients. In our simulations, we start the algorithm with initial guesses of spatially homogeneous media, an approach that is more realistic than having super accurate initial guesses. This is why we do not expect the cross-talk to disappear. In the same manner, the cross-talk between $\Gamma$ and $\sigma$ is mainly caused by the difference between the sensitivity of the data to $\mu:=\Gamma \sigma$ and $\sigma$. The internal data is very sensitive to $\mu:=\Gamma\sigma$ as it is directly proportional to the quantity but is less sensitive to $\sigma$ as the solution $U$ to the diffusion equation is not very sensitive to the change of $\sigma$. In Stage II of~\Cref{ALG:Three-Stage}, we eliminate $\mu$ to reconstruct $(\sigma, \gamma)$. This inverse problem is not very stable, leading to poor reconstructions of $\sigma$. This then leads to a poor reconstruction of $\Gamma$ in Stage III of the reconstruction. All that said, developing methods to reduce or eliminate cross-talks between coefficients is one of our plans for future studies. 
}

\paragraph{Experiment IV.} In the last numerical experiment, we compare our reconstruction algorithm with a one-stage coupling scheme based on a straightforward least-square implementation. This coupling of DOT and QPAT is done by combining the data mismatch terms in DOT and QPAT to have an objective functional that depends on all three coefficients:
\begin{equation}\label{EQ:One-Stage Min}
    \Phi(\gamma,\sigma,\Gamma)=\wt \Phi_{\rm PAT}(\gamma,\sigma,\Gamma)+\Phi_{\rm DOT}(\gamma,\sigma,\Gamma)+ \edited{ \cR(\gamma, \sigma,\Gamma)},
\end{equation}
where $\Phi_{\rm DOT}(\gamma,\sigma,\Gamma)$ is the same as in~\eqref{phi-dot}, and the PAT consists of a direct least-square data misfit:
\begin{equation}
\wt \Phi_{\rm PAT}(\gamma, \sigma,\Gamma) = \dfrac{1}{2}\sum_{j=1}^{N_s}\int_{\Omega} (H_j-H_j^*)^2 d\bx\,.
\end{equation}
The regularization functional $\cR(\gamma, \sigma,\Gamma)$ takes the same form as in~\eqref{reg}. The main point is that this approach aims to reconstruct all three coefficients simultaneously while our~\Cref{ALG:Three-Stage} reconstructs $(\sigma, \gamma)$ first and then $\Gamma$.

\begin{figure}[!htb]
\settoheight{\tempdima}{\includegraphics[width=.25\linewidth]{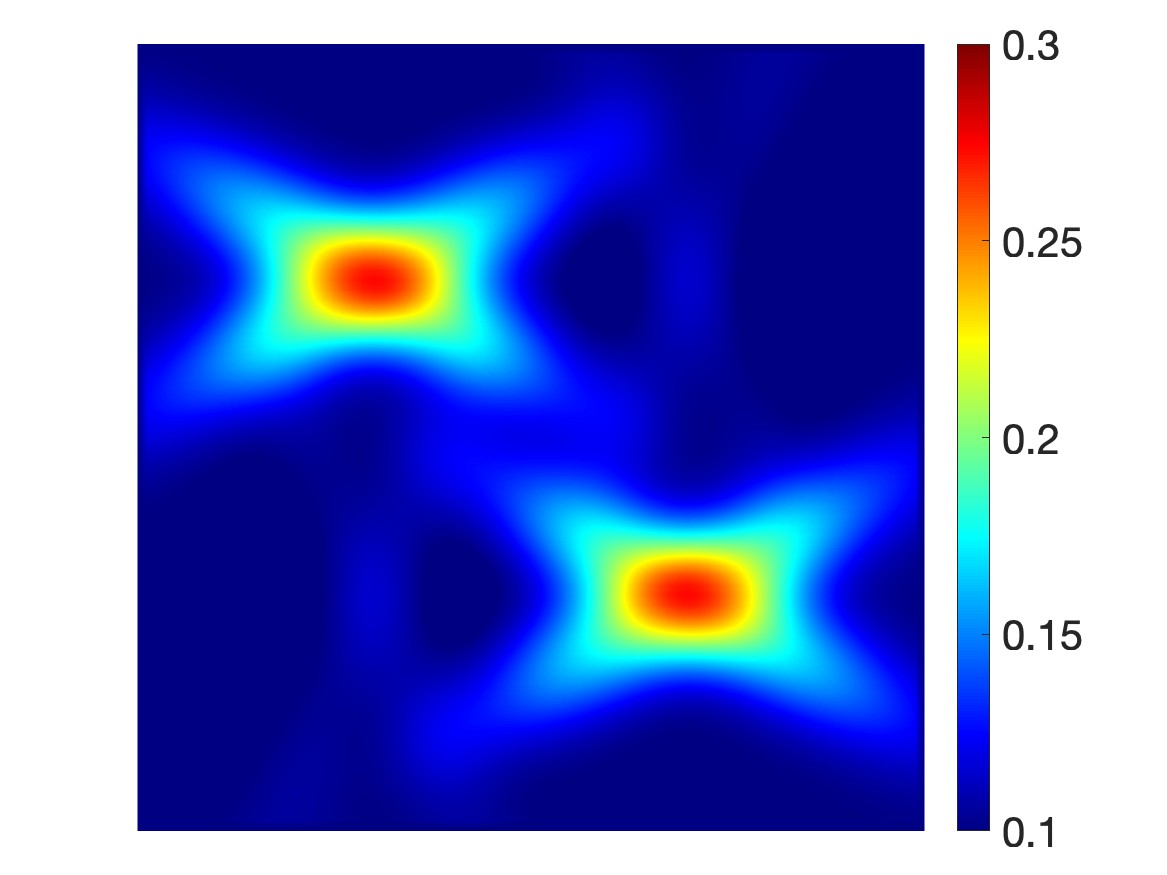}}%
\centering
\begin{tabular}{@{}c@{}c@{}c@{}c@{}}
\rowname{$\sigma$}
&\includegraphics[width=.3\linewidth]{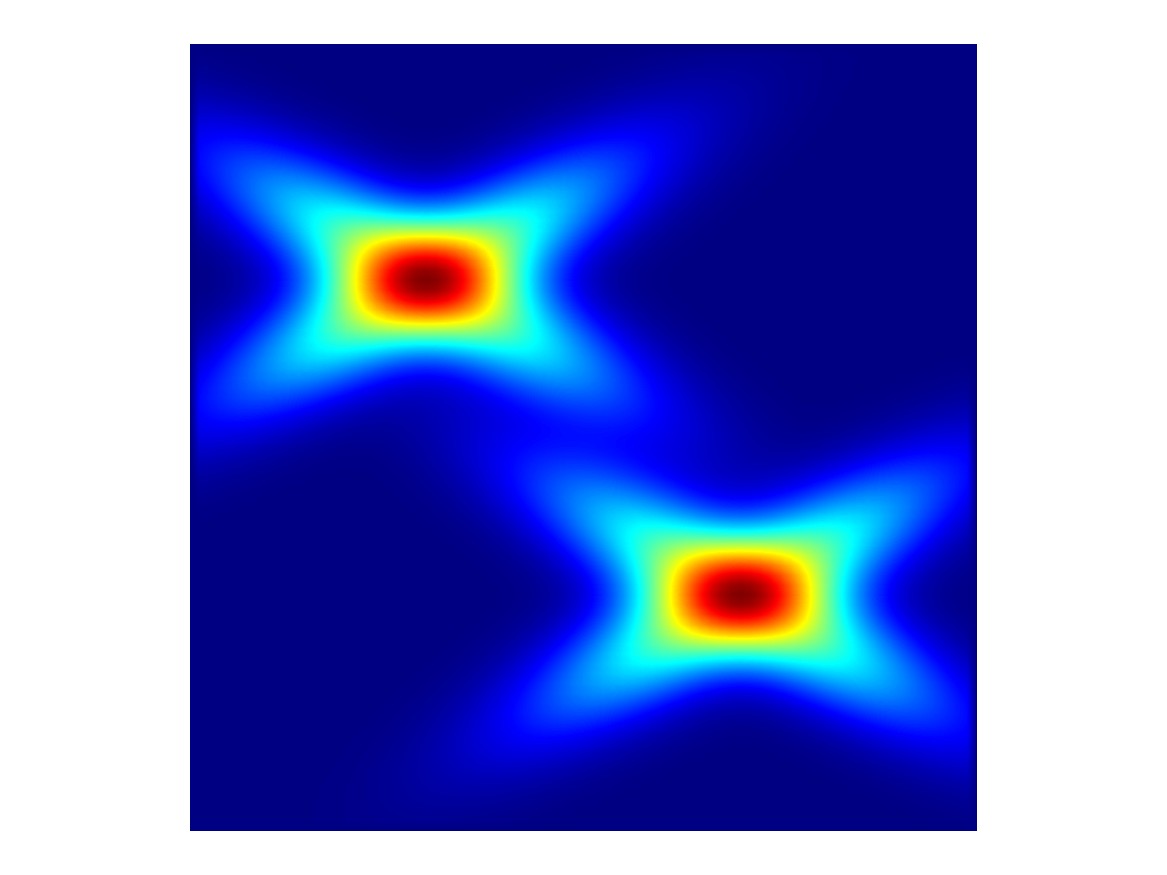}
&\includegraphics[width=.3\linewidth]{FinalFigures/DirectCoupleGaussian/sigmar.jpg}
& \includegraphics[width=.3\linewidth]{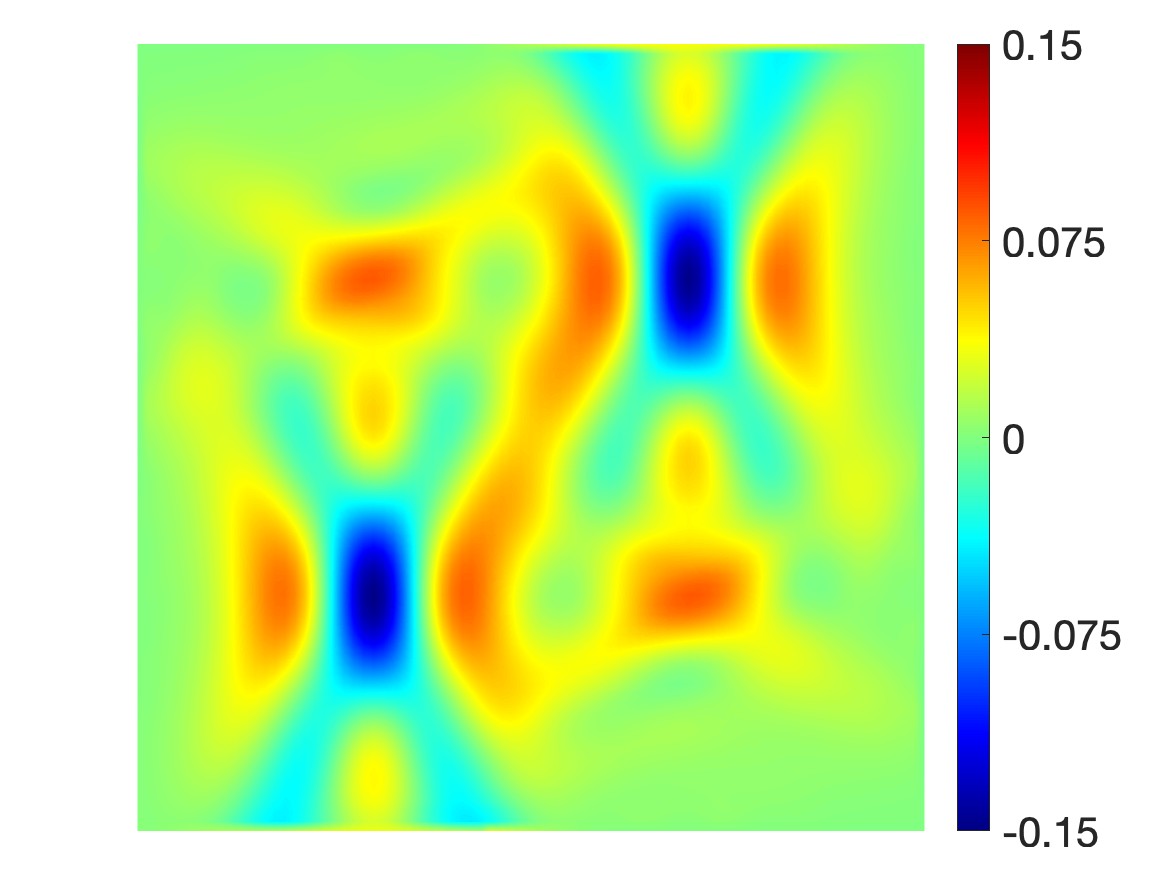}\\
\rowname{$\gamma$}
&\includegraphics[width=.3\linewidth]{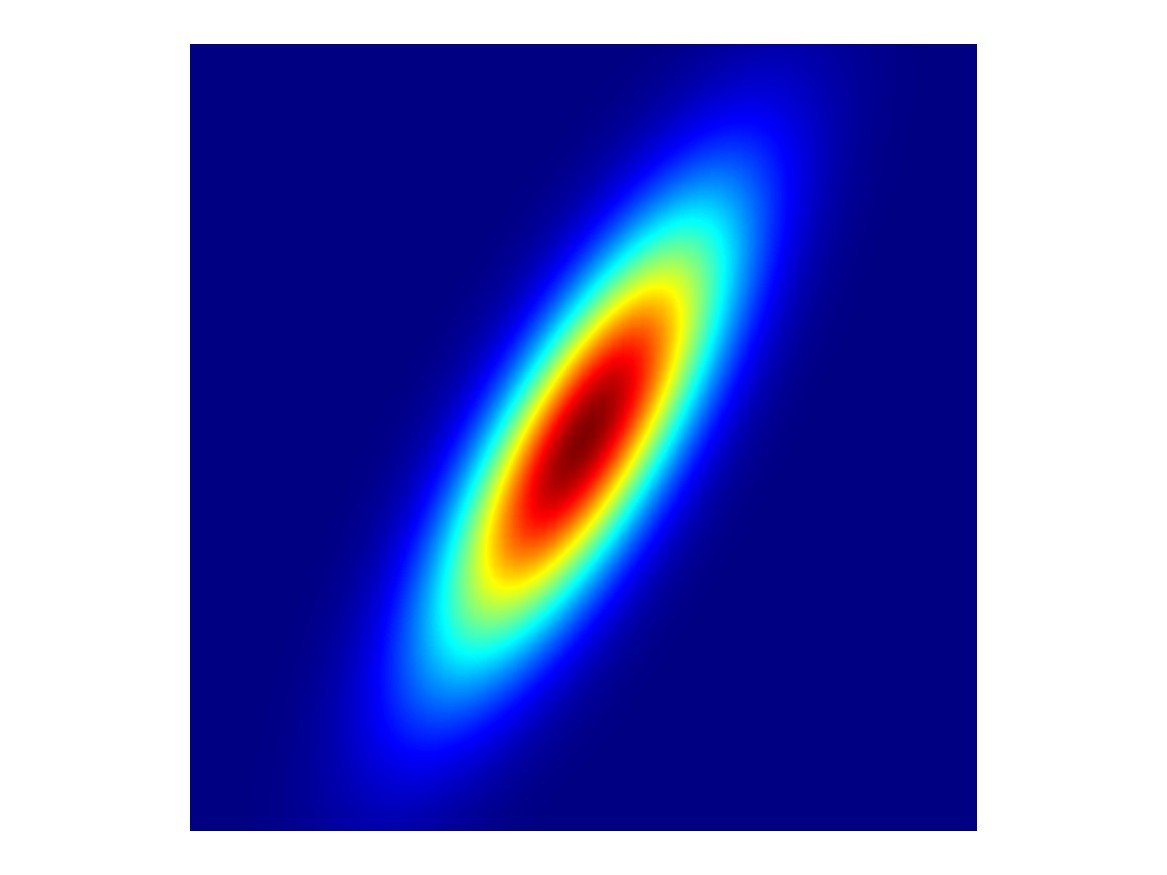}
&\includegraphics[width=.3\linewidth]{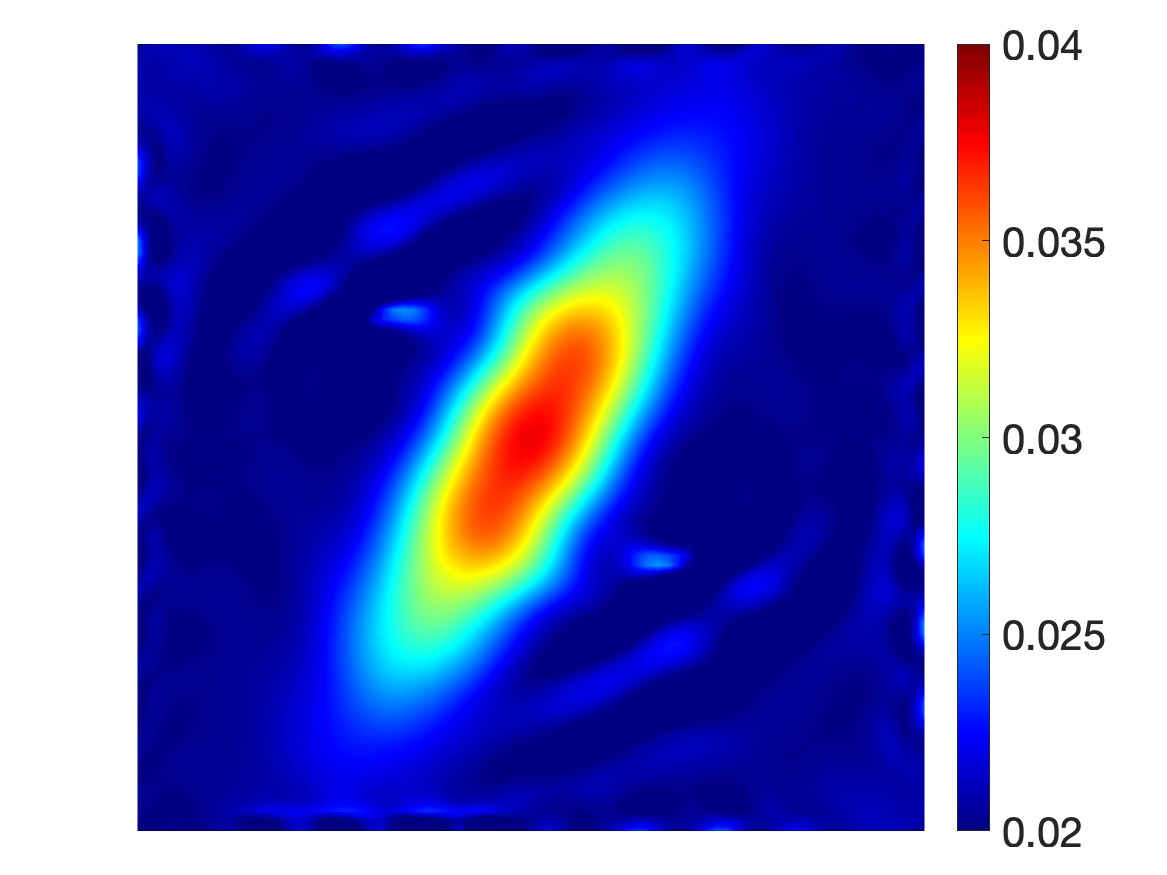}
&\includegraphics[width=.3\linewidth]{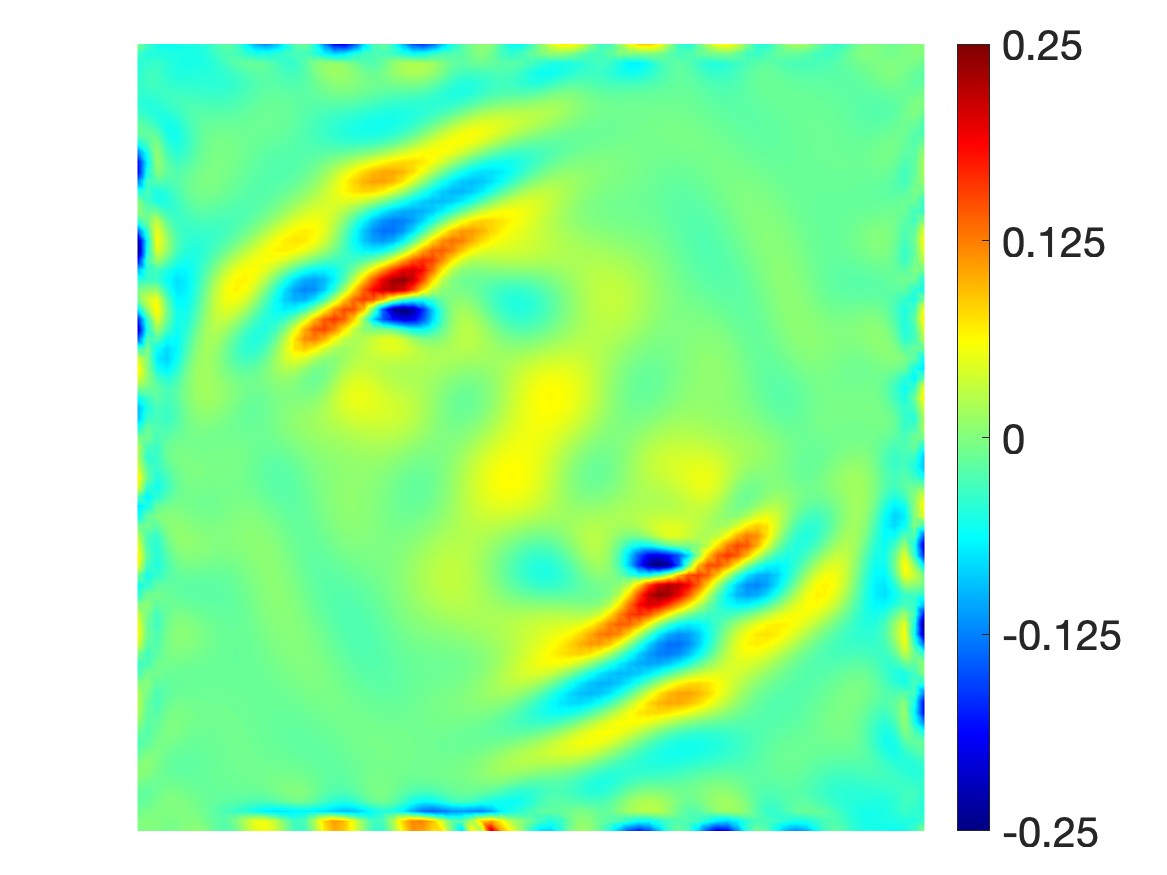}\\
\rowname{$\Gamma$} 
&\includegraphics[width=.3\linewidth]{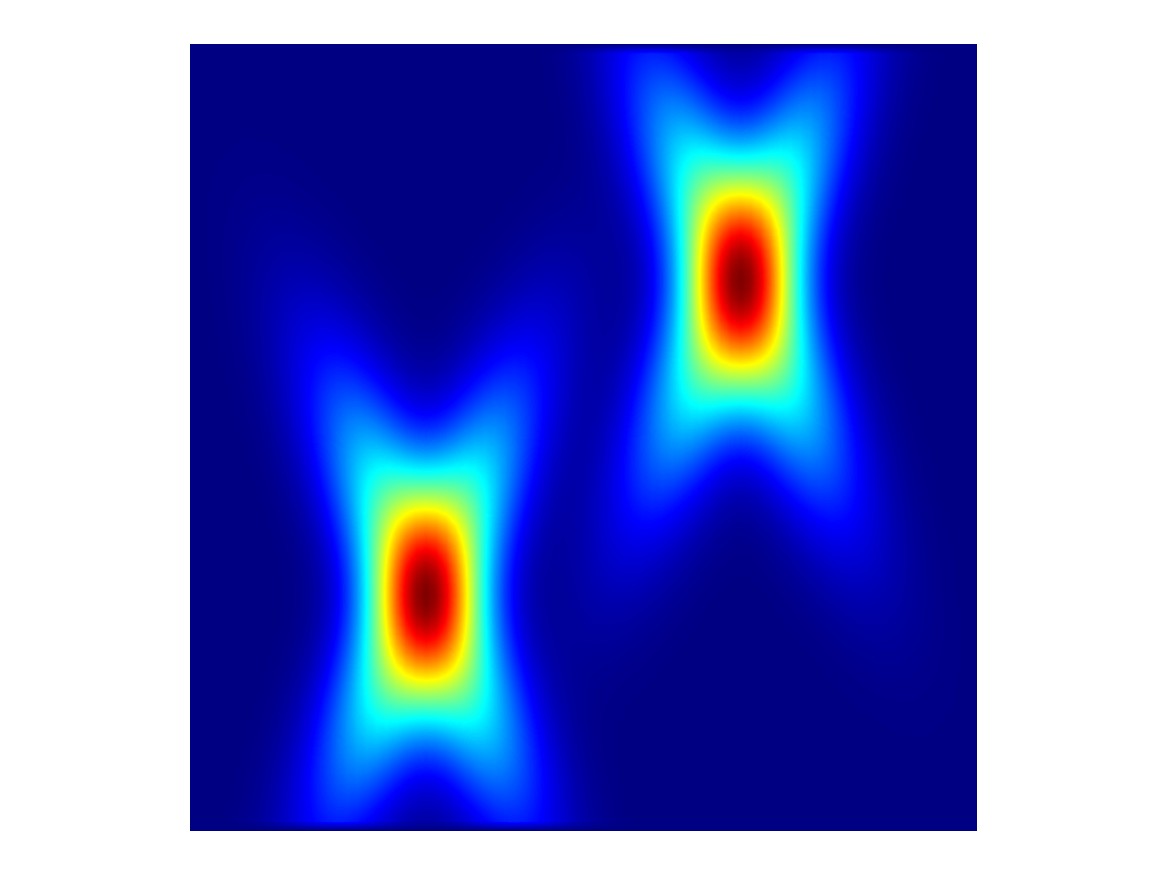}
&\includegraphics[width=.3\linewidth]{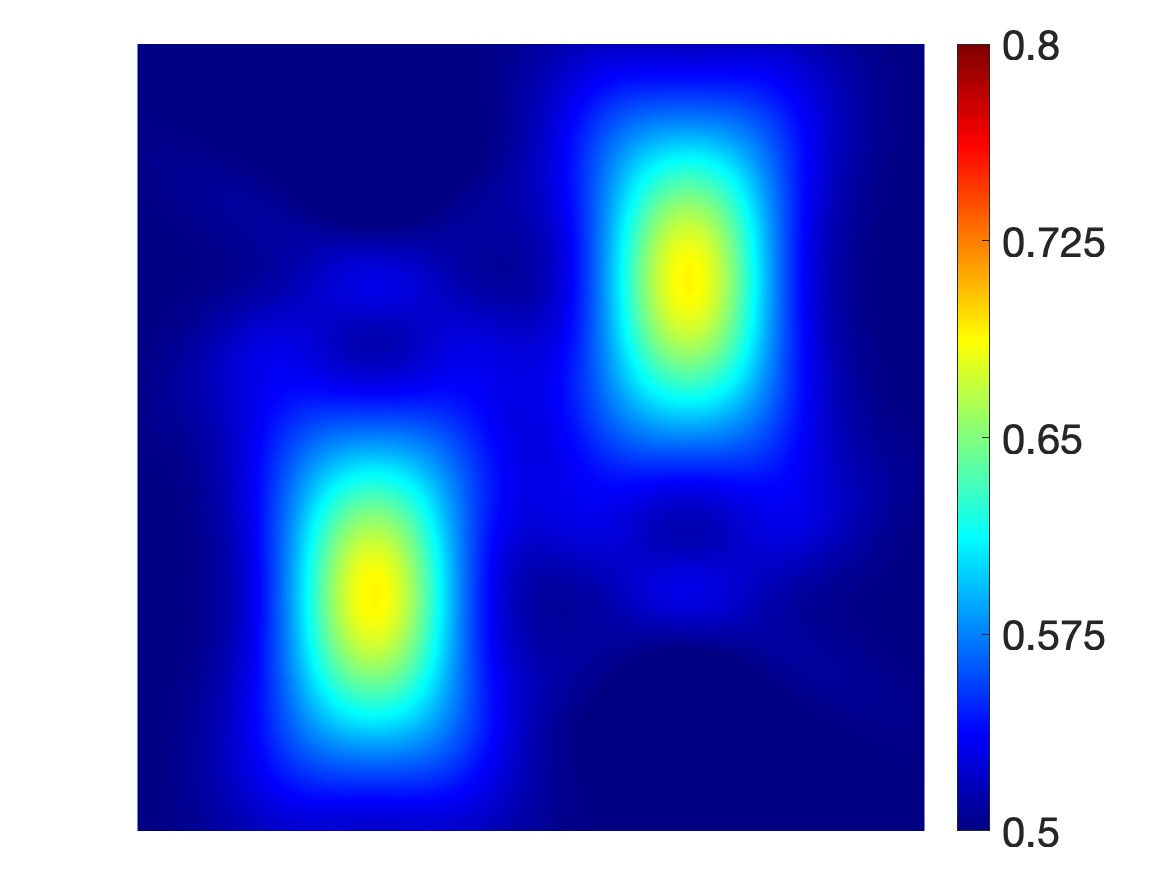}
&\includegraphics[width=.3\linewidth]{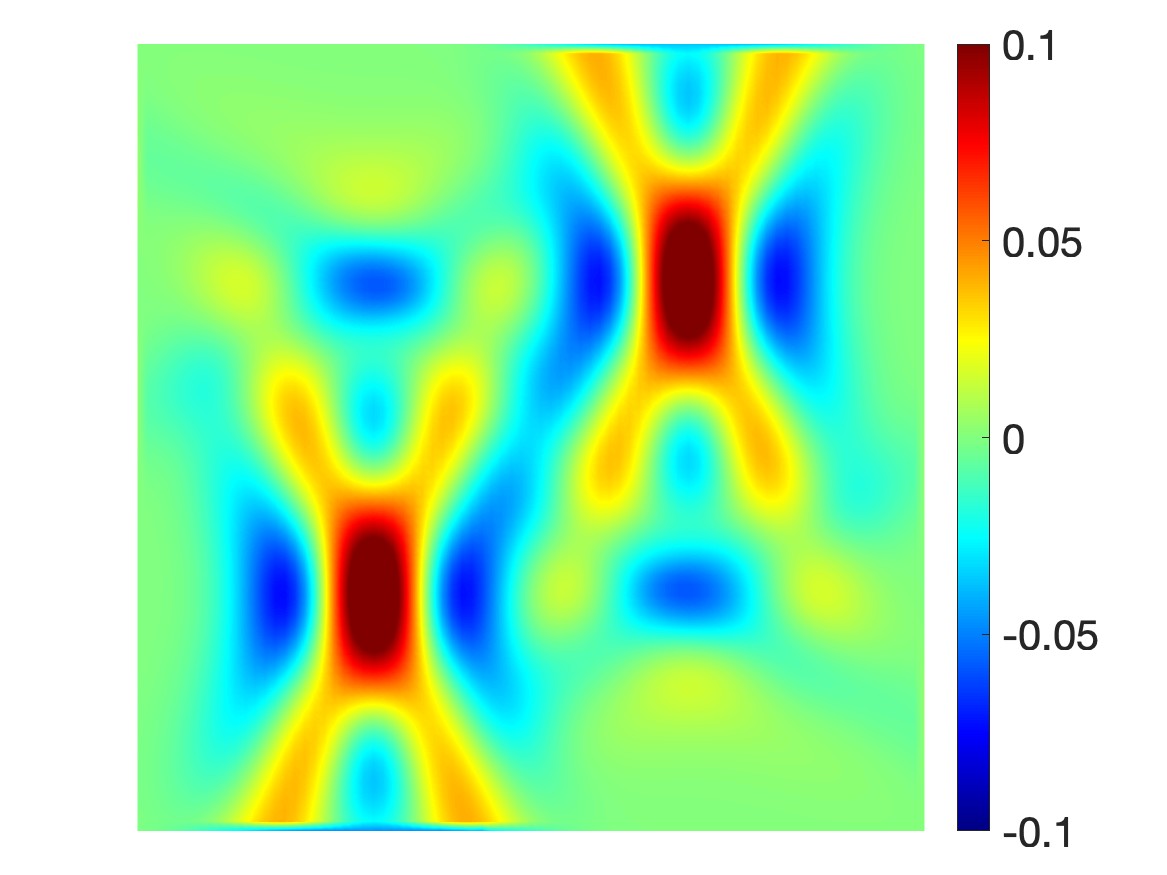}\\
& True coefficients & Reconstructions & Relative error
\end{tabular}
\caption{Reconstructions and relative errors with the single-stage method based on objective functional~\eqref{EQ:One-Stage Min} in Experiment IV.
}
\label{FIG:e4}
\end{figure}

We implemented this single-stage reconstruction algorithm based on the minimization of~\eqref{EQ:One-Stage Min}. In~\Cref{FIG:e4}, we show the reconstruction of the coefficient profiles given in~\eqref{EQ:oval} so that we could directly compare with the results in Experiment II. The reconstructions shown are done with noise-free synthetic data generated by the diffusion model with the Dirichlet boundary condition.

In general, we observe that the single-stage algorithm is less stable and more sensitive to the algorithm's starting point (that is, the initial guess of the coefficients). The results in~\Cref{FIG:e4} can be directly compared with the corresponding results for the Dirichlet boundary conditions in~\Cref{FIG:e2}) that we obtained with our~\Cref{ALG:Three-Stage}. The single-stage algorithm generated results with larger errors (approximately five times in general). This seems especially true for $\Gamma$. The reconstruction results lose entirely the details of the Gaussian mixture shapes. Cross-talk between coefficients, especially between $\sigma$ and $\Gamma$, is even more obvious in this case.

While the result in~\Cref{FIG:e4} is only a particular numerical example, what we observe in the numerical tests we performed is that the single-stage method is easier to be trapped in a local minimal after a few iterations. On the same set of experiments, our reconstruction procedure in~\Cref{ALG:Three-Stage} continuously reduces the value of the objective function until it reaches the stopping criterion we set. With the separation of the stages of the reconstruction, our method simplifies the overall optimization landscape to reduces the difficulty of simultaneous recovery of the three coefficients.  

\section{Conclusions and further remarks}
\label{SEC:Con}

\RED{We proposed a three-stage numerical image reconstruction method for an inverse coefficient problem in quantitative photoacoustic tomography coupled with additional optical current data supplemented from diffuse optical tomography. The objective is to simultaneously reconstruct the absorption, scattering, and Gr\"uneisen coefficients from both QPAT and DOT data. The additional optical data ensures the uniqueness result for all coefficients.
We formulate the inverse problem as a PDE-constrained optimization problem for the absorption and scattering coefficients and reconstruct the Gr\"uneisen coefficient by the averaged data. We also demonstrate that when the Gr\"uneisen coefficient is known,  the addition of the optical measurements allows to have a more balanced reconstruction of the scattering and absorption coefficients. Detailed numerical simulation results based on synthetic data with both smooth and sharp boundaries are presented to demonstrate the performance of the proposed method, compared with direct coupling methods using the least square misfit.}

\RED{The numerical results indeed show some promise in reconstructing three coefficients simultaneously. The key ingredient of the algorithm is the design of the joint objective functional \eqref{joint} to reduce the number of variables involved.} While the convexity of the optimization problem was not specifically addressed, it is generally assumed to be non-convex. The involvement of more variables tends to introduce greater concavity into the objective function. This is further supported by the results depicted in \Cref{FIG:e4}.
The coupled optimization is also believed to be predominantly localized. Our further experimental observations indicate that grid size significantly influences the identification of the minimizer: a coarser grid tends to yield more localized minimizers. Also, the results from Dirichlet boundary conditions are better than Robin because Robin involves one more discretization along the boundary.
Regardless of the run time and computational cost, theoretically, we can still find the global minimizer to the coupled problem without introducing a last step in \Cref{ALG:Three-Stage}. 

While the inclusion of an extra reconstruction step enhances our results, there remains a concern about cross-talk effects stemming from the joint design. As indicated in \eqref{joint}, not only is $\Gamma$ factored out, but $\sigma$ is also separated from the mismatched data. Finding an analytical solution to this coupled problem is challenging without a concrete form of the coefficients, even in a one-dimensional scenario. However, it is plausible that $\sigma$ impacts the forward solution $U$ predominantly by a negative exponential rate. Extracting $\sigma$ from the square also factors out the impact of $\sigma$ on the objective function and, therefore, makes it less sensitive to the gradient. Nonetheless, this should not be a concern as long as we have a relationship between $\Gamma$ and $\sigma$ or they are not close to each other in small medium \cite{CoLaArBe-JBO12,XuWa-RSI06,YaZhMaWa-JBO14}. Exploring alternative designs for the objective function that can simultaneously decrease the number of variables while preserving sensitivity would certainly be an interesting direction to continue. In recovering the coefficients whose shape has sharp boundaries, the errors are concentrated along the boundaries, which might be caused by the Tikhonov-type regularization. A possible remedy is to consider other regularization types that capture the sharp boundaries better, such as the total-variational type, which is an extension that we can further study. 

\section*{Acknowledgments}

This work is partially supported by the National Science Foundation through grants DMS-1937254 and DMS-2309802. \RED{We would like to thank the anonymous referees for their constructive comments that helped us improve the quality of this work.}

\appendix
\section{Further discussions on Stage II of~\texorpdfstring{\Cref{ALG:Three-Stage}}{}}
\label{SEC:Sensitivity}

The key ingredient of our three-stage method is the decomposition of the reconstruction of $(\sigma, \gamma, \Gamma)$ into the reconstructions of $(\sigma, \gamma)$ and then $\Gamma$. This decomposition simplified the complex landscape of the original optimization problem. We now discuss in a little more detail how this decomposition is achieved through the novel objective function in Stage II.

Let us denote by $(\Gamma^*, \gamma^*, \sigma^*)$ the true coefficients that generate the data. We can then write the integral for data mismatch in photoacoustic part, that is, $\Phi_{\rm PAT}$ in~\eqref{EQ:Phi PAT}, as
\begin{equation}\label{joint}
\begin{aligned}
    \Phi_{\rm PAT}(\Gamma, \gamma, \sigma) & = \frac{1}{2}\sum_{i=2}^{N_s}\sum_{j<i} \int_{\Omega} (H_i^*H_j-H_i H_j^*)^2 d\bx \\ & =	\frac{1}{2}\sum_{i=2}^{N_s}\sum_{j<i}
\int_{\Omega} \Gamma^2\Gamma^{*2}\sigma^2\sigma^{*2}(U_i^* U_j- U_i U_j^*)^2 d\bx.
\end{aligned}
\end{equation}
{\bf Independence of the global minimizer on $\Gamma$.} If we assume that the data $\{H_j\}_{j=1}^{N_s}$ are either noise-free or noisy but stay in the range of the forward operator $\Lambda_{\Gamma,\gamma,\sigma}^{\rm QPAT}: g\mapsto H_g$ in~\eqref{EQ:QPAT Data Map}, then~\eqref{joint} shows an obvious feature of the functional $\Phi_{\rm PAT}$: the global minimizer of $\Phi_{\rm PAT}$ is $(\gamma, \sigma)=(\gamma^*, \sigma^*)$, 
which is independent of the value of $\Gamma$ whenever $\Gamma>0$ and $\sigma>0$. This is simply because the mismatch is only determined by the difference between the solutions to the diffusion equations $\{U_j\}_{j=1}^{N_s}$, and such solutions are only dependent on $(\gamma, \sigma)$, not $\Gamma$. This is the main factor that allows us to view $\Phi_{\rm PAT}$ as a function of $(\gamma, \sigma)$ only in Stage II of the reconstruction algorithm.

While the global minimizer of $\Phi_{\rm PAT}$ does not depend on $\Gamma$, the landscape of $\Phi_{\rm PAT}$ as a functional of $\gamma$ and $\sigma$ indeed depends on $\Gamma$, because the factor $\Gamma^2(\bx)$ weights the data mismatch $(U_i^* U_j- U_i U_j^*)^2$ in~\eqref{joint} at different spatial locations differently. Therefore, in Stage II, one could intentionally choose a $\Gamma$ that weights the data mismatch in a desired manner if needed. We do not pursue this direction in our work but just fix $\Gamma$ to be the initial choice $\Gamma_0$ in Stage II and only optimize over $\sigma$ and $\gamma$.

\begin{figure}[!htb]
\settoheight{\tempdima}{\includegraphics[width=.45\linewidth]{FinalFigures/DirectCoupleGaussian/sigmar.jpg}}%
\centering
\begin{tabular}{@{}c@{}c@{}}
\includegraphics[height=.35\linewidth]{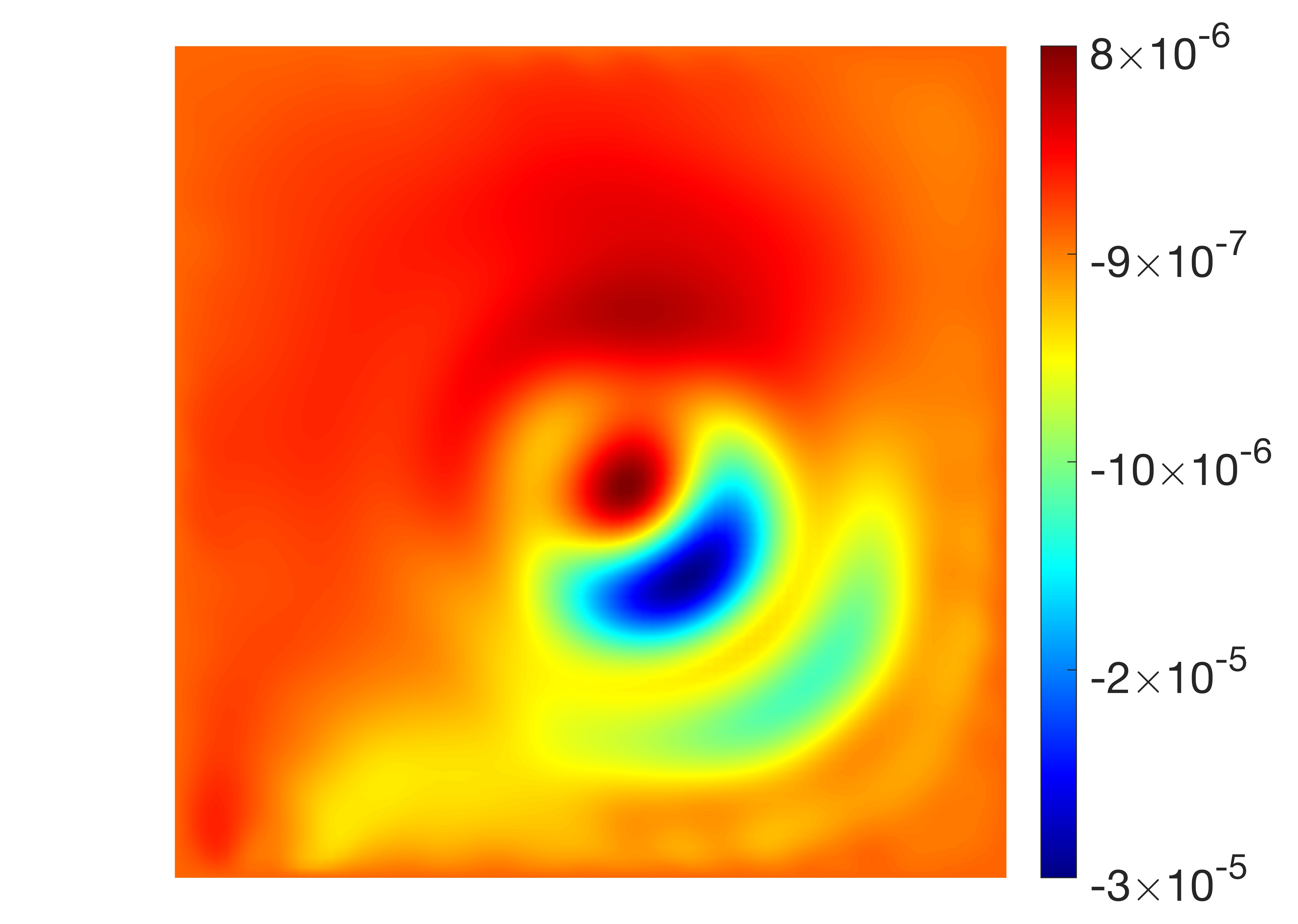}&
\includegraphics[height=.35\linewidth]{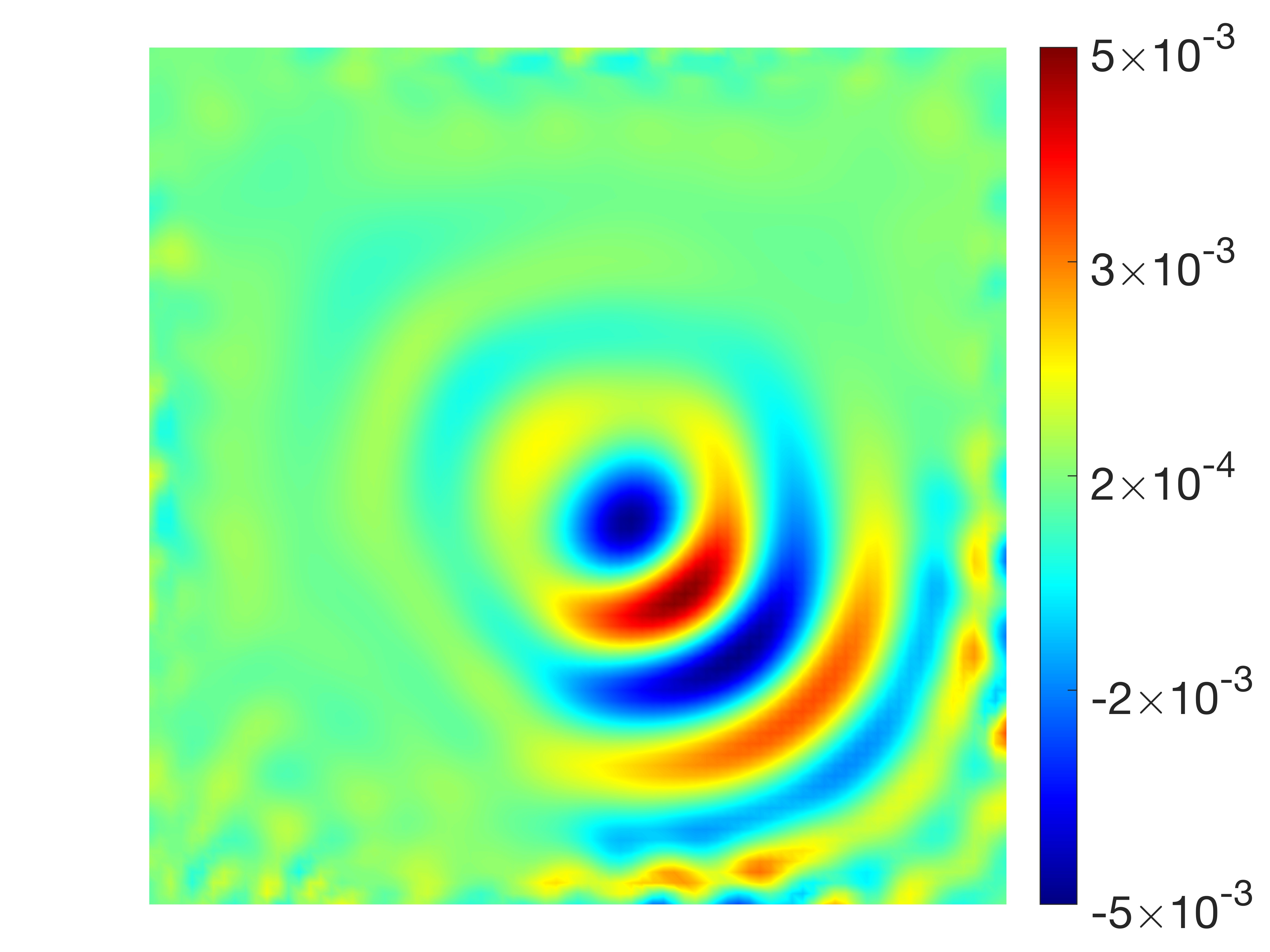}\\
$\frac{\sigma_1-\sigma_2}{\sigma_1}$ &
$\frac{\gamma_1-\gamma_2}{\gamma_1}$
\end{tabular}
\caption{The relative differences between the reconstructions of $\sigma$ (left) and $\gamma$ (right) in Stage II using different $\Gamma$.
}
\label{FIG:perturbation}
\end{figure}
To verify the above reasoning of insensitivity of the objective function to $\Gamma$ in Stage II, we show in~\Cref{FIG:perturbation} the differences between the reconstruction of the coefficients $(\sigma,\gamma)$ in Stage II with $\Gamma$ taking two different values $\Gamma_1(\bx) = 0.3\exp\big(-10|\bx-\bx_0|^2\big),\bx_0=(1.4, 0.6)$,
or 
$\Gamma_2(\bx) = 16$. The true coefficients are
\begin{equation}%
\begin{aligned}
\sigma^*(\bx) &=&& 0.3\exp \big(-10|(\bx-\bx_0)|^2\big), & \bx_1=(0.4, 1.5), \\
\gamma^*(\bx) &=&& 0.04\exp\big(-10|(\bx-\bx_1)|^2\big),&\bx_2=(1.0, 1.0), \\
\end{aligned}
\end{equation}

The two reconstructions are $(\sigma_1, \gamma_1)$ with $\Gamma$ fixed to be $\Gamma_1$ and $(\sigma_2, \gamma_2)$ with $\Gamma_2$. It is clear from the figures of the relative difference, $\dfrac{\sigma_1-\sigma_2}{\sigma_1}$ and $\dfrac{\gamma_1-\gamma_2}{\gamma_1}$, between the reconstructions that the value of $\Gamma$ has little influence on the reconstruction result. This is the numerical evidence of why we could eliminate $\Gamma$ in Stage II of~\Cref{ALG:Three-Stage}.

\paragraph{Elimination of $\Gamma$ via Normalization.} In fact, it is possible to completely eliminate the dependence of the objective function on $\Gamma$ in Stage II. This can be achieved by forming the ratio between data generated by different sources. More precisely, the following objective function for the photoacoustic part works:
\begin{multline}\label{EQ:Obj Ratio}
    \wh \Phi_{\rm PAT}(\Gamma, \gamma, \sigma) :=\frac{1}{2}\sum_{i=2}^{N_s}\sum_{j<i} \int_{\Omega} \Big(\frac{H_i}{H_j}-\frac{H_i^*}{H_j^*}\Big)^2 d\bx =\frac{1}{2}\sum_{i=2}^{N_s}\sum_{j<i} \int_{\Omega} \Big(\frac{H_iH_j^*-H_i^*H_j}{H_jH_j^*}\Big)^2 d\bx\,.
\end{multline}
Since the ratio $\frac{H_i}{H_j}=\frac{\Gamma\sigma U_i}{\Gamma\sigma U_j}=\frac{U_i}{U_j}$ is independent of $\Gamma$, the functional $\wh \Phi_{\rm PAT}$ is independent of $\Gamma$. The functional depends on $\gamma$ and $\sigma$ since the solutions $\{U_j\}_{j=1}^{N_s}$ depend on $\gamma$ and $\sigma$ through the diffusion equation. The second equality in~\eqref{EQ:Obj Ratio} shows that this new objective function is simply a normalized version of $\Phi_{\rm PAT}$ in~\eqref{EQ:Phi PAT}. However, this normalization is nonlinear as the model prediction $H_j$, not only the measured data $H_j^*$, is on the denominator. This makes the calculation of the gradient of the objective function a little more complicated. Using the normalized objective function requires the quantities on the denominator to be strictly positive. This limits the illuminating source we could use in practice. These are the main reasons we think~\eqref{EQ:Phi PAT} is more practically relevant.

{\small
\bibliography{BIB-REN}

\begin{thebibliography}{10}

\bibitem{Alessandrini-AA88}
{\sc G.~Alessandrini}, {\em Stable determination of conductivity by boundary
  measurements}, Applic. Analysis, 27 (1988), pp.~153--172.

\bibitem{AlVe-AAM05}
{\sc G.~Alessandrini and S.~Vessella}, {\em Lipschitz stability for the inverse
  conductivity problem}, Advances in Applied Mathematics, 35 (2005),
  pp.~207--241.

\bibitem{Arridge-IP99}
{\sc S.~R. Arridge}, {\em Optical tomography in medical imaging}, Inverse
  Probl., 15 (1999), pp.~R41--R93.

\bibitem{ArLi-OL98}
{\sc S.~R. Arridge and W.~R.~B. Lionheart}, {\em Nonuniqueness in
  diffusion-based optical tomography}, Opt. Lett., 23 (1998), pp.~882--884.

\bibitem{ArSc-IP09}
{\sc S.~R. Arridge and J.~C. Schotland}, {\em Optical tomography: forward and
  inverse problems}, Inverse Problems, 25 (2009).
\newblock 123010.

\bibitem{AsBeCaCaNa-DCDS24}
{\sc A.~Aspri, A.~Benfenati, P.~Causin, C.~Cavaterra, and G.~Naldi}, {\em
  Mathematical and numerical challenges in diffuse optical tomography inverse
  problems}, Discrete and Continuous Dynamical Systems-S, 17 (2024),
  pp.~421--461.

\bibitem{Bal-IP09}
{\sc G.~Bal}, {\em Inverse transport theory and applications}, Inverse
  Problems, 25 (2009).
\newblock 053001.

\bibitem{BaRe-IP11}
{\sc G.~Bal and K.~Ren}, {\em Multi-source quantitative {PAT} in diffusive
  regime}, Inverse Problems, 27 (2011).
\newblock 075003.

\bibitem{BaRe-IP12}
\leavevmode\vrule height 2pt depth -1.6pt width 23pt, {\em On multi-spectral
  quantitative photoacoustic tomography in diffusive regime}, Inverse Problems,
  28 (2012).
\newblock 025010.

\bibitem{BaUh-IP10}
{\sc G.~Bal and G.~Uhlmann}, {\em Inverse diffusion theory of photoacoustics},
  Inverse Problems, 26 (2010).
\newblock 085010.

\bibitem{Beard-IF11}
{\sc P.~Beard}, {\em Biomedical photoacoustic imaging}, Interface Focus, 1
  (2011), pp.~602--631.

\bibitem{BeCaQu-arXiv24}
{\sc A.~Benfenati, P.~Causin, and M.~Quinteri}, {\em A modular deep
  learning-based approach for diffuse optical tomography reconstruction},
  arXiv:2402.09277,  (2024).

\bibitem{BeFrVe-SIAM21}
{\sc E.~Beretta, E.~Francini, and S.~Vessella}, {\em Lipschitz stable
  determination of polygonal conductivity inclusions in a two-dimensional
  layered medium from the dirichlet-to-neumann map}, SIAM Journal on
  Mathematical Analysis, 53 (2021), pp.~4303--4327.

\bibitem{CoLaArBe-JBO12}
{\sc B.~Cox, J.~G. Laufer, S.~R. Arridge, and C.~Beard}, {\em Quantitative
  spectroscopic photoacoustic imaging: a review}, J. Biomed. Optics, 17 (2012).
\newblock 061202.

\bibitem{CoLaBe-SPIE09}
{\sc B.~T. Cox, J.~G. Laufer, and P.~C. Beard}, {\em The challenges for
  quantitative photoacoustic imaging}, Proc. of SPIE, 7177 (2009).
\newblock 717713.

\bibitem{ElMiSc-M2AS17}
{\sc P.~Elbau, L.~Mindrinos, and O.~Scherzer}, {\em Inverse problems of
  combined photoacoustic and optical coherence tomography}, Mathematical
  Methods in the Applied Sciences, 40 (2017), pp.~505--522.

\bibitem{Foschiatti-JMAA24}
{\sc S.~Foschiatti}, {\em Lipschitz stability estimate for the simultaneous
  recovery of two coefficients in the anisotropic schr{\"o}dinger type equation
  via local cauchy data}, Journal of Mathematical Analysis and Applications,
  531 (2024), p.~127753.

\bibitem{GaOsZh-LNM12}
{\sc H.~Gao, S.~Osher, and H.~Zhao}, {\em Quantitative photoacoustic
  tomography}, in Mathematical Modeling in Biomedical Imaging II: Optical,
  Ultrasound, and Opto-Acoustic Tomographies, H.~Ammari, ed., vol.~2035 of
  Lecture Notes in Mathematics, Springer, 2012, pp.~131--158.

\bibitem{GaScWi-JMIV15}
{\sc T.~Glatz, O.~Scherzer, and T.~Widlak}, {\em Texture generation for
  photoacoustic elastography}, Journal of Mathematical Imaging and Vision, 52
  (2015), pp.~369--384.

\bibitem{HaNeRa-IP15}
{\sc M.~Haltmeier, L.~Neumann, and S.~Rabanser}, {\em Single-stage
  reconstruction algorithm for quantitative photoacoustic tomography}, Inverse
  Problems, 31 (2015).
\newblock 065005.

\bibitem{HaPuArTa-IPI24}
{\sc N.~H{\"a}nninen, A.~Pulkkinen, S.~Arridge, and T.~Tarvainen}, {\em
  Estimating absorption and scattering in quantitative photoacoustic tomography
  with an adaptive monte carlo method for light transport}, Inverse Problems
  and Imaging,  (2024), pp.~0--0.

\bibitem{HiPiUlUl-Book08}
{\sc M.~Hinze, R.~Pinnau, M.~Ulbrich, and S.~Ulbrich}, {\em Optimization with
  PDE constraints}, vol.~23, Springer Science \& Business Media, 2008.

\bibitem{ImYa-IP98}
{\sc O.~Y. Imanuvilov and M.~Yamamoto}, {\em Lipschitz stability in inverse
  parabolic problems by the {Carleman} estimate}, Inverse Problems, 14 (1998),
  pp.~1229--1245.

\bibitem{Isakov-Book06}
{\sc V.~Isakov}, {\em Inverse {Problems} for {Partial} {Differential}
  {Equations}}, Springer-Verlag, New York, second~ed., 2006.

\bibitem{KaMoPuTa-BOE24}
{\sc J.~Kangasniemi, M.~Mozumder, A.~Pulkkinen, and T.~Tarvainen}, {\em
  Stochastic {Gauss-Newton} method for estimating absorption and scattering in
  optical tomography with the monte carlo method for light transport},
  Biomedical Optics Express, 15 (2024), pp.~4925--4942.

\bibitem{KlHi-IP03}
{\sc A.~D. Klose and A.~H. Hielscher}, {\em {Quasi-Newton} methods in optical
  tomographic image reconstruction}, Inverse Probl., 19 (2003), pp.~387--409.

\bibitem{Kuchment-MLLE12}
{\sc P.~Kuchment}, {\em Mathematics of hybrid imaging. a brief review}, in The
  Mathematical Legacy of Leon Ehrenpreis, I.~Sabadini and D.~Struppa, eds.,
  Springer-Verlag, 2012.

\bibitem{KuKu-HMMI10}
{\sc P.~Kuchment and L.~Kunyansky}, {\em Mathematics of thermoacoustic and
  photoacoustic tomography}, in Handbook of Mathematical Methods in Imaging,
  O.~Scherzer, ed., Springer-Verlag, 2010, pp.~817--866.

\bibitem{LiWa-PMB09}
{\sc C.~Li and L.~Wang}, {\em Photoacoustic tomography and sensing in
  biomedicine}, Phys. Med. Biol., 54 (2009), pp.~R59--R97.

\bibitem{MaRe-CMS14}
{\sc A.~V. Mamonov and K.~Ren}, {\em Quantitative photoacoustic imaging in
  radiative transport regime}, Comm. Math. Sci., 12 (2014), pp.~201--234.

\bibitem{Mandache-IP01}
{\sc N.~Mandache}, {\em Exponential instability in an inverse problem for the
  {S}chr\"{o}dinger equation}, Inverse Probl., 17 (2001), pp.~1435--1444.

\bibitem{Nachman-AM96}
{\sc A.~Nachman}, {\em Global uniqueness for a two-dimensional inverse boundary
  value problem}, Ann. Math., 143 (1996), pp.~71--96.

\bibitem{NyPuTa-BOE17}
{\sc O.~Nyk\"{a}nen, A.~Pulkkinen, and T.~Tarvainen}, {\em Quantitative
  photoacoustic tomography augmented with surface light measurements}, Biomed.
  Optics. Exp., 8 (2017), pp.~4382--4395.

\bibitem{PaSc-IP07}
{\sc S.~K. Patch and O.~Scherzer}, {\em Photo- and thermo- acoustic imaging},
  Inverse Problems, 23 (2007), pp.~S1--S10.

\bibitem{RaFa-PRR23}
{\sc J.~Radford and D.~Faccio}, {\em Information transport and limits of
  optical imaging in the highly diffusive regime}, Physical Review Research, 5
  (2023), p.~L022008.

\bibitem{Ren-CiCP10}
{\sc K.~Ren}, {\em Recent developments in numerical techniques for
  transport-based medical imaging methods}, Commun. Comput. Phys., 8 (2010),
  pp.~1--50.

\bibitem{ReBaHi-SIAM06}
{\sc K.~Ren, G.~Bal, and A.~H. Hielscher}, {\em Frequency domain optical
  tomography based on the equation of radiative transfer}, SIAM J. Sci.
  Comput., 28 (2006), pp.~1463--1489.

\bibitem{ReGaZh-SIAM13}
{\sc K.~Ren, H.~Gao, and H.~Zhao}, {\em A hybrid reconstruction method for
  quantitative photoacoustic imaging}, SIAM J. Imag. Sci., 6 (2013),
  pp.~32--55.

\bibitem{Rondi-AAM06}
{\sc L.~Rondi}, {\em A remark on a paper by {Alessandrini and Vessella}},
  Advances in Applied Mathematics, 36 (2006), pp.~67--69.

\bibitem{SaTaCoAr-IP13}
{\sc T.~Saratoon, T.~Tarvainen, B.~T. Cox, and S.~R. Arridge}, {\em A
  gradient-based method for quantitative photoacoustic tomography using the
  radiative transfer equation}, Inverse Problems, 29 (2013).
\newblock 075006.

\bibitem{Scherzer-Book10}
{\sc O.~Scherzer}, {\em Handbook of Mathematical Methods in Imaging},
  Springer-Verlag, 2010.

\bibitem{SiJaPr-EPJP24}
{\sc H.~G. Siddalingaiah, R.~P.~K. Jagannath, and G.~R. Prashanth}, {\em
  Randomized recursive techniques for image reconstruction in diffuse optical
  tomography}, The European Physical Journal Plus, 139 (2024), p.~584.

\bibitem{SiTh-JBO19}
{\sc M.~S. Singh and A.~Thomas}, {\em Photoacoustic elastography imaging: a
  review}, J. Biomed. Opt., 24 (2019), pp.~040902--040902.

\bibitem{SuHoSuQi-PMB23}
{\sc Z.~Sun, Y.~Hou, M.~Sun, and Q.~Meng}, {\em Quantitative photoacoustic
  tomography with light fluence compensation based on radiance monte carlo
  model}, Physics in Medicine \& Biology, 68 (2023), p.~065009.

\bibitem{SyUh-AM87}
{\sc J.~Sylvester and G.~Uhlmann}, {\em Global uniqueness theorem for an
  inverse boundary value problem}, Ann. Math., 125 (1987), pp.~153--169.

\bibitem{TaCo-JBO24}
{\sc T.~Tarvainen and B.~Cox}, {\em Quantitative photoacoustic tomography:
  modeling and inverse problems}, Journal of Biomedical Optics, 29 (2024),
  pp.~S11509--S11509.

\bibitem{Uhlmann-BMS14}
{\sc G.~Uhlmann}, {\em Inverse problems: seeing the unseen}, Bulletin of
  Mathematical Sciences, 4 (2014), pp.~209--279.

\bibitem{Wang-DM04}
{\sc L.~V. Wang}, {\em Ultrasound-mediated biophotonic imaging: a review of
  acousto-optical tomography and photo-acoustic tomography}, Disease Markers,
  19 (2004), pp.~123--138.

\bibitem{Wang-IEEE08}
\leavevmode\vrule height 2pt depth -1.6pt width 23pt, {\em Tutorial on
  photoacoustic microscopy and computed tomography}, IEEE J. Sel. Topics
  Quantum Electron., 14 (2008), pp.~171--179.

\bibitem{XuWa-RSI06}
{\sc M.~Xu and L.~V. Wang}, {\em Photoacoustic imaging in biomedicine}, Rev.
  Sci. Instr., 77 (2006).
\newblock 041101.

\bibitem{YaChWaSh-IEEE24}
{\sc F.~Yang, Z.~Chen, P.~Wang, and Y.~Shi}, {\em Phase-domain photoacoustic
  mechanical imaging for quantitative elastography and viscography}, IEEE
  Trans. Biomed. Eng., 71 (2024).

\bibitem{YaZhMaWa-JBO14}
{\sc D.-K. Yao, C.~Zhang, K.~Maslov, and L.~V. Wang}, {\em Photoacoustic
  measurement of the gr{\"u}neisen parameter of tissue}, Journal of biomedical
  optics, 19 (2014), pp.~017007--017007.

\bibitem{YiYaWaWaGuCaZhHe-BOE24}
{\sc H.~Yi, R.~Yang, Y.~Wang, Y.~Wang, H.~Guo, X.~Cao, S.~Zhu, and X.~He}, {\em
  Enhanced model iteration algorithm with graph neural network for diffuse
  optical tomography}, Biomedical Optics Express, 15 (2024), pp.~1910--1925.

\bibitem{Zemp-AO10}
{\sc R.~J. Zemp}, {\em Quantitative photoacoustic tomography with multiple
  optical sources}, Applied Optics, 49 (2010), pp.~3566--3572.

\bibitem{Zimmermann-IP23}
{\sc P.~Zimmermann}, {\em Inverse problem for a nonlocal diffuse optical
  tomography equation}, Inverse Problems, 39 (2023), p.~094001.

\end{thebibliography}
\bibliographystyle{siam}
}

\end{document}